\documentclass[11pt]{article}
\usepackage{newclude}
\usepackage{bbm}
\usepackage{fullpage}
\usepackage{amsfonts}
\usepackage{bbding}
\usepackage{graphicx,float}
\usepackage{amsfonts,amssymb,amsmath,amsthm,amscd}
\usepackage{mathrsfs}
\usepackage{xcolor}
\usepackage{empheq}
\usepackage{adforn}
\usepackage{cancel}
\usepackage{xr-hyper}
\usepackage{xr}
\usepackage{mdframed}
\usepackage{mathdots}
\usepackage{enumerate}
\usepackage{lmodern}
\usepackage{mdframed}
\usepackage{stmaryrd}
\usepackage{blindtext}
\usepackage[all, knot]{xy}
\usepackage{leftidx}
\usepackage{accents}
\usepackage{appendix}
\usepackage{authblk}
\usepackage{cases}

\definecolor{slightblue}{rgb}{.8, .8, 1}
\definecolor{hair}{RGB}{100,225,190}
\definecolor{ruby}{RGB}{220,50,120}
\definecolor{grass}{RGB}{150,220,110}
\definecolor{ceruleanblue}{rgb}{0.16, 0.32, 0.75}
\definecolor{deepcarmine}{rgb}{0.66, 0.13, 0.24}
\definecolor{otterbrown}{rgb}{0.4, 0.26, 0.13}
\definecolor{sapphire}{rgb}{0.03, 0.15, 0.4}

\usepackage[colorlinks=true,  linkcolor=otterbrown, citecolor=sapphire,  urlcolor=hair!80!black]{hyperref}

\usepackage[T1]{fontenc}

\newtheorem{theorem}{Theorem}[section] \newtheorem{lemma}[theorem]{Lemma}
\newtheorem{proposition}[theorem]{Proposition}

\theoremstyle{definition} 
\newtheorem{definition}[theorem]{Definition}

\newtheorem{remark}[theorem]{Remark} \numberwithin{equation}{section}
\numberwithin{figure}{section}

\newcommand{\Cb}{\mathbb{C}}

\newcommand{\Hb}{\mathbb{H}}

\newcommand{\Pb}{\mathbb{P}}
\newcommand{\Qb}{\mathbb{Q}}
\newcommand{\Rb}{\mathbb{R}}

\newcommand{\Tb}{\mathbb{T}}

\newcommand{\Zb}{\mathbb{Z}}

\newcommand{\Ac}{\mathcal{A}}
\newcommand{\Bc}{\mathcal{B}}

\newcommand{\Dc}{\mathcal{D}}

\newcommand{\Fc}{\mathcal{F}}
\newcommand{\Gc}{\mathcal{G}}
\newcommand{\Hc}{\mathcal{H}}

\newcommand{\Jc}{\mathcal{J}}
\newcommand{\Kc}{\mathcal{K}}

\newcommand{\Pc}{\mathcal{P}}
\newcommand{\Qc}{\mathcal{Q}}
\newcommand{\Rc}{\mathcal{R}}
\newcommand{\Sc}{\mathcal{S}}
\newcommand{\Tc}{\mathcal{T}}
\newcommand{\Uc}{\mathcal{U}}
\newcommand{\Vc}{\mathcal{V}}
\newcommand{\Wc}{\mathcal{W}}
\newcommand{\Xc}{\mathcal{X}}
\newcommand{\Yc}{\mathcal{Y}}
\newcommand{\Zc}{\mathcal{Z}}

\usepackage{sidecap,subcaption}
\usepackage[section]{placeins}
\usepackage{wrapfig}
\usepackage{lpic}

\sidecaptionvpos{figure}{c}

\newcommand{\wt}{\widetilde}

\DeclareFontFamily{OMX}{yhex}{}
\DeclareFontShape{OMX}{yhex}{m}{n}{<->yhcmex10}{}
\DeclareSymbolFont{yhlargesymbols}{OMX}{yhex}{m}{n}
\DeclareMathAccent{\wideparen}{\mathord}{yhlargesymbols}{"F3}

\makeatletter
\newcommand*\rel@kern[1]{\kern#1\dimexpr\macc@kerna}
\newcommand*\wb[1]{%
	\begingroup
	\def\mathaccent##1##2{
		\rel@kern{0.8}%
		\overline{\rel@kern{-0.8}\macc@nucleus\rel@kern{0.2}}%
		\rel@kern{-0.2}%
	}%
	\macc@depth\@ne
	\let\math@bgroup\@empty \let\math@egroup\macc@set@skewchar
	\mathsurround\z@ \frozen@everymath{\mathgroup\macc@group\relax}%
	\macc@set@skewchar\relax
	\let\mathaccentV\macc@nested@a
	\macc@nested@a\relax111{#1}%
	\endgroup
}
\makeatother

\newcommand{\dist}{\mathrm{dist}}

\title{Sharp asymptotics for arm probabilities in critical planar percolation}
	\date{Preliminary draft, \today}
	\author[1]{Hang Du\thanks{\href{mailto:duhang@pku.edu.cn}{duhang@pku.edu.cn}}}
	\author[1]{Yifan Gao\thanks{\href{mailto:gaoyif@pku.edu.cn}{gaoyif@pku.edu.cn}}}
	\author[1]{Xinyi Li\thanks{\href{mailto:xinyili@bicmr.pku.edu.cn}{xinyili@bicmr.pku.edu.cn}}}
	\author[2]{Zijie Zhuang\thanks{\href{mailto:zijie123@wharton.upenn.edu}{zijie123@wharton.upenn.edu}}}
	\affil[1]{Peking University}
	\affil[2]{University of Pennsylvania}
	
\begin{document}
	\maketitle
	\begin{abstract}
In this work, we consider critical planar site percolation on the triangular lattice and derive sharp estimates on the asymptotics of the probability of half-plane $j$-arm events for $j\geq 1$ and planar (polychromatic) $j$-arm events for $j>1$. These estimates greatly improve previous results and in particular answer (a large part of) a question of Schramm ({\it ICM Proc.}, 2006). 
\end{abstract}
	\section{Introduction}
Percolation is without doubt one of the most studied statistical mechanics models in probability. As an ideal playground for the study of phase transitions and criticality, it has received considerable attention from probabilists and statistic physicists in the past 60+ years. Despite its simple setup, it is the source of many fascinating yet difficult mathematical problems, with some already well answered and many more still very far from being solved. Starting from the beginning of 21st century, there have been a lot of breakthroughs in the study of a particular case of this model, namely the critical planar percolation on the triangular lattice. In the ground-breaking work by Smirnov \cite{smirnov2001critical}, it is shown that the both crossing probabilities and the exploration process have conformally invariant scaling limits that can be given through Cardy's formula and described as a (chordal) Schramm-Loewner evolution with parameter 6 (${\rm SLE}_6$) respectively. Later, in \cite{camia2006two} Camia and Newman give the characterization of the full scaling limit, via the collection of non-self-crossing loops that is known now as the conformal loop ensemble (CLE). To keep the introduction concise, we refer readers to the classical book \cite{grimmett1999percolation} and a recent survey \cite{duminil2018sixty} for more on the history and some recent progresses of this model.

Critical exponents are central notions in understanding the behavior of critical models. In the case of percolation, many of them can be derived from the so-called arm exponents, i.e.\ exponents in the power-law decay of the probability of arm events (as the mesh size tends to 0), which we will introduce shortly. For critical planar percolation on the triangular lattice, Smirnov and Werner calculate the precise values of the half-plane arm exponents and (polychromatic) $j$-arm exponents with $j>1$ in \cite{MR1879816}, using Kesten's scaling relation from \cite{MR879034}, Cardy's formula as well as relations between interfaces and ${\rm SLE}_6$. Roughly in the same time, Lawler, Schramm and Werner obtain the one-arm exponent in the plane in \cite{lawler2002one}, also via SLE-related calculations. We will briefly recall their results in Claim (3) of Lemma \ref{prior_1} and \eqref{eq:planeonearm} respectively. 

However, these asymptotics (with an $o(1)$ in the exponent) are not completely satisfactory. One particular reason is that in dealing with scaling limits involving microscopic quantities from critical planar percolation, such as the pivotal, cluster or interface measures from \cite{MR3073882} (see in particular Theorems 4.3, 5.1 and 5.5, ibid.) and the natural parametrization of the interface (see Theorem 1.4 in \cite{HLS}), instead of renormalizing with an explicit factor $\eta^{\alpha}$ where $\eta$ is the mesh size and $\alpha$ the corresponding exponent, one normalizes by something implicit, namely the asymptotic arm probability (which corresponds to $a_j$ in this paper); see also Remark 4.10 in \cite{MR3073882} for more on this issue.

Thus, it is natural to ask if there are precise estimates of the arm events, just as what have been obtained for other critical models, e.g., sharp asymptotics of the one-point function of the loop-erased random walk (LERW); see \cite{kenyon2000asymptotic}, \cite{li2019onepoint} and \cite{lawler2019four} for the two-, three- and four-dimensional cases respectively.

In fact, much before the scaling limit results mentioned in the paragraph above are obtained, Schramm already asks the question in the proceedings of ICM 2006 if it is possible to improve the estimates. In Problem 3.1 of \cite{schramm2007conformally}, he points out that ``it would be especially nice to obtain estimates that are sharp up to multiplicative constants''.

In a few special cases\footnote{Where the exponent is an integer, namely the half-plane 2- and 3-arm exponents and planar 5-arm exponent.}, both the value and up-to-constants estimates can be obtained without relating to SLE's; see Lemma 5 of \cite{kesten1998almost}, Theorem 24 of \cite{nolin2008near} or the first exercise sheet of \cite{werner2007lectures}.
Improvements for other arm  asymptotics are much more difficult. Mendelson, Nachmias and Watson obtain in \cite{mendelson2014rate} a power-law rate for the Cardy's formula\footnote{Also independently obtained roughly at the same time by Binder, Chayes and Lei in \cite{binder2015rate}.} which yields a slight improvement of half-plane one-arm asymptotics (see \eqref{eq:MNW}), although still not as strong as up-to-constants estimates.

In this paper, we answer a large part\footnote{The exception is the planar one-arm case. See Remark \ref{rmk:planeonearm} for more detail.} of Schramm's question by deriving sharp estimates for the probability of half-plane arm events and planar $j$-arm events with $j>1$, with power-law error bounds in most cases, which gives much more than what he asks for.

\subsection{Main results}
We start with necessary notation. 
Let $\Tb$ denote the triangular lattice where each face is an equilateral triangle and $\Tb^*$ denote the dual graph. We consider the critical Bernoulli percolation on $\Tb^*$ in which each hexagon is colored red (=open) or blue (=closed) independently with equal probability. For $j\geq 1$, let $\Bc_j(r,R)$ denote the half-plane $j$-arm event that there exist $j$ disjoint crossings with alternating colors of the semi-annulus $A^+(r,R)$ in the upper half-plane with inner and outer radii $r$ and $R$ respectively (see Section \ref{subsec:arm-event} for the precise definition of arm events and Section \ref{subsec:setting} for the discretization of domains). We will also consider a variant of half-plane arm events, denoted by $\Hc_j(r,R)$, which corresponds to $j$-alternating arms from the segment $[-r,r]$ to the semi-circle $C_R^+$ within the half-disk $B_R^+$. Write $b_j(r,R): = \Pb[\Bc_j(r,R)]$ and define $h_j$ similarly. In the seminal work \cite{MR1879816}, it is showed that for any $j \geq 1$, the sequence $b_j(r,R)$ has a power-law decay in $R$ with exponent 
\begin{equation}\label{eq:betadef}
\beta_j=j(j+1)/6,\quad\mbox{such that as $R \to \infty$, }\quad b_j(r,R)=R^{-\beta_j+o(1)}\,.
\end{equation}
In \cite{mendelson2014rate} the following improvement of $b_1$ is obtained:
\begin{equation}\label{eq:MNW}
b_1(1,n)=e^{O(\sqrt{\log \log n})} n^{-1/3}=(\log n)^{O(1/\sqrt{\log\log n})}n^{-1/3}.
\end{equation}

We are now ready to state our main results for half-plane arm probabilities, which give sharp asymptotics with power-law error bounds for both $b_j$ and $h_j$.
\begin{theorem}
\label{thm:newes-half}
For any $j \geq 1$, 
there exist constants $c(j,r)>0$ and $C_b(j,r),C_h(j,r)\geq 0$ such that for any real $n>r$,
\begin{equation}
\label{eq:newes-half}
b_j(r,n)=C_b n^{-\beta_j} \Big(1+O(n^{-c})\Big);\quad h_j(r,n)=C_h n^{-\beta_j}\Big(1+O(n^{-c})\Big)\,.
\end{equation}
\end{theorem}
Note that in the cases of $j=1,2,3$, $C_b(j,1),C_h(j,1)>0$, hence it makes sense to pick $r=1$, and our theorem applies to classical half-plane arm events out of one hexagon. The same asymptotics also hold for any fixed ``inner initial configuration''. In this case, the constants depend on $j$ and the initial configuration.

\medskip

We now turn to (polychromatic) planar $j$-arm events. Again, see Sections \ref{subsec:setting} and \ref{subsec:arm-event} for conventions and precise definitions. We start with the classical $j$-arm event $\Pc_j(r,R)$, which is defined as the event that there exist $j$ disjoint crossings of the annulus $A(r,R)$ and not all of the same color (except in the case $j=1$). We will also work with several variants of arm events, with requirements on color sequences, location constraints on the outer ending points of arms and/or the connectedness of certain pairs of arms. In particular, for $j>1$, we consider $\Xc_j$ and $\Ac_j$, which are variants satisfying the requirement on color sequences, with and without location constraints respectively; see \eqref{def:a} and below for precise definition, as well as $\Yc_j$ are $\Zc_j$, which are defined similarly but with the extra requirement on the connectedness; see \eqref{def:y} and below for precise definitions. In the cases $j=2,3,4,5,6$, $\Ac_j(1,R)$ are the classical alternating $j$-arm events.

Write $p_j(r,R) = \Pb[\Pc_j(r,R)]$ and define $a_j$, $x_j$, $y_j$ and $z_j$ similarly. In  \cite{MR1879816}, it is showed that for any $j \geq 2$, the sequence $p_j(r,R)$ has a power-law decay in $R$ with critical exponent 
\begin{equation}\label{eq:alphadef}
\alpha_j = (j^2-1)/12,\quad\mbox{such that as $R \to \infty$, }\quad p_j(r,R)=R^{-\alpha_j+o(1)}\,.
\end{equation}

We now state our main results for planar arm probabilities which give sharp asymptotics with power-law bounds for variants $x_j$ and $y_j$ and as a consequence sharp asymptotics without explicit error bounds for $a_j$ and up-to-constants estimates for $p_j$.
\begin{theorem}
\label{thm:newes-whole}
For any $j \geq 2$, 
there exist constants $c(j,r)$ and $C_x(j,r),C_y(j,r)\geq 0$ such that for all real $n>r$,
\begin{equation}
\label{eq:newes-whole}
x_j(r,n)=C_x n^{-\alpha_j} \Big(1+O(n^{-c})\Big)\,;\quad y_j(r,n)=C_y n^{-\alpha_j} \Big(1+O(n^{-c})\Big)\,.
\end{equation}
\end{theorem}
Note that when $j=2,3,4,5$, $C_x(j,1),C_y(j,1)>0$, and also when $j=6$, $C_x(6,1)>0$. Hence in these cases (note that special care is needed for $x_6(1,n)$; see Remark \ref{rmk:a6} for discussions) our results hold for arm events out of a single hexagon. 
As an immediate consequence of Theorem \ref{thm:newes-whole} and the up-to-constants asymptotic equivalence of different arm events (see Claim (2) of Lemma \ref{prior_1}),
\begin{equation}\label{eq:roughasymp}
p_j(r,n) \asymp a_j(r,n) \asymp n^{-\alpha_j}.
\end{equation}

Moreover, we are able to obtain sharper estimates (albeit without a power-law error bound) than \eqref{eq:roughasymp} for the alternating arm probability $a_j(r,n)$. For simplicity, we only give them for arm events out of one hexagon for $j\leq 5$ (but similar asymptotics hold for arm events with other inner initial configurations, including the case of $a_6(1,n)$). See Remark \ref{rmk:a6} for more discussions.
\begin{theorem}
\label{thm:newes2-whole}
For $j=2,3,4,5$,  there exists some $C_a(j)>0$ such that for all real $n>1$,
\begin{equation}
\label{eq:newes2-whole}
 a_j(1,n)= C_a n^{-\alpha_j}\Big(1+o(1)\Big).
\end{equation}
\end{theorem}
A direct application of Theorem \ref{thm:newes2-whole} for $j=2,4$ is that one can replace renormalizing factors in the scaling limit of pivotal and interface measures by precise powers of the mesh size in Theorems 4.3 and 5.5 of \cite{MR3073882} and that in the natural parametrization in Theorem 1.4 of \cite{HLS}. Another application is the improvement for the asymptotics of the correlation or characteristic length, which is also a central object in the study of near-critical and dynamical percolation. In particular, our asymptotics on planar 4-arm events imply that various versions of the correlation length for planar critical percolation on triangular lattice are up-to-constants equivalent to $|p-1/2|^{-4/3}$. For a more thorough account on the correlation length, see e.g., Section 7 of \cite{nolin2008near}.

\subsection{Comments}\label{subsec:comments}
In this subsection, we briefly comment on the proof and discuss possible directions for generalization.

We start with the strategy of the proof. As Schramm has already points out immediately below his question in \cite{schramm2007conformally}, the crux of the matter lies in ``the passage between the discrete and continuous setting'' and he poses a related question on obtaining ``reasonable estimates for the speed of convergence of the discrete processes''. 

This question is solved recently by Binder and Richards in \cite{binder2020convergence}, which constitutes the PhD thesis \cite{richards2021convergence} of Richards; see also \cite{binder2021rate} for an extended abstract.  More precisely, they verify that the framework developed in \cite{viklund2015convergence} for the power-law convergence rate of random discrete models towards SLE indeed works for percolation, the harmonic explorer and the FK-Ising model. In particular, they obtain a power-law convergence rate for the exploration process of planar critical percolation up to some stopping time. See Section \ref{subsec:exploration} for detailed discussions of their results.

We now explain how we derive sharp asymptotics for arm events out of the power-law convergence. For better illustration we will discuss the classical half-plane arm-events (i.e., Theorem~\ref{thm:newes-half} for $b_j$) in more detail and only briefly point out necessary modifications in other cases.

At first glance, sharp asymptotics for arm events should follow naturally once we relate the discrete exploration path to ${\rm SLE}_6$. However, a few obstacles prevent us from applying the power-law convergence rate directly. 

The first one is that such results only give information down to the mesoscopic level (in terms of the mesh size), while in this work we endeavor to reach the microscopic level. To overcome this difficulty, one resorts to coupling techniques developed for critical planar percolation, which allow us to ``decouple'' the microscopic ``initial configurations'' (if we figuratively regard arm events as the consequence of  outward explorations) with the macroscopic boundary conditions. Inspired by the arguments of Theorem 1.2 in \cite{li2019onepoint} for the one-point function of 3D LERW, the ideas above can be crystallized in the following  proportion estimates, from which Theorem~\ref{thm:newes-half} for $b_j$ can be directly derived.

Denote $r_\star(j)=\min\{r\in \mathbb{N}:\star_j(r,n)>0\mbox{ for all sufficiently large }n\}$ for $\star\in\{\operatorname{b},\operatorname{h},\operatorname{p},\operatorname{a},\operatorname{x},\operatorname{y},\operatorname{z}\}$. It is clear that when dealing with some arm probability $\star_j(r,n)$, it suffice to focus on the case $r\ge r_\star(j)$. Note that $r_b(j)\le r_h(j)$ and $r_a(j)=r_x(j)\le r_y(j)=r_z(j)$; in addition, $r_h(j)=1$ for $ j\le 3$ and $r_y(j)=1$ for $j\le 5$.
\begin{proposition}
\label{prop:ratioestimate-half}
Given $j \geq 1$, for any $r\geq r_b(j)$ and $m\in (1.1,10)$,
\begin{equation}
\label{eq:ratioestimate-half}
\frac{b_j(r,m^2 n)}{b_j(r,mn)}=\frac{b_j(r,mn)}{b_j(r,n)} \Big(1+O(n^{-c})\Big)\,,
\end{equation}
where $O(n^{-c})$ is independent of the choice of $m$.
\end{proposition}

The second one is that the classical arm event $\Bc_j$ is not the ideal choice to be described by the exploration path in a discretized domain. More precisely, working with $\Bc_j$ directly involves dealing with a sequence of semi-annuli with shrinking inner radii. Although in principle one should be able to produce similar power-law convergence in these varying domains as a perturbation of the half-disk, a more sensible choice is to make use of couplings of percolation configurations conditioned on arm events (see the discussion on the proof of Proposition~\ref{prop:ratioestimate-half} below), and work with the variant $\Hc_j$ instead; see Proposition \ref{prop:super-strong} for more detail.

The third one is that the result of Binder and Richards provides convergence rate only for the exploration up to some stopping time, not for the whole path (although it is believed to hold for the latter as well). To overcome this difficulty, we enlarge the domain in which we compare percolation explorations and ${\rm SLE}_6$, so that the segment before the stopping time already suffices to provide comparison of arm probabilities; see Section \ref{subsec:4.1} for more detail.

We now discuss the proof of Proposition~\ref{prop:ratioestimate-half} in more detail. The main ingredients are Propositions~\ref{prop:couple-h-half} and \ref{prop:couple-bh-half}, in which we derive comparisons of conditioned arm probabilities, which allow us to compare arm probabilities between different boundary conditions and also to compare different types of arm events. These estimates are established through coupling results established in Propositions \ref{prop:one-coupling} and \ref{prop:coup-1} for $j=1$ and $j>1$ respectively as well as the power-law convergence of the exploration proces in the form of Proposition \ref{prop:couple-SLE-half}.

A crucial ingredient for these couplings is the separation lemma, which, first developed by Kesten in \cite{MR879034}, plays a key role in establishing quasi-multiplicativity properties of arm probabilities. We will discuss it (and several variants) in detail in Sections \ref{subsec:seplemmahalf} and \ref{subsec:seplemmawhole}. It is also worth mentioning that although not directly needed in the proof, we also obtain a super-strong separation lemma for the half-plane setup, which confirms a conjecture by Garban, Pete and Schramm in \cite{MR3073882} and is of independent interest.

\smallskip

We now turn to the plane case. Similarly, Theorem~\ref{thm:newes-whole} for $y_j$ can be derived from the following 
proportion estimates. 
\begin{proposition}
\label{prop:ratioestimate-whole}
Given $j \geq 2$,  for any $r\geq r_y(j)$ and $m\in (1.1,10)$,
\begin{equation}
\label{eq:ratioestimate-whole}
\frac{y_j(r,m^2 n)}{y_j(r,mn)}=\frac{y_j(r,mn)}{y_j(r,n)} \Big(1+O(n^{-c})\Big)\,,
\end{equation}
where $O(n^{-c})$ is independent of the choice of $m$.
\end{proposition}
This proposition follows from the same argument as that of Proposition~\ref{prop:ratioestimate-half} as the specific definition of $\Yc_j$ allows us to relate it to the exploration process.  
Note that in order to overcome the third obstacle as in the half-plane case, we will also enlarge the domain accordingly; see Section \ref{subsec:4.3} for more detail. The claim for $x_j$ follows from a coupling result that relates $y_j$ and $x_j$; see Proposition \ref{prop:couple-a-whole}.

Theorem~\ref{thm:newes2-whole}, whose proof is inspired by Proposition 4.9 of \cite{MR3073882}, is a corollary of Theorem \ref{thm:newes-whole}, Proposition \ref{prop:couple-A-whole-inner} (which relates the conditional arm probabilities for $\Xc_j$ and $\Ac_j$), and Claim (4) of Lemma \ref{prior_1} (which gives the convergence of macroscopic arm probabilities without a speed).

\begin{remark}\label{rmk:planeonearm} We now briefly mention some open questions and discuss possible generalizations.

\smallskip

\noindent 1) Planar 1-arm events are crucial objects to understand critical percolation clusters. In \cite{lawler2002one} it is shown that the 1-arm probabilities satisfy
\begin{equation}\label{eq:planeonearm}
p_1(r,R)=R^{-5/48+o_R(1)}\quad\mbox{ as }R\to\infty.
\end{equation}
It is then a natural question to wonder if our asymptotics also hold in this case. In fact, the key difficulty in extending our arguments to the one-arm case lies in the fact that planar one-arm events cannot be described by a single chordal exploration path. Instead, they can be described either by the collection of all interface loops or by the so-called radial exploration process, which correspond to ${\rm CLE}_6$ or the radial ${\rm SLE}_6$ respectively in the scaling limit. While we believe power-law convergence analogous to Theorem 4.1.11 in \cite{richards2021convergence} should also hold for these objects, solid arguments leading to such results still seem rather out of reach for the moment.

\smallskip

\noindent 2) Greater difficulties exist for the (planar) monochromatic arm-events. In \cite{beffara2011monochromatic}, Beffara and Nolin show the existence of monochromatic arm exponents for critical planar percolation using quasi-multiplicativity arguments. Although couplings of conditional monochromatic arm events should exist in some form (if so, it must be proved in a way different than ours), the lack of an explicit characterization through the exploration path is the main difficulty that prevents us from obtaining any new results.

\smallskip

\noindent 3) One may also wonder if our results hold for other lattices and/or types of percolation, in particular the critical bond percolation on $\Zb^2$. In this direction, we believe that all our coupling arguments (including the separation lemmas) will work with very little adaptation, while the biggest obstacle is definitely the passage from discrete to continuum. See \cite{duminil2020rotational} for some recent progress in this direction.

\smallskip

\noindent 4) Another interesting direction is to  obtain sharp asymptotics also for the FK-Ising arm events and the one-point function of the harmonic explorer (as well as the convergence in natural parametrization to ${\rm SLE}_4$), for Binder and Richards also establish power-law convergence for these two models in \cite{binder2020convergence}. We plan to investigate them in future works.
\end{remark}

Finally, we explain how this article is organized. In Section \ref{sec:notation}, we settle the setup and introduce various key notions and preliminary results, including arm events, faces,  separation lemmas, and the power-law convergence rate of the exploration paths. We give without proof various coupling results concerning arm events in Section \ref{sec:coupling} which are crucial to the main arguments and derive out of these couplings some relations on conditioned arm probabilities. Section \ref{sec:compare} is dedicated to the proof of the main theorems in which we also include the proofs of Propositions \ref{prop:ratioestimate-half} and \ref{prop:ratioestimate-whole}. In Section \ref{sec:proofcoupling}, we unfold the proofs for various coupling results postponed from Section \ref{sec:coupling}. A few technical proofs for preliminary results in Section \ref{sec:notation} are given in the Appendices.

\medskip

\noindent {\bf Acknowledgments}: 
Hang Du, Yifan Gao and Xinyi Li thank National Key R\&D Program of China (No.\ 2021YFA1002700 and No.\ 2020YFA0712900) and NSFC (No.\ 12071012) for support. Hang Du is also partially supported by the elite undergraduate training program of School of Mathematical Sciences at Peking University. Zijie Zhuang is partially supported by NSF grant DMS-1953848. Xinyi Li also thanks Daisuke Shiraishi for inspiring discussions.

	\section{Notation and preparatory results}
\label{sec:notation}
This section is dedicated to setup and preparatory results. In Section \ref{sec:prelimnotations} we introduce overall notation. We  then fix the setup and recall some basic tools for percolation in Sections \ref{subsec:setting}--\ref{subsec:basictools} resp. In Section \ref{subsec:arm-event}, we give the definition for various arm events, state some existing results their its asymptotics and give a ``functional equation'' type result  on sequences that will be useful for obtaining sharp asymptotics.  In Sections \ref{subsec:seplemmahalf} and \ref{subsec:seplemmawhole}, we state and prove the separation lemmas in the half-plane and plane resp., which are crucial to the coupling argument in this work. Finally in Section \ref{subsec:exploration}, we recall the power-law convergence rate of the rescaled exploration path towards ${\rm SLE}_6$ from \cite{richards2021convergence} and prove a variant tailored for this work.

\subsection{Notation and conventions}\label{sec:prelimnotations}
Let $\Cb$ stand for the complex plane and $\Hb=\{ x+iy: y>0\}$ stand for the upper half-plane. For convenience of notation we regard them as $\mathbb{R}^2$ and a subset thereof and use both sets of notation interchangeably.
Let $B(x,R)=\{z: |z-x|< R\}$ denote the ball of radius $R$ around $x$ and $C(x,R)=\partial B(x,R)$. For $0<r<R$, we write $A(x,r,R):=\{ z: r<|z-x|< R \}$ for the annulus of radii $r<R$ around $x$. We will omit ``$x$'' when $x=0$, i.e., abbreviate $B(0,R),C(0,R),A(0,r,R)$ as $B_R, C_R, A(r,R)$. We will add $+$ in the superscript to indicate the quantities are defined in the upper half-plane $\Hb$. For example, if $A\subseteq\Cb$, write $A^+$ for the set $A\cap\Hb$. Then, $B^+(x,R), C^+(x,R), A^+(x,r,R), B^+_R,C^+_R,A^+(r,R)$ are the corresponding subsets of $\Hb$. Without further specification, all the sequences we consider are indexed by real numbers, not only integers. 

We write $a_n=n^{\alpha+o(1)}$ if for all $\epsilon>0$, $n^{\alpha-\epsilon} \leq a_n \leq n^{\alpha+\epsilon}$ for all $n$ large enough. We write $a_n=O(n^{\alpha})$ if there exists $C>0$ such that $|a_n| \leq Cn^{\alpha}$ for all $n$, and write $a_n\asymp b_n$ if $a_n,b_n>0$ for all $n$ large enough, and $a_n/b_n$ is bounded both from above and below. We always write $j$ for the number of arms and $\eta$ for the scale of lattice. Without further specification, all the arcs are written counterclockwise. We write $\operatorname{d}(x,y)$ for the Euclid distance of $x,y\in \mathbb{R}^2$. Given a real number $x$, we write $\lfloor x \rfloor:=\max_{z\in\mathbb{Z}}\{z\leq x\}$ for the integer part of $x$.

Finally, let us explain our convention concerning constants. Constants like $\varepsilon$, $\delta$, $c$, $c'$, $C$ or $C'$  may change from line to line while those with a subscript like $c_1$ are kept fixed throughout the paper. All constants are universal except those marked at the first occurrence.
\subsection{Setting for percolation}\label{subsec:setting}
Let $\Tb$ be the triangular lattice embedded in the complex plane $\Cb$ where each face is an equilateral triangle of side length $1$. More precisely, the set of sites (vertices) in $\Tb$ is given by $T:=\Zb+e^{i\pi/3}\Zb$, and two sites are neighbors if they are at distance $1$.

We consider the critical site percolation on $\Tb$, defined through declaring each site $v\in T$ as open or closed equally with probability $\frac12$, independently of the other sites. More precisely, let $\{0,1\}^{T}$ be the sample space of configurations $(\omega_v)_{v\in T}$, where $\omega_v=1$ if $v$ is open, and $\omega_v=0$ if $v$ is closed. Write $\Pb$ for the product measure with parameter $\frac12$ on $T$. For better illustration, a site will be colored red (resp.\ blue) if it is open (resp.\ closed). In studying scaling limits we will also be interested in percolation on the rescaled lattice. For $\eta>0$, let $\eta\Tb$ be $\Tb$ rescaled by $\eta$. With slight abuse of notation, we will also use $\Pb$ to denote the product measure with parameter $\frac12$ on $\eta T$. However, we will fix $\eta=1$ in the majority of this work and  emphasize the mesh size only when\footnote{Namely in Subsections \ref{subsec:exploration}, \ref{subsec:4.1} and \ref{subsec:4.3}.} we do rescale the lattice.

We also view the percolation on the triangular lattice $\Tb$ as a random coloring of {\bf hexagons} (i.e.\ faces) on $\Tb^*$, which is dual to $\Tb$ such that each vertex (resp.\ edge) on $\Tb$ corresponds to a hexagonal face (resp.\ an edge) on $\Tb^*$.
We say two hexagons are \textbf{neighbors} (or adjacent) if they have a common edge. A \textbf{path} on $\Tb^*$ is a sequence of neighboring hexagons and they are all distinct, and we call it a \textbf{loop} if the two ends of the path are also neighbors. 
We will also consider paths on $\Tb^*$ as a graph, which we refer to as {\bf b-paths} to distinguish them from paths of hexagons we have just defined\footnote{In the language of this work, arms are paths while interfaces are b-paths.}. We say two interfaces are \textbf{adjacent} if there is a hexagon on $\Tb^*$ touched by both of them. 

Let $\partial_1$ and $\partial_2$ be two disjoint parts of the boundary of a hexagonal domain $D$  (in many cases two opposite sides of a topological quadrangle or inner and outer boundaries of a topological annulus). A \textbf{crossing} from $\partial_1$ to $\partial_2$ in $D$ is a monochromatic path of hexagons in $\eta\Tb^*$ such that the two ends intersect $\partial_1$ and $\partial_2$ respectively, and all other hexagons are inside $D$. A crossing will also be called an \textbf{arm} when we consider arm events in what follows.

We now turn to the issues related to discretization of domains. By default, we use the conventions in \cite{camia2006two} and refer readers thereto for more details. However, to avoid inconsistency for the discretization of (half)-disks, (semi-)annuli and (semi-)circles that are ubiquitous in this work, we give a special rule for the discretization of these objects; see the last paragraph of this subsection for details. We also remark that as will be discussed in Remark \ref{rmk:insensitivity}, our results are indeed quite stable against different choices of discretization. Hence our convention here is in fact merely an inessential technical matter.

We start with the default convention. A set of hexagons $E$ is called connected or simply connected if so is the set obtained by embedding $E$ in $\Cb$. If $E$ is a simply connected set of hexagons, we write $\Delta E$ for the outer boundary of $E$, i.e., the set of hexagons that are not in $E$ but adjacent to hexagons in $E$, and we write $\partial E$ for the topological boundary of $E$ by viewing $E$ as a domain in $\Cb$. The bounded simply connected set $E$ of hexagons is called a \textbf{Jordan set} if $\Delta E$ is a loop.

	Given a Jordan domain $D$ in $\Cb$, we write $D_{[\eta]}$ for the $\eta$-approximation of $D$, i.e., the largest Jordan set of $\eta\Tb^*$ that is contained in $D$. For $a \in \partial D$, we write $a_\eta$ for the vertex of $\eta\Tb^*$ in $\partial D_{[\eta]}$ such that it is closest to $a$ and it is on the common edge shared by two adjacent hexagons in $\Delta D_{{[\eta]}}$ (if there is a tie, choose one arbitrarily). When $\eta=1$, we omit $\eta$ or $[\eta]$ from the subscript unless we wish to emphasize that this is a discretized object.

We now turn to special discretization rules for points and intervals on the real axis, circles, disks and annuli centered at the origin. Note that if one applies the rules in the paragraph above to e.g., $B^+_r$ and $A(r,R)$ which are two touching objects in the continuum, the resulting hexagonal domains will have non-empty ``cracks'' in-between. Hence, we apply the following special rule for them.

We first designate, in an arbitrary way, the discrete circles\footnote{With slight abuse of notation we still denote by $C_R$ the discretized object; same for other objects in this paragraph.} $C_R$ for each $R\geq 1$ as a loop of b-paths on $\Tb^*$ that disconnects the origin from infinity which is at most $O(1)$-away from $C_R$, such that $C_1$ is the boundary of the hexagon containing the origin and there do not exist $R_1<R_2$ such that $C_{R_1}$ crosses $C_{R_2}$. Then, for each $1\leq r<R$ we regard $B_R$ as the hexagonal domain encircled by $C_R$ and $A(r,R)$ as $B_R\setminus B_r$ respectively. Similarly, we also designate for each $R\geq 1$ the discrete semi-circle $C^+_R$ as a b-path that disconnect $0_{1}$ from infinity in $\Hb_{[1]}$ with the special rule that $C^+_1$ is the upper three edges of the hexagonal face in $\Tb^*$ containing $0_1$. We then define half-disks $B_R^+$ and semi-annuli $A^+(r,R)$ in a similar fashion. With slight abuse of notation, we write $[-R,R]$ for the b-path that forms the part of the boundary of $\Hb_{[1]}$ encircled by $C^+_R$ and define other intervals accordingly. When working on the discretization with respect to the rescaled lattice $\eta\Tb$, we first counter-rescale to $\eta=1$, apply the rules above and then rescale everything back to $\eta\Tb$.

	\subsection{Some classical tools}\label{subsec:basictools}
	The aim of this subsection is to review some classical tools in percolation theory, namely, the Harris-FKG inequality, the BK-Reimer inequality and the Russo-Seymour-Welsh theory. These tools will be extensively used throughout this paper. We refer readers to \cite{grimmett1999percolation} for more details.
	
	\begin{definition}[Monotone events]
		There is a natural partial order on the set of configurations $\{ 0,1 \}^T$ given by $\omega_1\le \omega_2$ if and only if $\omega_1(v)\le \omega_2(v)$ for all $v\in T$. An event $A$ is called increasing if $1_{A}(\omega_1)\le1_{A}(\omega_2)$ whenever $\omega_1(v)\le \omega_2(v)$, where $1_A$ is the indicator function of $A$. An event $A$ is called decreasing if its complement $A^c$ is increasing.
	\end{definition}
	
	\begin{lemma}[Harris-FKG inequality]
		If $A$ and $B$ are both increasing or decreasing events, then
		\[
		\Pb[A\cap B]\ge \Pb[A] \Pb[B].
		\]
	\end{lemma}
	
	\begin{lemma}[BK-Reimer inequality]\label{lem:Reimer}
		Let $A\Box B$ denotes the disjoint occurrence of the events $A$ and $B$, meaning that $\omega\in A\Box B$ if there exist disjoint sets of sites $E$ and $F$ (possibly depending on $\omega$) with the property that one can verify that $\omega\in A$ (resp. $\omega\in B$) by looking at sites in $E$ (resp. $F$) only. Then, we have
		\[
		\Pb[A\Box B]\le \Pb[A] \Pb[B].
		\]
	\end{lemma}
	
	\begin{lemma}[Russo-Seymour-Welsh (RSW) Theorem]
		For any topological quadrangle (i.e., a Jordan domain with four marked points on the boundary) $D$ in $\Cb$ with two opposite sides $\partial_1$ and $\partial_2$, there exist constants $\eta_0>0$ and $c=c(D)>0$ such that for all $\eta\le\eta_0$, with probability at least $c$, there is an open crossing in $\eta\Tb^*$ from $\partial_1$ to $\partial_2$ in $D_{[\eta]}$ (which we call a quad crossing).
	\end{lemma}

The combination of the Harris-FKG inequality and RSW theorem allows us to ``glue'' arms of the same color to create paths that have macroscopic geometric restrictions with uniformly positive probability. E.g., a loop that disconnects the annulus $A(R,2R)$ can be constructed through gluing of quad crossings. We will colloquially refer to these arguments as ``RSW-FKG gluing'' or ``applications of the RSW theory''.

\subsection{Arm events and asymptotics}\label{subsec:arm-event}
In this subsection, we briefly review the definition of classical arm events in planar percolation and introduce variants that we will specifically use in this work. We then review existing asymptotics of the arm probabilities and finally lay out a ``functional equation'' lemma on sequences that is tailored for the derivation of sharp asymptotics in our work.

\begin{figure}
	\centering
	\includegraphics[height=.2\textwidth]{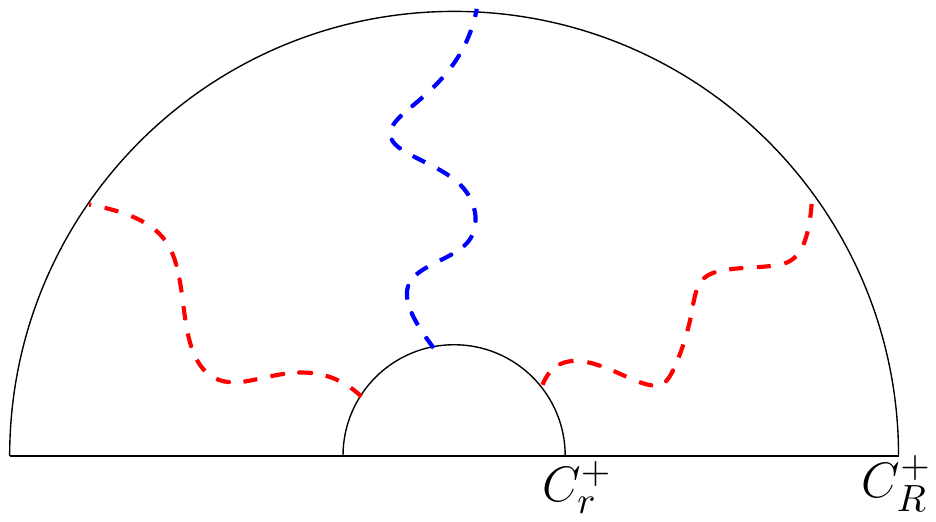}\qquad \qquad
	\includegraphics[height=.2\textwidth]{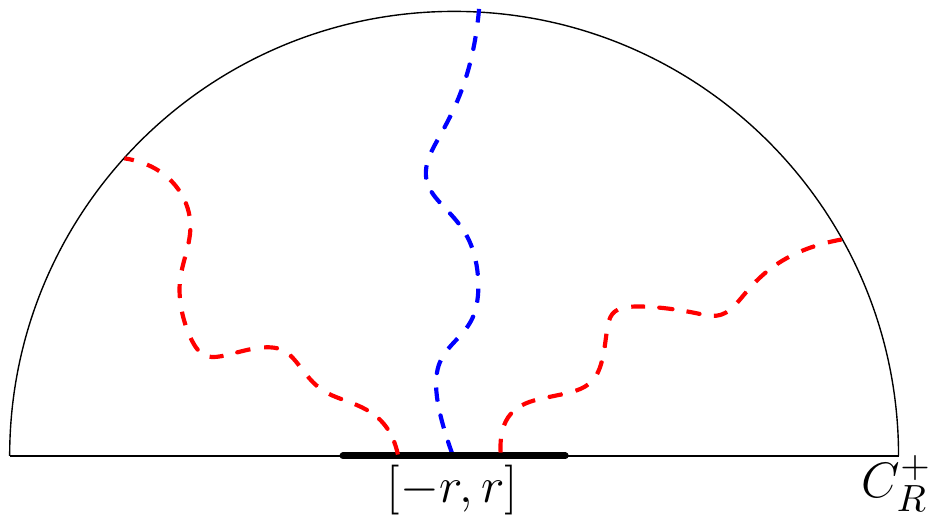}
	\caption{{\bf Left:} The event $\Bc_3(r,R)$. {\bf Right:} The event $\Hc_3(r,R)$. The semi-annulus $A^+(r,R)$ and the line segment $[-r,r]$ (in bold black) should be viewed as b-paths discretized according to the convention given at the end of Section \ref{subsec:setting}. Dashed red and blue curves represent arms (composed of hexagons) with respective colors.}
	\label{fig:BH}
\end{figure}

We first introduce two types of half-plane arm events.
We fix $\eta =1 $ throughout these definitions. In the following, the symbols $B^+_R, C^+_R$, $A^+(r,R)$, $[-r,r]$, $C_r$ and $A(r,R)$ refer to discrete objects introduced at the end of Section \ref{subsec:setting}. For $j\geq 1$, let $\Bc_j(r,R)$ and $\Hc_j(r,R)$ stand for the half-plane $j$-arm events from scale $r$ to scale $R$ in the semi-annulus and half-disk respectively, defined as 
\begin{equation}
	\label{def:b}
	\Bc_j(r,R) = \{ \mbox{There are } j \mbox{ disjoint arms from } C^+_r \mbox{ to }C^+_R \mbox{ in } A^+(r,R) \mbox{ of alternating colors} \};
\end{equation}
\begin{equation}
	\label{def:h}
	\Hc_j(r,R) = \{ \mbox{There are } j \mbox{ disjoint arms from } [-r,r] \mbox{ to }C^+_R  \mbox{ in } B^+_R  \mbox{ of alternating colors} \},
\end{equation}
where \textbf{alternating colors} means that the color pattern is red, blue, red, \ldots, in counterclockwise order. See Figure~\ref{fig:BH} for an illustration.
\begin{remark}\label{rem:HjnotBj} We now briefly comment on the definitions above.

\smallskip

\noindent 1) Among the two definitions above, $\Bc_j(r,R)$ is the classical half-plane arm event in the literature, however it is not the convenient setup for us to relate to the scaling limit of the exploration process.  To overcome this difficulty, we introduce $\Hc_j(r,R)$, which perfectly solves this problem. See Section \ref{subsec:4.1}, in particular Lemma~\ref{lemma4.2-1} for how the  definition of $\Hc_j(r,R)$ comes into play.

\smallskip
	
\noindent 2) As discussed in Remark 2 in \cite{MR1879816}, in the half-plane case it is not restrictive to study arm events of alternating colors thanks to a ``color-switching'' trick. In fact, the probability of $\Bc_j(r,R)$ (or $\Hc_j(r,R)$) remains the same for any color sequence.
\end{remark}

\medskip

\begin{figure}
	\begin{subfigure}{.24\textwidth}
		\centering
		\includegraphics[width=0.9\linewidth]{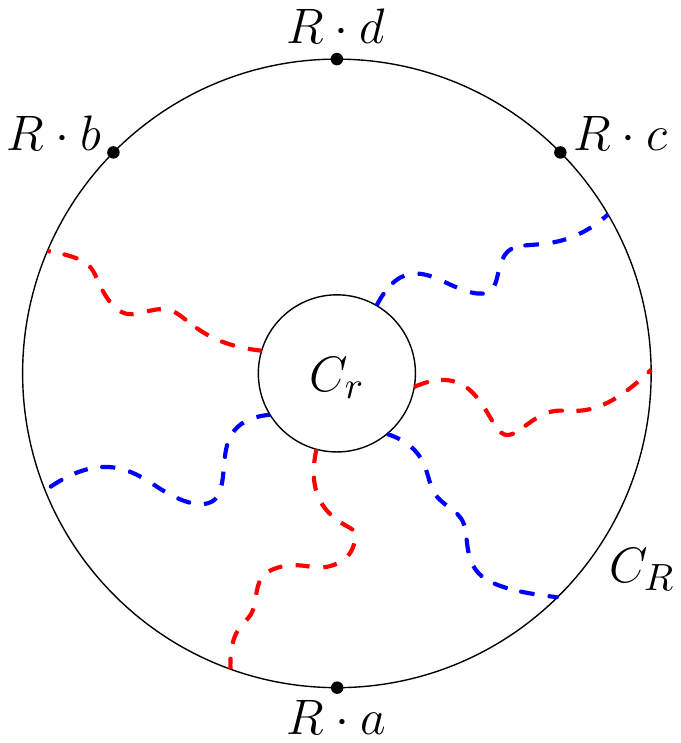}
		\caption{$\Xc_6(r,R)$.}
	\end{subfigure}
	\begin{subfigure}{.24\textwidth}
		\centering
		\includegraphics[width=0.9\linewidth]{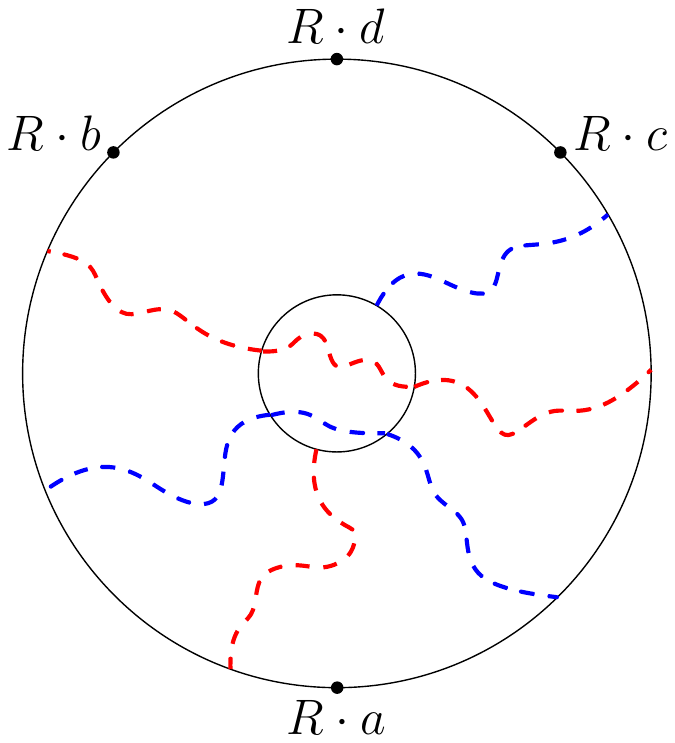}
		\caption{$\Yc_6(r,R)$.}
	\end{subfigure}
	\begin{subfigure}{.24\textwidth}
		\centering
		\includegraphics[width=0.9\linewidth]{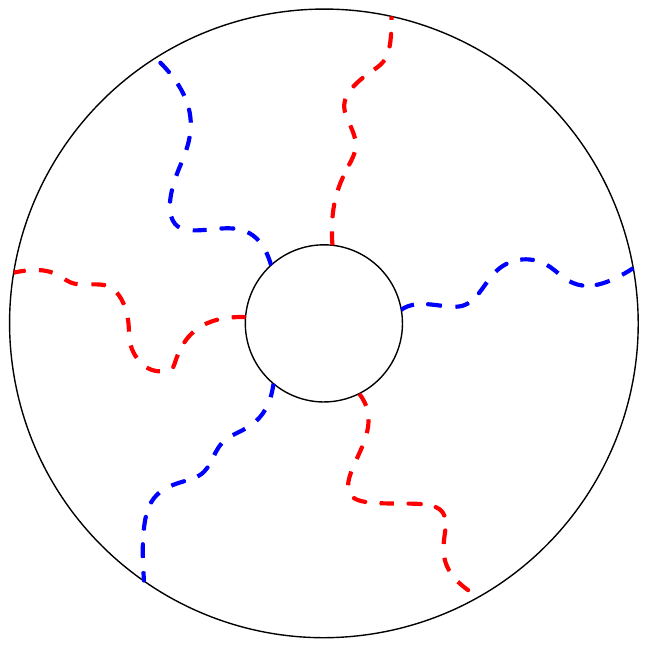}
		\caption{$\Ac_6(r,R)$.}
	\end{subfigure}
	\begin{subfigure}{.24\textwidth}
		\centering
		\includegraphics[width=0.9\linewidth]{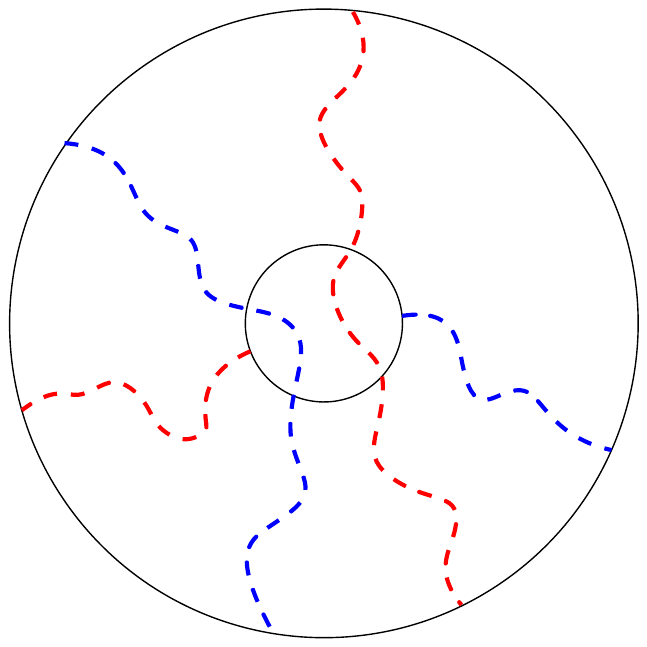}
		\caption{$\Zc_6(r,R)$.}
	\end{subfigure}
	\caption{Illustration of the events $\Xc_j(r,R)$, $\Yc_j(r,R)$, $\Ac_j(r,R)$ and $\Zc_j(r,R)$ when $j=6$. The annuli should be viewed as hexagonal domains on $\Tb^*$ under the convention given at the end of Section~\ref{subsec:setting}. The dashed red and blue curves represent arms of the respective colors.}
	\label{fig:XY}
\end{figure}

We now introduce the (polychromatic) arm events in the plane. For $j \geq 2$, the classical $j$-arm event from scale $r$ to $R$ in the plane is given by 
\begin{equation}
	\label{def:p}
	\Pc_j(r,R) = \{ \mbox{There are } j \mbox{ disjoint arms from } C_r \mbox{ to } C_R \mbox{ in } A(r,R), \mbox{ not all of the same color} \} \,.     
\end{equation}
However, in this work we are going to mainly work with the following variants. We refer readers to Figure~\ref{fig:XY} for an illustration. Fix four points $a=(0,-1)$, $b=(-1/2,\sqrt{3}/2),c=(1/2,\sqrt{3}/2)$ and $d=(0,1)$ on $C_1$. Let
$$
l(j)=\lfloor j/4\rfloor+\lfloor(j+1)/4\rfloor+1\,; \quad r(j)=\lfloor(j+2)/4\rfloor+\lfloor(j+3)/4\rfloor-1\,.
$$
(Note that $l(j)+r(j)=j$.) We call $j$ disjoint arms from $C_r$ to $C_R$ in $A(r,R)$ are ``with the prescribed pattern'' if there are $l(j)$ of them from $C_r$ to $R\cdot\wideparen{ba}$ (the counterclockwise arc from $R\cdot b$ to $R\cdot a$) with color pattern red, blue, red \ldots (counted clockwise from the point $R\cdot a$); and the rest $r(j)$ ones from $C_r$ to $R\cdot \wideparen{ac}$ with color pattern blue, red, blue \ldots (counted  counterclockwise from $R\cdot a$). Then, we define the following variants of planar $j$-arm events $\Xc_j, \Yc_j$ as
\begin{equation}
	\label{def:a}
	\Xc_j(r,R) = \{ \mbox{There exist } j \mbox{ disjoint arms from } C_r \mbox{ to }C_R \mbox{ in } A(r,R) \mbox{ with the prescribed pattern} \} \,,
\end{equation}
\begin{equation}
	\label{def:y}
	\begin{split}
		\Yc_j(r,R) &= \{ \mbox{There exist } j \mbox{ disjoint arms from } C_r \mbox{ to }C_R \mbox{ in } A(r,R) \mbox{ with the prescribed pattern} \\
		& \qquad  \mbox{and furthermore the left }(k+1) \mbox{-th arm and the right } k \mbox{-th arm are connected }\\
		& \qquad \quad \mbox{ in } B_R \text{ by a path of their (common) color for all } 1 \leq k \leq (l(j)-1) \wedge r(j) \}.
	\end{split}
\end{equation}
We will also consider two other variants $\Ac_j$ and $\Zc_j$, 
which correspond to the events $\Xc_j$ and $\Yc_j$ respectively but with no constraints on the endpoints on $C_R$. When $j$ is even, $\Ac_j$ is the arm event with alternating color sequences, referred to in the literature as the {\bf alternating arm events}. Note that $\Ac_j$ is different from $\Pc_j$ except for $j=2$, although they are comparable; See Claim (2) of Lemma \ref{prior_1}.

\begin{figure}
	
\end{figure}

\medskip

We now turn to arm probabilities. We write 
\begin{equation}
	\label{def:sequence}
b_j(r,R)= \Pb[ \Bc_j(r,R)]\,,  \mbox{ and similarly  $h_j$, $p_j$, $x_j$, $y_j$ and $a_j$ for arm probabilities}.
\end{equation}
We now collect several quick facts and review some classical estimates on the asymptotics of arm probabilities. For simplicity, we will fix $r\ge r_h(j)$ and write $h_j(n)$ as a shorthand for $h_j(r,n)$ (and resp.\ for other arm probabilities) throughout this subsection.

We now collect several classical facts on the asymptotics of arm probabilities. In the following lemma, we show that in both half-plane and plane cases, the probabilities of variants arm events are up-to-constants equal to each other.
\begin{lemma}
	\label{prior_1}
The half-plane arm probabilities are comparable for $j\ge 1$ and so is the case for planar arm events and $j\ge 2$. Namely
 \begin{itemize}
 \item[(1).] $h_j(n) \asymp b_j(n)$;
 \item[(2).] $p_j(n) \asymp x_j(n) \asymp a_j(n)$ and for $r>r_0(j)$ such that $y_j(n)>0$ for all large $n$, we further have $a_j(n) 
 	\asymp y_j(n)$.
 \end{itemize}
Moreover, with the arm exponents defined in \eqref{eq:betadef} and \eqref{eq:alphadef} resp.,
 \begin{itemize}
 \item[(3).] $b_j(R_1,R_2) = (R_1/R_2)^{-\beta_j+o(1)}$ as $R_1\ge r_b(j)$ and $R_1/R_2\to 0$. Also for $j \geq 2$, $p_j(R_1,R_2) = (R_1/R_2)^{-\beta_j+o(1)}$ as $R_1 \ge r_p(j)$ and $R_1/R_2 \rightarrow 0$. The same holds for $h_j,a_j,x_j,y_j$.
\end{itemize}
Finally, the macroscopic arm probabilities have the following asymptotics:
\begin{itemize}
 \item[(4).] There exist $f_j,g_j:(0,1)\to\mathbb{R^+}$ such that  for $0<\epsilon<1$,
$$
\lim_{n\to\infty} x_j(\epsilon n,n)=f_j(\epsilon)\quad \mbox{ and }\quad \lim_{n\to\infty} a_j(\epsilon n,n)=g_j(\epsilon).
$$
 \end{itemize}
\end{lemma}
The proofs for Claims (1) and (2) are classical applications of separation lemmas (which we will introduce shortly in this section) and will be postponed to Appendix~\ref{sec:proof-prior}. Claim (3) is the main result in \cite{MR1879816}, while Claim (4) for $a_j$ is Lemma 2.9 of \cite{MR3073882} while the case for $x_j$ follows from essentially same arguments.

\begin{remark}
\label{rmk:prior}
The power of (3) lies in the fact that $R_1$ can also go to infinity with $R_2$ as long as $R_1/R_2$ goes to zero. In the context, we will use two particular forms frequently. One is that for fixed $\alpha \in (0,1)$
\begin{equation}\label{eq:rmk_prior}
 h_j(n^\alpha ,n) \asymp b_j(n^\alpha ,n)=n^{-\beta_j(1-\alpha)+o(1)}\;\mbox{ and }\; x_j(n^\alpha ,n) \asymp y_j(n^\alpha ,n) \asymp p_j(n^\alpha ,n)=n^{-\alpha_j(1-\alpha)+o(1)};
\end{equation}
The other is that there exists $c>0$ such that when $R_1/R_2 \rightarrow \infty$
$$b_3(R_1,R_2) \leq (R_1/R_2)^{1+c} \quad \mbox{and} \quad p_6(R_1,R_2) \leq (R_1/R_2)^{2+c}\,,$$  because $\beta_3 = 2>1$ and $\alpha_6 = 35/12>2$.
\end{remark}
In the following lemma, we show that the arm probabilities are ``almost'' monotone. We also postpone its proof to Appendix~\ref{sec:proof-prior}.
\begin{lemma}
	\label{prior_3} For any $j\geq 1$
	\begin{equation*}
		\lim_{\epsilon \rightarrow 0} \liminf_{n \rightarrow \infty} \inf_{n \leq s \leq t \leq (1+\epsilon)n} \frac{h_j(t)}{h_j(s)} \geq 1\,.
	\end{equation*}
The same inequalities also hold for $b_j$ and $x_j$, $y_j$, $p_j$ and $a_j$ (for $j\geq 2$).
\end{lemma}

The following lemma is a ``functional equation'' type result on sequences and it shows that we can obtain Theorem~\ref{thm:newes-half} from Proposition~\ref{prop:ratioestimate-half} and Theorem~\ref{thm:newes-whole} from Proposition~\ref{prop:ratioestimate-whole} together with the preliminary estimates Lemma~\ref{prior_1} and \ref{prior_3} (which imply that the sequence of various types of arm probabilities satisfy Assumption (2) in the following lemma). Its proof can be found in Appendix \ref{sec:B}.
\begin{lemma}
\label{lem:sequence}
Consider a set of positive real numbers $\{ a_n: n\in \mathbb{R}\cap[1,\infty)\}$ such that:
\begin{enumerate}
\item[(1).] $\frac{a_{m^2 n}}{a_{m n}} =\frac{a_{m n}}{a_n}\Big(1+O(n^{-c})\Big)$ for all $m \in (1.1,10)$, where the constants in $O(n^{-c})$ is independent of $m$ in the given range.
\item[(2).] $\lim_{\epsilon \rightarrow 0} \liminf_{n \rightarrow \infty} \inf_{n \leq s \leq t \leq (1+\epsilon)n} \frac{a_t}{a_s}\geq 1$.
\end{enumerate}
Then, there exist $0<C<\infty$ and $-\infty<\alpha<\infty$ such that $a_n =Cn^{\alpha}\Big(1+O(n^{-c})\Big)$. 
\end{lemma}

\subsection{Faces}\label{subsec:faces}
When working on couplings of arm events, a crucial concept we will need throughout the arguments is {\bf faces}, which forms a special example of stopping sets. The notion of stopping set is in some sense a two-dimensional version of the stopping time. 

\begin{figure}
	\centering
	\includegraphics[height=.2\textwidth]{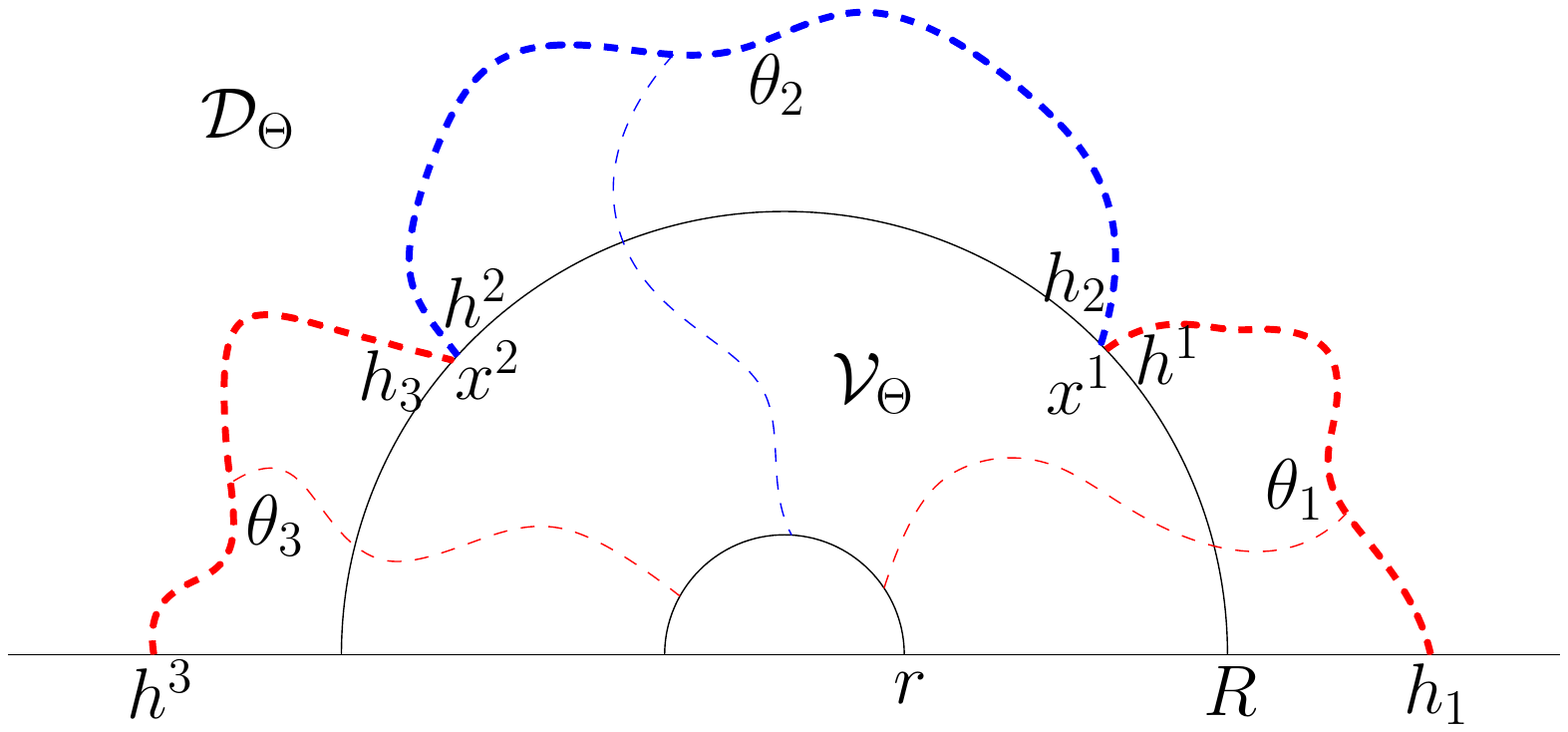}\quad
	\includegraphics[height=.2\textwidth]{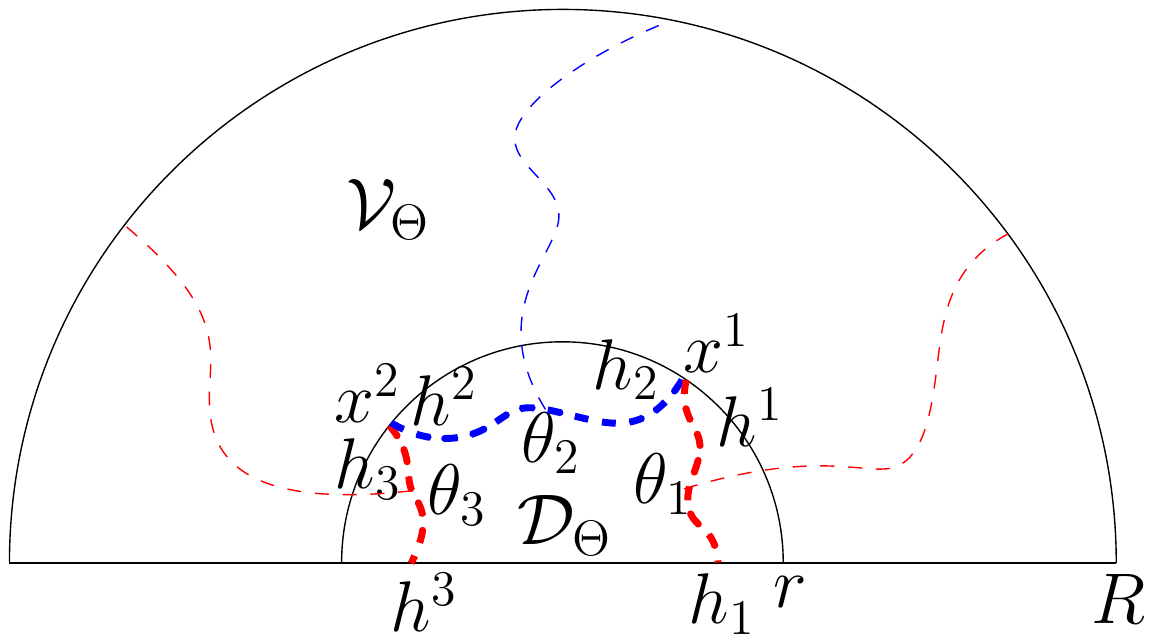}
	\caption{{\bf Left:} A configuration of outer faces around $C_R^+$. {\bf Right:} A configuration of inner faces around $C_r^+$. In both cases, $j=3$ and the configuration of faces $\Theta=\{\theta_1,\theta_2,\theta_3\}$ is depicted in heavy dashed curves with corresponding colors. We have indicated the relative position of $\Dc_{\Theta}$ and $\Vc_{\Theta}$. The events $\Bc_j^{\Theta}(r,R)$ occur when the arms indicated by normal dashed curves with colors exist. Moreover, $h_i$ and $h^i$ are the end-hexagons of $\theta_i$, and $x^1$ and $x^2$ are the endpoints $\Theta$.}
	\label{fig:hface}
\end{figure}

We begin with faces in the half-plane. Recall the conventions given at the end of Section~\ref{subsec:setting} for discretization. We call a set of paths $\Theta=\{ \theta_1,\cdots,\theta_j \}$ a configuration of \textbf{outer faces} around $C_r^+$ if the following conditions are satisfied:
\begin{itemize}
	\item For $1\le i\le j$, if $i$ is odd (resp.\ even), then $\theta_i$ is a red (resp.\ blue) path from $h_i$ to $h^i$ such that $\theta_i \subset \Hb \setminus B_r^+$.
	\item For $1\le i\le j$, the end-hexagons $h_i$'s and $h^i$'s satisfy the condition that $h_1, h^j$ touch $\partial\Hb_{\eta}$, and $h^1,h_2\dots,h^{j-1},h_j$ touch $C_r^+$. 
	Furthermore, $(h_1,h^1,\cdots,h_j,h^j)$ are listed in counterclockwise order, and $h^i$ is adjacent to $h_{i+1}$ for $1\le i\le j-1$. 
\end{itemize} 
See Figure~\ref{fig:hface} for an illustration. The paths $\theta_1,\dots,\theta_j$ are called the \textbf{outer faces} of $\Theta$. 
	The set of points $x^i:=h^i\cap h_{i+1}\cap C_r^+$ for all $1\le i\le j-1$ are called the \textbf{endpoints} of $\Theta$.
	The concepts of configurations of \textbf{inner faces} around $C_r$ are defined similarly except that we require instead that $\theta_i\subset B_r^+$ for all $i$. In both cases, each path will be called a face. We will also use $\Theta$ to denote the union of all hexagons in $\theta_i,1\le i\le j$ with a slight abuse of notation. For $\Theta$ a configurations of outer faces, we write ${\bf \Dc_{\Theta}}$ and ${\bf \Vc_{\Theta}}$ for the connected components of $\Hb_{\eta}\backslash\Theta$ whose boundary contains  $\infty$ and $0$ respectively\footnote{In these definitions, the initials D and V stand for ``discovered'' and ``vacant'', signifying the status of these regions when we explore inwards; same for inner faces. See also the definition for stopping sets.}. When $\Theta$ represents a configuration of inner faces, we exchange the role of $\infty$ and $0$ in the definition accordingly. See Figure \ref{fig:hface} for an illustration.

Suppose $r<R$. If $\Theta$ is a configuration of outer faces around $C_R^+$, then we write $\mathbf{\Bc_j^{\Theta}(r,R)}$ (resp. $\mathbf{\Hc_j^{\Theta}(r,R)}$) for the event that each outer face is connected to $C_r^+$ (resp. $[-r,r]$) by an arm in $\Vc_{\Theta}$ of the same color. The same definition applies to a configuration of inner faces $\Theta$ around $C_r^+$, with each face connected to $C_R^+$ by an arm also in $\Vc_{\Theta}$.  
When $\Theta$ is specified, there would be no confusion in whether $\mathbf{\Bc_j^{\Theta}(r,R)}$ refers to a configuration of inner or outer faces.

\medskip

In a similar fashion, in the plane, we define the configuration of outer or inner faces around a circle as a circular chain of $j$ (an even integer) paths of alternating colors. More precisely, we call the set of paths
	$\Theta=\{\theta_1, \cdots,\theta_j\}$ a configuration of outer (resp.\ inner) faces around $C_r$ if
	\begin{itemize}
		\item $\theta_1,\cdots,\theta_j$ have alternating colors and they are contained in $\Cb\setminus B_r$ (resp.\ $B_r$).
		\item The end-hexagons of these $j$ paths, $(h_1,h^1,\cdots,h_j,h^j)$, touching $C_r$, are listed in counterclockwise order (still, $h_i$ and $h^i$ are the two ends of $\theta_i$).
		Moreover, $h^i$ is adjacent to $h_{i+1}$ for all $1\le i\le j$ (where we set $h_{j+1}:=h_1$).
	\end{itemize}
\begin{figure}
	\centering
	\includegraphics[height=.3\textwidth]{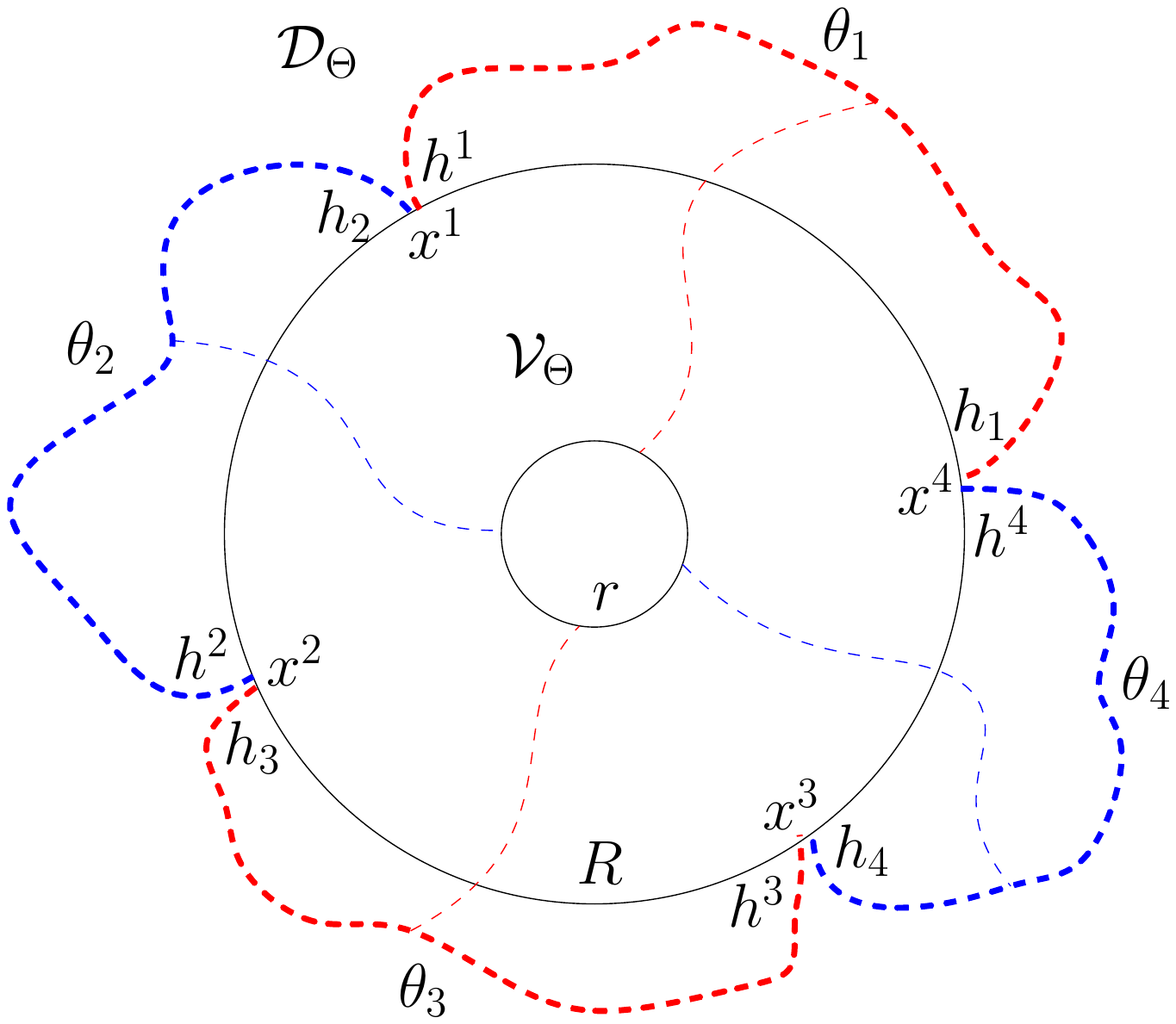}\quad\quad\quad\quad
	\includegraphics[height=.3\textwidth]{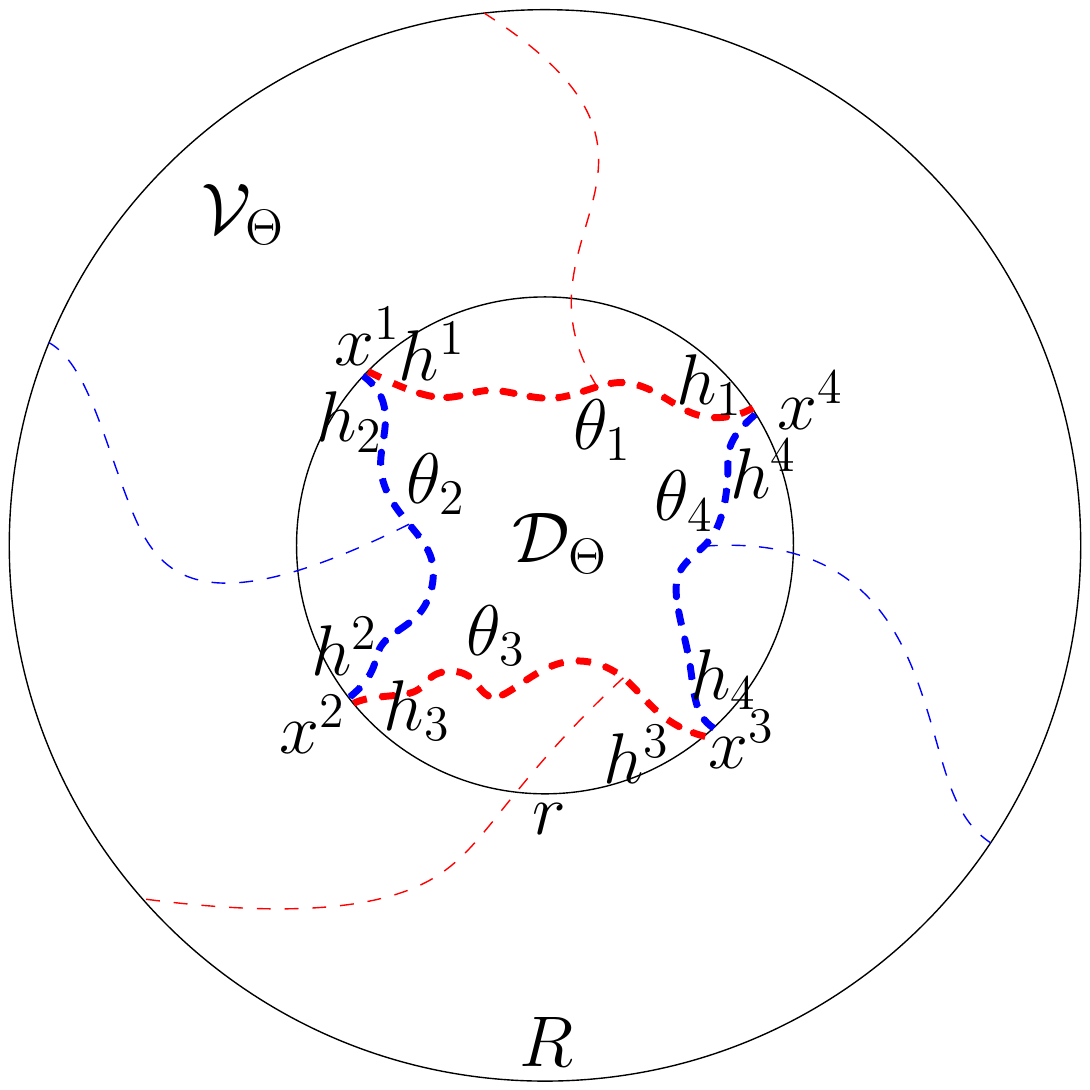}
	\caption{{\bf Left}: A configuration of outer faces around $C_R$. {\bf Right}: A configuration of inner faces around $C_r$. In both cases, $j=4$ and the configuration of faces $\Theta=\{\theta_1,\ldots,\theta_4\}$ is depicted in heavy dashed curves with corresponding colors. We have indicated the relative position of $\Dc_{\Theta}$ and $\Vc_{\Theta}$. The events $\Ac_j^{\Theta}(r,R)$ occur when the arms indicated by normal dashed curves with colors exist. Moreover, $h_i$ and $h^i$ are the end-hexagons of $\theta_i$, and $x^1,\ldots,x^4$ are the endpoints of $\Theta$.}
	\label{fig:face}
\end{figure}
See Figure~\ref{fig:face} for an illustration.
	Again, the paths $\theta_i$'s are called the outer (resp.\ inner) faces of $\Theta$, and we can define the end-points $x^1,\cdots,x^j$ analogously (by setting $x^i:=h^i\cap h_{i+1}\cap C_r$). Similar to the half-plane case, for any configuration of outer faces $\Theta$, let $\Dc_{\Theta}$ and $\Vc_{\Theta}$ denote the connected components of $\Cb\setminus\Theta$ which contains $\infty$ and $0$ respectively. If $\Theta$ stand for a configuration of inner faces instead, we then exchange the role of $0$ and $\infty$ in the definition accordingly. For $0<r<R$, if $\Theta$ is a configuration of outer faces around $C_R$, we write $\mathbf{\Ac_j^{\Theta}(r,R)}$ for the event that each outer face is connected to $C_r$ by an arm of the corresponding color in $\Dc_{\Theta}$. The corresponding events for inner faces are defined in a similar way.

By a {\bf stopping set}, we mean a random configuration of inner or outer faces $\Theta$ which can be determined by an ``one-sided'' exploration process. That is to say, there exists an exploration process such that when the exploration stops, we are able to determine $\Theta$ from the configuration in $\Dc_\Theta$ only, and all the hexagons in $\Vc_{\Theta}$ are left unexplored.

\subsection{Separation lemmas in the half-plane}\label{subsec:seplemmahalf}
In this subsection, we give different versions of the separation lemma in the half-plane. 
Their counterparts in the plane will be given in the next subsection separately.

The separation lemma for planar percolation first appears in Kesten’s work \cite{MR879034}. This is the key tool that allows one to deal with conditional arm probabilities and obtain quasi-multiplicativity for arm events. It is worth mentioning that its counterpart in the setup of Brownian motion established by Lawler in \cite{MR1386292} also has multiple implications in the study of special points of Brownian motion. 

We will state the separation lemma in terms of interfaces. By an \textbf{interface} we mean a b-path that separates clusters of different colors in some domain. 

We now introduce the notion of ``quality'' to measure how separate the interfaces are.	Suppose $u<v$.
	Let $\Gamma$ be a set of interfaces crossing $A^+(u,v)$ from $C^+_u$ to $C^+_v$. Suppose that $\Gamma$ contains $j\ge 1$ interfaces and let $\{x^1,\cdots,x^j\}$ be the collection of endpoints of these interfaces on $C^+_v$ listed in counterclockwise order. Then $\Gamma$ has an \textbf{exterior quality} defined as 
	\begin{equation}\label{eq:ex-Q}
		Q_{\rm ex}(\Gamma):=\frac{1}{v}\operatorname{d}(v,x^1)\wedge \operatorname{d}(x^1,x^2)\wedge\cdots\wedge\operatorname{d}(x^j,-v).
	\end{equation}
	In an analogous way, we define the \textbf{interior quality} for interfaces and denote it by $Q_{\rm in}$. 	We can also define the quality\footnote{There is no need to add subscripts since in this case the type of quality can be read from the type of faces.} $Q$ for a configuration of faces $\Theta$ in a similar fashion. We call $\Theta$ to be \textbf{well-separated} if $Q(\Theta)>j^{-1}$.
	
	The separation lemma roughly shows that conditioned on the interfaces reach a long distance without intersecting each other, then their endpoints will separate at some macroscopic distance with a universal positive probability. We now give the first version of separation lemma, which involves no initial configuration of faces, and is rather standard. This version can be viewed as the half-plane counterpart of \cite[Lemma A.4]{MR2630053} and its proof follows from the same line, so we omit it.
	
	\begin{lemma}[The Standard Separation lemma]\label{lem:H-sep}
		For any $j\ge 2$, there exists $c(j)>0$ such that for any $100j\le 2r\le R$, the following hold. 
		\begin{itemize}
			\item Let $\Gamma$ be the set of interfaces crossing $A^+(r,R)$ connecting $C_r^+$ and $C_R^+$, then,
			\begin{equation}\label{eq:H-sep-B}
				\Pb[Q_{\rm in}(\Gamma)\wedge Q_{\rm ex}(\Gamma)>j^{-1}\mid \Bc_j(r,R)]>c.
			\end{equation}
			\item Let $\Gamma$ be the set of interfaces in $B_R^+$ connecting $[-r,r]$ and $C_R^+$, then 
			\begin{equation}\label{eq:H-sep-H}
				\Pb[Q_{\rm ex}(\Gamma)>j^{-1}\mid \Hc_j(r,R)]>c.
			\end{equation}
		\end{itemize}
		
	\end{lemma}
	
	The following strong separation lemma is the half-plane version of \cite[Proposition A.1]{MR3073882}. In fact, Proposition \ref{prop:sepa} is called ``strong'' because it holds for any given initial configuration of faces which might be of very bad quality, not because it implies Lemma \ref{lem:H-sep} (in fact it does NOT -- note that in this strong version we do require the initial configuration to be $j$ faces, which does not appear in the conditioning in Lemma \ref{lem:H-sep}).
	For the sake of completeness, we will prove it in Appendix \ref{sec:A}.
	
	\begin{proposition}[Strong separation lemma]\label{prop:sepa}
		For any $j\ge 2$, there exists $c(j)>0$ such that for any $100j\eta\le 2r\le R$, the following hold.
		\begin{itemize}
			\item {\rm (Outward)} For any configuration of inner faces $\Theta=\{ \theta_1,\cdots,\theta_j \}$ with endpoints $\{ x^1,\cdots, x^{j-1} \}$ around $C_{r}^+$, let $\Gamma$ be the $(j-1)$-tuple of interfaces in  $\Vc_{\Theta}$ which start from $x^1,\cdots, x^{j-1}$ respectively until they reach $C^+_R$. Then
			\begin{equation}\label{eq:ex-sep}
				\Pb[Q_{\rm ex}(\Gamma)>j^{-1} \mid \Bc_j^{\Theta}(r,R)]>c.
			\end{equation}
			\item {\rm (Inward)} In a similar fashion, the same claim also holds for any configuration of outer faces $\Theta$ with $C^+_r$, $C^+_R$ and $Q_{\rm ex}$ replaced by $C^+_R$, $C^+_r$ and $Q_{\rm in}$ respectively.
		\end{itemize}
	\end{proposition}

In this work, we will also need a slightly strengthened version of the strong separation lemma stated below. Compared with Proposition \ref{prop:sepa}, it lifts the restriction of requiring exactly $j$ faces as the initial configuration, although it still requires an a priori upper bound on the number of faces. For a color configuration $\omega_0$ on $\Tb^*$ (or only on a subset $D$ of $\Tb^*$), we denote by $\{\omega_D=\omega_0\}$ the event that the color configuration inside $D$ coincides with $\omega_0$.
	\begin{proposition}[Slightly stronger separation lemma]\label{prop:half-sep-2}
		For any $K\ge j\ge 2$, there exists $c=c(j,K)>0$ such that for any $10j\cdot 2^K\le 2^Kr\le 2^KR_1\le R_2\le R$, letting $\Gamma$ be the set of interfaces crossing $A^+(R_1,R_2)$, the following holds:
		\begin{itemize}
			\item {\rm (Outward)} For any configuration of inner faces $\Theta$ around $C_{R_1}^+$ with no more than $K$ faces, and any color configuration $\omega_0$ that coincides with $\Theta$ and satisfies $\Pb[\Bc_j(r,R)\mid \omega_{\Dc_{\Theta}}=\omega_0]>0$,
			\begin{equation}\label{eq:Q-ex-1}
				\Pb\Big[|\Gamma|=j-1\text{ and }Q_{\rm ex}(\Gamma)>j^{-1} \mid \Bc_j(r,R),\omega_{\Dc_\Theta}=\omega_0\Big]>c\,.
			\end{equation}
		
			\item 
			{\rm (Inward)} In a similar fashion, for any configuration of outer faces $\Theta$ around $C_{R_2}^+$ with no more than $K$ faces, and  any suitable color configuration $\omega_0$, 
			\begin{equation}\label{eq:Q-in-2}
				\mbox{The same inequality holds with $Q_{\rm ex}$ replaced by $Q_{\rm in}$.}
			\end{equation}
		\end{itemize}
		The claims \eqref{eq:Q-ex-1} and \eqref{eq:Q-in-2} above also hold if $\Bc_j(r,R)$ is replaced by $\Hc_j(r,R)$.
	\end{proposition}
	\begin{proof}
		We focus on \eqref{eq:Q-ex-1}  as all other cases are almost the same. Given any configuration of inner faces $\Theta$ around $C^+_{R_1}$ with $k\ (j\le k\le K)$ faces, and any configuration $\omega_0$ that coincides with $\Theta$ and satisfies $\Pb[\Bc_j(r,R)\mid \omega_{\Dc_{\Theta}}=\omega_0]>0$, denote 
$$
\Bc=\Bc_j(r,R),\;\Qc=\{\omega_{\Dc_\Theta}=\omega_0\}\mbox{ and }\Uc=\{|\Gamma|=j-1\text{ and }Q_{\text{ex}}(\Gamma)>j^{-1}\}
$$
for conciseness. Our aim is to show $\Pb[\Uc\mid \Bc\cap \Qc]>c(j,K)$ for some constant $c(j,K)>0$.
		
		We first deal with the case $R_2=R$ and we proceed by induction on $K$. The case $K=j$ reduces to \eqref{eq:ex-sep} since when $\Theta$ has exactly $j$ faces, $\Pb[\,\cdot\mid \Bc\cap \Qc]$ is identical to $\Pb[\,\cdot\mid \Bc_j^\Theta(r,R)]$ in $\Vc(\Theta)$. Now, for some $K>j$, assume the result (when $R_2=R$ ) is true for $(K-1)$ and $\Theta$ has exactly $K$ faces. Let $\Gamma_1$ be the set of interfaces starting from the $(K-1)$ end-points of $\Theta$, truncated at their first visit of $\mathbb{R}\cup C^+_{2R_1}$. Consider the event 
		$$
		\Bc_{K}^\Theta(R_1,2R_1)=\big\{ \mbox{all of these $(K-1)$ interfaces reach $C_{2R_1}^+$}\big\}.
		$$
		Assume $\Bc_{K}^\Theta(R_1,2R_1)$ fails, then once we explore the $(K-1)$ interfaces starting from the $(K-1)$ endpoints and stop exploring immediately upon reaching $\Rb\cup C_{2R_1}^+$, we will see that some of them merge together or hit the real line before reaching $C_{2R_1}^+$. Thus the exploring process induces a configuration of inner faces $\Theta'$ around $C_{2R_1}^+$ with no more than $(K-1)$ faces. Hence, the induction hypothesis gives 
		\begin{equation}\label{eq:cor-sep-Kc^c}
			\Pb\left[\Uc\mid \Bc,\Qc,\left(\Bc_K^\Theta(R_1,2R_1)\right)^c\right]>c(j,K-1)\,.
		\end{equation}  
		On the other hand, assume $\Bc_{K}^\Theta(R_1,2R_1)$ happens. Clearly $$\Bc\cap \Bc_K^\Theta(R_1,2R_1)\subset \Bc_{K}^\Theta(R_1,2R_1)\cap \Bc_j(4R_1,R),$$ and also $\Bc_j(4R_1,R)$, $\Bc_{K}^\Theta(R_1,2R_1)$ and $\Qc$ are mutually independent, so 
		\begin{equation}\label{eq:cor-sep-1}
			\Pb\Big[\Bc,\Bc_K^\Theta(R_1,2R_1)\mid \Qc\Big]\le \Pb\Big[\Bc_{K}^\Theta(R_1,2R_1)\Big]\Pb\Big[\Bc_j(4R_1,R)\Big]\,.
		\end{equation}
		Furthermore, let $\Gamma_2$ be the set of interfaces crossing $A^+(4R_1,R)$, by FKG-RSW gluing,
		\begin{equation}\label{eq:cor-sep-2}
			\Pb\Big[\Uc,\Bc,\Bc_K^\Theta(R_1,2R_1)\mid \Qc\Big]\ge c_1 \Pb\Big[\Bc_{K}^\Theta(R_1,2R_1),Q_{\text{ex}}(\Gamma_1)>K^{-1}\Big]\Pb\Big[\Bc_j(4R_1,R),Q_\text{in}(\Gamma_2)>j^{-1}\Big]\,
		\end{equation}
		for some $c_1=c_1(j,K)>0$ (here the condition $\Pb[\Bc_j^\Theta|\Qc]>0$ assures the existence of feasible connecting pattern).
		By \eqref{eq:H-sep-B} and \eqref{eq:ex-sep} (taking $j$ therein as $K$), we get
		\begin{align*}
			\Pb\big[\Bc_K^\Theta(R_1,2R_1),Q_\text{ex}(\Gamma_1)>K^{-1}\big]&\;>c_2\Pb\big[\Bc_{K}^\Theta(R_1,2R_1)\big]\mbox{ and }\\
			\Pb\big[\Bc_j(4R_1,R),Q_\text{in}(\Gamma_2)>j^{-1}\big]&\;>c_3\Pb\big[\Bc_j(4R_1,R)\big]\,
		\end{align*}
		for some positive constants $c_2(K),c_3(j)$.
		This combined with \eqref{eq:cor-sep-1} and \eqref{eq:cor-sep-2} yields
		\begin{equation}\label{eq:cor-sep-Kc}
			\Pb\left[\Uc\mid \Bc,\Qc,\Bc_{K}^\Theta(R_1,2R_1)\right]=\frac{\Pb[\Vc,\Bc,\Bc_{K}^\Theta(R_1,2R_1)\mid \Qc]}{\Pb[\Bc,\Bc_{K}^\Theta(R_1,2R_1)\mid \Qc]}>c_1c_2c_3=c_4(j,K)\,.
		\end{equation}
		Now pick $c(j,K)=c(j,K-1)\wedge c_4(j,K)$, then $\Pb[\Uc\mid \Bc,\Qc]>c(j,K)$ follows from \eqref{eq:cor-sep-Kc^c}, \eqref{eq:cor-sep-Kc} and the total probability formula. This proves \eqref{eq:Q-ex-1} for the case $R_2=R$. 
		
		For the general case of \eqref{eq:Q-ex-1}, by the quasi-multiplicativity (see \cite[(3.14)]{MR2630053}) and what we just proved, we have for some constants $c_5,c_6,c_7>0$ depending only on $j,K$,
		\begin{equation*}
			\Pb[\Uc,\Bc\mid \Qc]> c_5\Pb[\Uc,\Bc_j(r,R_2)\mid \Qc]\Pb[\Bc_j(R_2,R)]> c_6\Pb[\Bc_j(r,R_2)\mid \Qc]\Pb[\Bc_j(R_2,R)]> c_7\Pb[\Bc\mid \Qc]\,.
		\end{equation*}
		This shows $\Pb[\Uc\mid \Bc,\Qc]>c_7(j,K)$, as required.
	\end{proof}
	
	\begin{remark}
		\label{rmk:constant-not-essential}
		The constant $2$ in the condition $2r\leq R$ in the statement of Lemma~\ref{lem:H-sep} and Proposition~\ref{prop:sepa} is not essential. Indeed, similar results still hold for any constant greater than $1$ when $r$ is large enough, and the proof is just a verbal change of the proof given in Appendix~\ref{sec:A}. Then it is clear from the above proof that Proposition~\ref{prop:half-sep-2} still holds for large $r$ when the constant $2^K$ is replaced by any constant greater than $1$. 
	\end{remark}

Finally, in the half-plane, we are also able to obtain a much stronger version of the separation lemma, which first appears in the name of ``super strong separation lemma'' as Conjecture A.4 in \cite{MR3073882} in a slightly different form. This ``ultimate'' version no longer requires any a priori upper bound on the number of faces and allows us to explore the configuration annulus by annulus when we construct couplings in Section \ref{sec:coupling}, which is simpler and more natural than the exploration of faces.

	\begin{proposition}[Super strong separation lemma]\label{prop:super-strong}
		For any $j\ge 2$, there exist $M,c>0$ depending only on $j$, such that for any $10M^2j\le M^2r\le R$ and $Mr\le u\le R/M$, if we let $\Gamma$ be the set of interfaces crossing $A^+(r,R)$, then the following holds.
		\begin{itemize}
			\item {\rm (Outward)} For any percolation configuration $\omega_0$ which satisfies $\Pb\left[\Bc_j(r,R)\mid \omega_{B_u^+}=\omega_0\right]>0$, 
			\[
			\Pb\left[|\Gamma|=j-1\mbox{ and }Q_{\text{ex}}(\Gamma)>j^{-1}\mid \Bc_j(r,R),\omega_{B_u^+}=\omega_0\right]>c\,.
			\]
			\item {\rm (Inward)} In a similar fashion, for any suitable color configuration $\omega_0$, the same claim follows with $\omega_{B_u^+}$ and $Q_{\text{ex}}$ replaced by $\Hb\setminus B_u^+$ and $Q_{\text{in}}$ respectively.
		\end{itemize}
		In addition, similar results also hold if $\Bc_j(r,R)$ is replaced by $\Hc_j(r,R)$.
	\end{proposition}
	
	With this powerful tool in hand, we are able to derive the coupling results and estimates in the half-plane in a relatively straightforward (and natural) way. Although we believe that the super strong separation lemma also holds for the plane case, unfortunately we could not find a proof due to the lack of a natural ordering for interfaces in the plane. In order to state and prove all results in parallel, even in the half-plane case we choose to use only the weaker version Proposition~\ref{prop:half-sep-2}, which can be extended to the plane case (see Proposition~\ref{prop:C-sep}) without difficulty. In what follows, we provide a sketch of proof for the half-plane super strong separation lemma.

	\begin{figure}
	\begin{subfigure}{.33\textwidth}
		\centering
		\includegraphics[width=\linewidth]{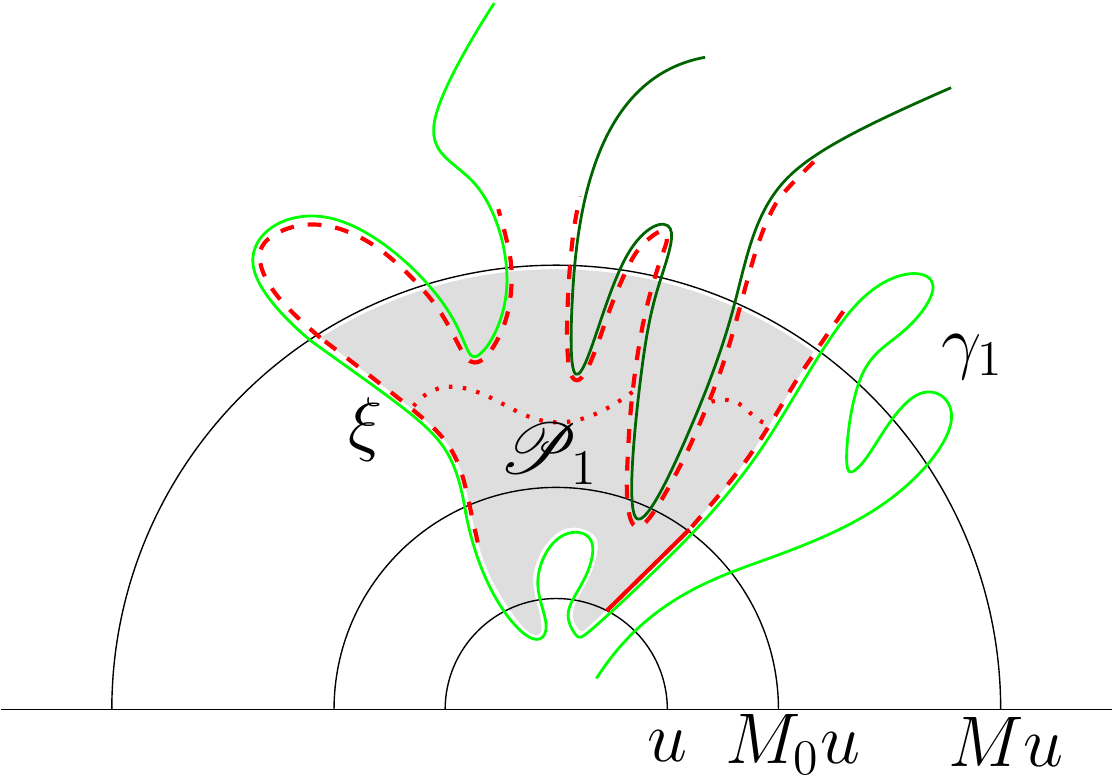}
		\caption{}
	\end{subfigure}
	\begin{subfigure}{.33\textwidth}
		\centering
		\includegraphics[width=\linewidth]{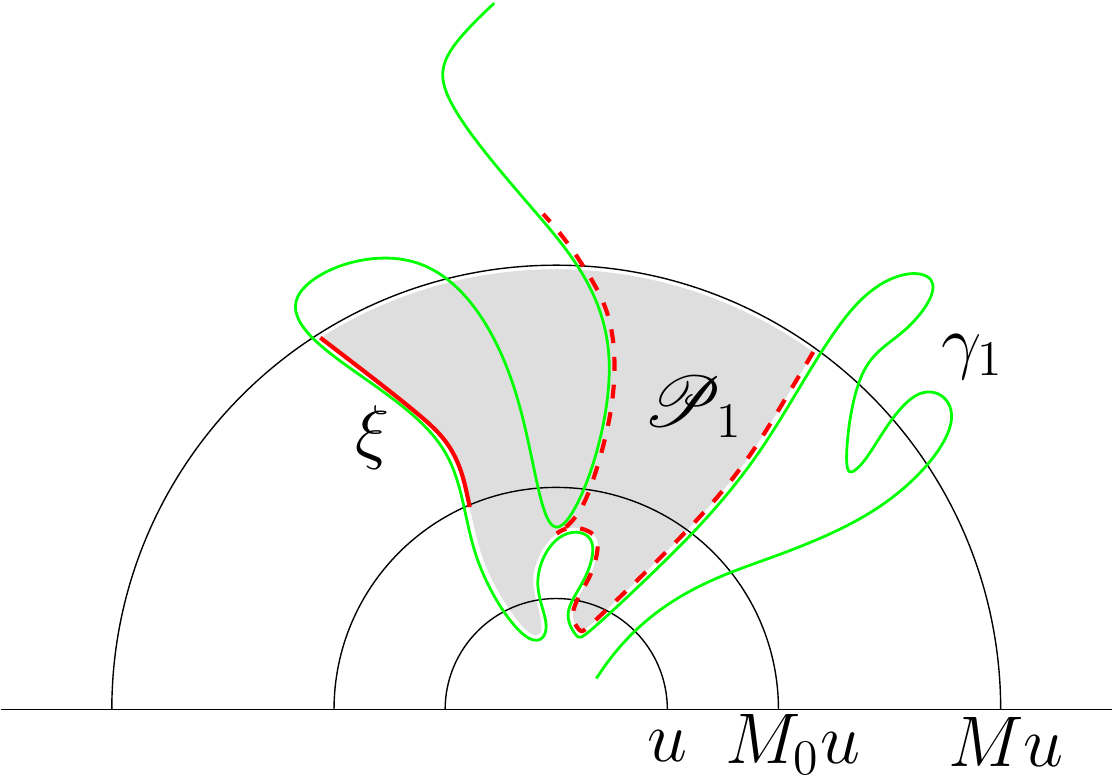}
		\caption{}
	\end{subfigure}
	\begin{subfigure}{.33\textwidth}
		\centering
		\includegraphics[width=\linewidth]{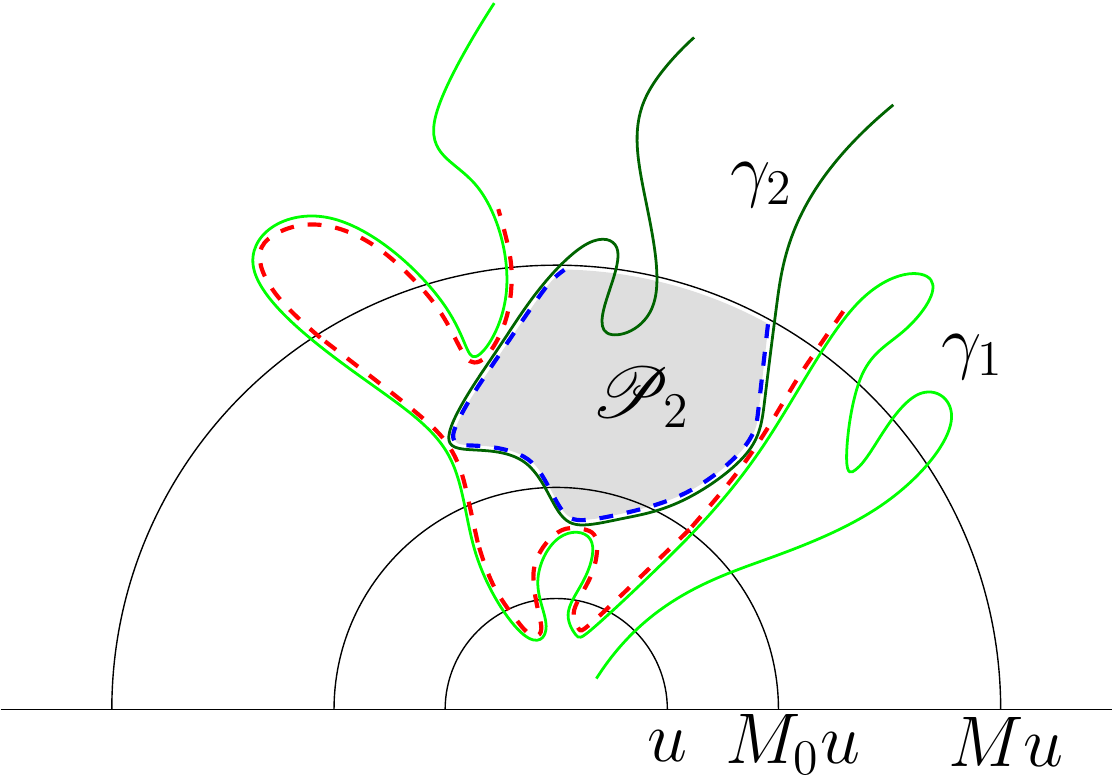}
		\caption{}
	\end{subfigure}
	\caption{Three cases for the proof of Proposition \ref{prop:super-strong}. \textbf{(a)}: No touchings for $\mathscr{P}_1$. Here and in case {(b)}, $\gamma_1$ is the first (in counterclockwise order) global interface that crosses $A^+(u,Mu)$ more than once, depicted in green, while $\xi$ is the leftmost crossing of $\gamma_1$ from $u$ to $M u$. The pocket $\mathscr{P}_1$ induced from $\gamma_1$ down to scale $u$ is given by the gray region, which is enclosed by $\xi$ and the part of $\gamma_1$ proceeding $\xi$ from the last visit on $C^+_{Mu}$. The figure illustrates the case that neither the part of $\gamma_1$ following $\xi$ nor other (global) interfaces \textbf{touch}  $\mathscr{P}_1$ (i.e., adjacent to $\mathscr{P}_1$ inside $A^+(u,M_0u)$). By RSW theory with positive probability there are red crossings (in dotted) connecting left and right sides of $\mathscr{P}_1$, on which there is an extra red arm (in solid) that does not use the hexagons ensuring $\Bc_j(r,R)$. \textbf{(b)}: The interface $\gamma_1$ touches $\mathscr{P}_1$. In this case, the red arm (in solid) neighboring $\xi$ is an extra arm which plays the same role as that in case (a). \textbf{(c)}: Another interface touches $\mathscr{P}_1$. We use $\gamma_2$ to denote this new interface which is depicted in dark green. Then $\gamma_2$ induces a pocket $\mathscr{P}_2$ inside $\mathscr{P}_1$ down to scale $M_0u$, depicted in gray (erasing $\mathscr{P}_1$). In the pocket $\mathscr{P}_2$, we can iterate the arguments as before by considering existence (or non-existence thereo) of touchings of $\mathscr{P}_2$ inside $A^+(M_0u,M_0^2u)$.}
	\label{fig:sss}
\end{figure}	
	
	\begin{proof}[Sketch of Proof for Proposition~\ref{prop:super-strong}] We only discuss the outward case. The crux is to show the following fact: there is a constant $c'(j)>0$, such that for any $Mr\le u\le R/M$ and any percolation configuration $\omega_0$ satisfying $\Pb[\Bc_j(r,R)\mid \omega_{B_u^+}=\omega_0]>0$ 
		if we let $\Gamma_u$ be the set of interfaces crossing $A^+(u,Mu/2)$, then $\Pb[|\Gamma_u|=j-1\mid \Bc_j(r,R),\omega_{B_u^+}=\omega_0]>c'$. 
		Once this is true, Proposition~\ref{prop:super-strong} follows by exploring all interfaces crossing $A^+(u,Mu/2)$ then applying the strong separation lemma (Proposition~\ref{prop:sepa}) in $A^+(Mu/2,Mu)$.

		To prove the fact, we note that the event $\{|\Gamma_u|>j-1\}$ implies that there is a ``pocket'' down to scale $u$ induced from some interface (see $\mathscr{P}_1$ in (a) of Figure~\ref{fig:sss}). Basically, the existence of such an interface pocket together with some FKG-RSW gluing technique leads to the disjoint occurrence of $\Bc_j(r,R)$ and a rare one-arm event crossing a scale $M_0$ (see (a) and (b) of Figure~\ref{fig:sss}). If this is the case, then applying BK-Reimer's inequality yields the desired result. However, the interfaces might touch one another, imposing difficulty to find an extra arm and apply BK-Reimer's inequality (see (c) of Figure~\ref{fig:sss}). To tackle this issue, we need an iterative argument to deal with the nesting cases. Thanks to the simply connecting property of the zero-removed half-plane, we have a natural ordering for all these interfaces so the nesting process is indeed finite in some sense. With this observation, the desired fact can be proved by a union bound if we set $M_0$ large enough and choose $M=2M_0^j$ such that the iterative arguments end in $(j-1)$ steps. This concludes the proof.
	\end{proof}

\subsection{Separation lemmas in the plane}\label{subsec:seplemmawhole}

This subsection is devoted to separation lemmas  in the plane. As the arguments are similar to the half-plane case, we only give the statements in Proposition \ref{prop:C-sep} and omit the proofs.

To begin with, we give the definition of the quality of interfaces, just as in the half-plane case. For $u<v$, let $\Gamma$ be a set of percolation interfaces that start from $C_u$ and end at the first visit of $C_v$. Suppose that $\Gamma$ contains $j$ interfaces with endpoints $x^1,\dots,x^j$ on $C_v$, then we define the \textbf{exterior quality} of $\Gamma$ as
	\begin{equation}
		Q_{\text{ex}}(\Gamma):=\frac{1}{v}\operatorname{d}(x^1,x^2)\wedge\operatorname{d}(x^2,x^3)\wedge\cdots\wedge\operatorname{d}(x^j,x^1).
	\end{equation}
Again, we define the \textbf{interior quality} for interfaces in the same fashion and denote it by $Q_{\rm in}$ as well as the quality for a configuration of faces $\Theta$ and denote it by $Q(\Theta)$.	If $\Theta$ satisfies $Q(\Theta)>j^{-1}$, we call it {\bf well-separated}. Recall the notation $\Vc(\Theta)$ and $\Dc(\Theta)$ from Section \ref{subsec:faces}.

There are counterparts of Lemma~\ref{lem:H-sep} and Proposition~\ref{prop:sepa} in the plane case, but things become subtler since the number of interfaces crossing an annulus must be even. Lemma A.4 and Proposition A.1 in \cite{MR3073882} only consider the case when $j$ is even, but the same argument with little change yields similar results for odd $j$. More precisely, we set $J=2\lfloor j/2\rfloor$ for any $j\ge 2$, then there exists $c(j)>0$ such that for any $10j\le2r\le R$, and any configuration of outer faces $\Theta=\{\theta_1,\dots,\theta_J\}$ around $C_R$, let $\Gamma$ be the set of interfaces crossing $A(r,R)$, then
\begin{equation}
	\Pb[Q_\text{in}(\Gamma)>J^{-1}\mid \Ac^\Theta_j(r,R)]>c.
\end{equation}
And similar result holds for any configuration of inner faces with $J$ faces. From this we can prove (by almost identical arguments) the counterpart of Proposition~\ref{prop:half-sep-2} in the plane case.
\begin{proposition}\label{prop:C-sep}
	For any $K\ge j\geq 2$, there exists $c=c(j,K)>0$ such that for any $10j\cdot 2^K\leq 2^Kr\leq 2^KR_1\leq R_2\le R$, let $\Gamma$ be the set of interfaces crossing $A(R_1,R_2)$, then the following holds. For any configuration of outer faces $\Theta$ around $C_{R_2}$ with no more than $K$ faces, any color configuration $\omega_0$ coincides with $\Theta$ and satisfies $\Pb[\Xc_j(r,R)\mid \omega_{\Dc_{\Theta}}=\omega_0]>0$,
	\begin{equation}\label{eq:SSL_1}
		\Pb[|\Gamma|=J\text{ and }Q_\text{ex}(\Gamma)>J^{-1}\mid \Xc_j(r,R),\omega_{\Dc_\Theta}=\omega_0]>c\,.
	\end{equation}
	Moreover, the same claim also holds for $\Yc_j(r,R)$.
\end{proposition}
The first claim \eqref{eq:SSL_1} can be deduced similarly as before. For the second claim, we note that although conditioning on $\Yc_j(r,R)$ involves more requirements on connectedness, it could be resolved by similar gluing technique given in the proof of Lemma 3.3 in \cite{MR3073882}. We omit the details.

\subsection{Percolation exploration process and scaling limits}\label{subsec:exploration}
In this subsection, we define the percolation exploration process and discuss its power-law convergence towards ${\rm SLE}_6$ obtained by Binder and Richards. In particular, we give in Proposition~\ref{prop:couple-SLE} a variant of their results tailored for our work.  
We will also show that the exploration process is stable under perturbations of the boundary in Proposition~\ref{prop:stable-boundary}.  

Given a Jordan domain $\Omega$ and two boundary points $a, b$, for any $0<\eta<1$, recall that $\Omega_\eta$ is a connected subgraph of $\Tb_\eta$ (where we assume that $\eta$ is small enough) and $a_\eta$ and $b_\eta$ are two vertices on the boundary of $\Omega_\eta$. In addition, $a_\eta$ and $b_\eta$ divide the boundary of $\Omega_\eta$ into two connected parts, and we will assign red and blue colors to the external vertices adjacent to these two parts, respectively, according to the context. Define the exploration process $\gamma_\eta$ as a directed path from $a_\eta$ to $b_\eta$ on $\Tb^*_\eta$ which is the interface between the two red and blue clusters containing the boundary.
 
The Schramm-Loewner evolution (SLE), first introduced by Schramm in \cite{schramm2000scaling}, is a family of conformally invariant random curves that serves as the canonical candidate for the scaling limit of curves, in particular interfaces of various critical models in 2D. In the case of critical planar percolation, Smirnov proves in \cite{smirnov2001critical} that the scaling limit of the exploration process $\gamma_\eta$ is indeed the chordal ${\rm SLE}_6$ path from $a$ to $b$ in $\Omega$, which we denote by $\gamma$ below. 

As discussed in Section \ref{subsec:comments}, Binder and Richards improve the above convergence by giving a power-law rate for curves up to some stopping time. We now describe their results in more detail.
Given $U \subset \Omega$, let $T_\eta$ be the first time that $\gamma_\eta$ enters $U_\eta$ and $T$ be the first time that $\gamma$ enters $U$. Define the distance $d$ between two paths $l_1:[0,t_1] \rightarrow \Rb^2$ and $l_2:[0,t_2] \rightarrow \Rb^2$ as
\begin{equation*}
d(l_1,l_2)=\inf_{\theta}\sup_{0 \leq s \leq t_1} \left|l_2(\theta(s))-l_1(s)\right|,
\end{equation*}
where $\theta$ runs over all increasing homeomorphisms from $[0,t_1]$ to $[0,t_2]$. We call $\Omega$ a {\bf nice domain} if $\Omega$ is a domain with a piece-wise smooth boundary. Note that the domains we consider in Sections \ref{subsec:4.1} and \ref{subsec:4.3}  (both denoted by $\Omega$) are nice domains.

In Theorem 4.1.11 of \cite{richards2021convergence}, it is proved that under some coupling of $\gamma$ and $\gamma_\eta$, they are close upon specific stopping times with high probability. We collect their results here.
\begin{theorem}[\cite{richards2021convergence}]
\label{thm:richards}
When $\Omega$ is a nice domain and $U = B(b,\epsilon) \cap \Omega$ for some $\epsilon>0$, there exists $u = u(\Omega, U)>0$ such that under some coupling of $\gamma$ and $\gamma_\eta$
\begin{equation*}
\Pb\left[d\left(\gamma_\eta|_{[0,T_\eta]}, \gamma|_{[0,T]}\right) > \eta^u \right] = O(\eta^u) \quad \mbox{as }\eta \rightarrow 0\,.
\end{equation*}
\end{theorem}

Combined with estimates on critical percolation, we can prove the following proposition which states that the time can be taken as the hitting time of a nice domain containing the endpoint in our setup.

\begin{proposition}
\label{prop:couple-SLE}
For a nice domain $\Omega$ and some nice domain $U\subset\Omega$ that contains a neighbor of $b$, there exists $u = u(\Omega)>0$ such that under some coupling of $\gamma$ and $\gamma_\eta$ (which we will later call {\bf good coupling}),
\begin{equation}
\label{eq:couple-SLE}
\Pb\left[d\left(\gamma_\eta|_{[0,T_\eta]}, \gamma|_{[0,T]}\right) > \eta^u \right]<O(\eta^u) \,.
\end{equation}
\end{proposition}
\begin{proof}
Fix $\epsilon>0$ such that $U \cap \Omega \supset B(b,\epsilon) \cap \Omega$. We write $T'_\eta$ for the first time that $\gamma_\eta$ enters $B(b_\eta,\epsilon)$ and $T'$ for the first time that $\gamma$ enters $B(b,\epsilon)$. Note that the stopping times in Theorem~\ref{thm:richards} can be chosen as $T_\eta'$ and $T'$. Hence, there exists $c>0$ such that under some coupling $\Pb$ of $T_\eta$ and $T$, 
\begin{equation}
\label{eq:couple-SLE-2}
\Pb\left[d\left(\gamma_\eta|_{[0,T'_\eta]}, \gamma|_{[0,T']}\right) >\eta^c \right]<O(\eta^c) \,.
\end{equation}
We now show that $\Pb$ is indeed a good coupling. Since $U_{[\eta]} \cap \Omega_{[\eta]} \supset B(b_\eta,\epsilon) \cap \Omega_{[\eta]}$ and $U \cap \Omega \supset B(b,\epsilon) \cap \Omega$, we have $T_\eta \leq T_\eta'$ and $T \leq T'$. Now, it is sufficient to prove that $\gamma_\eta$ and $\gamma$ respectively hit $U_\eta$ and $U$ almost at the same time and place with high probability. Write
\begin{align*}
\Fc=\big\{&\mbox{there does not exist $x$ in $\partial U_\eta$ such that $\gamma_\eta$ enters $B(x, \eta^c)$ but does not hit $U_\eta$ }\\ &\qquad \mbox{between the last entering and first leaving times of $B(x, \eta^{c'})$, and the same for $\gamma$}\big\}\,,
\end{align*}
here $c'$ is a constant to be chosen and $c>c'>0$. If $\Fc^c$ happens, we have a half-plane $3$-arm event from a $\eta^c$-ball on $\partial U_{[\eta]}$ to distance $\eta^{c'}$\footnote{Careful readers will find that the boundary is not absolutely straight so these are not half-plane $3$-arm events defined before and hence we cannot directly apply estimates from Remark~\ref{rmk:prior}. However, it is relatively easy to show that the probability that there exist polychromatic arms crossing a $(\pi + 0.01)$-cone of radii $R_1$ and $R_2$ is also smaller than $(R_1/R_2)^{1+c}$. Another issue is that the boundary here can also be tilted, but similar arguments hold for a tilted $(\pi+0.01)$-cone as well thanks to the conformal invariance of the scaling limit. We will also meet similar issues later, but they can be handled in the same way.}, or, at a non-smooth point or a point $\eta^c$ close to $\partial \Omega_{[\eta]}$, a $1$-arm event in the half-plane squared. So, $\Pb[\Fc^c] \leq O(\eta^{-c}) \times O(\eta^{(1+\delta)(c-c')})+O(1) \times O(\eta^{\delta (c-c') })$. Thus, we can choose $c'$ small such that $\Pb[\Fc^c] \leq O(\eta^{c''})$. Next, we show that 
$$\Fc \cap \left\{ d\left(\gamma_\eta|_{[0,T'_\eta]}, \gamma|_{[0,T']}\right) \leq \eta^c \right\} \subset \left\{ d\left(\gamma_\eta|_{[0,T_\eta]}, \gamma|_{[0,T]}\right) \leq \eta^{c'}\right\}.$$ Assume that $d\left(\gamma_\eta|_{[0,T'_n]}, \gamma|_{[0,T']}\right) \leq \eta^c$. When $\gamma_\eta$ hits $U_{[\eta]}$, the process $\gamma$ also gets $\eta^c$ close to $U$. Then, $\gamma$ will hit $U$ before traveling a distance of $\eta^{c'}$, otherwise there will be a point $x$ on $\partial U$ such that $\gamma$ enters $B(x,\eta^c)$ but does not hit $U$ before leaving $B(x,\eta^{c'})$ which contradicts $\Fc$. Similarly, when $\gamma$ hits $U$, the process $\gamma_\eta$ will hit $U_{[\eta]}$ before leaving distance $\eta^{c'}$. So, $\gamma_\eta$ and $\gamma$ hit respectively $U_{[\eta]}$ and $U$ almost at the same time and place, and their $d$-distance is smaller than $\eta^{c'}$. This completes the proof of the above relationship. Therefore, under this coupling
$$
\Pb\left[d\left(\gamma_\eta|_{[0,T_\eta]}, \gamma|_{[0,T]}\right) >\eta^{c'}\right] \leq \Pb\left[d\left(\gamma_\eta|_{[0,T'_n]}, \gamma|_{[0,T']}\right) >\eta^c \right] +\Pb[\Fc^c] \leq O(\eta^c )+O(\eta^{c''})\,.
$$
Letting $u= \min \{c,c',c'' \}$ yields \eqref{eq:couple-SLE} as desired.
\end{proof}
We will also use a variant of this proposition, i.e., Proposition~\ref{prop:couple-SLE-half}, in the proof of Proposition~\ref{prop:ratioestimate-half}.

The next proposition states that the exploration process is insensitive to the change of boundaries. More precisely, when the boundary changes by $\eta^\delta$, then the exploration process changes by at most $\eta^{c(\delta)}$ with probability at least $(1-\eta^{c(\delta)})$.   
\begin{proposition}
\label{prop:stable-boundary}
Suppose $\Omega$ is a nice domain and $\Omega_{[\eta]}$ be the discretization of $\Omega$ defined in Section \ref{subsec:setting}. Let $\Omega'_{\eta}$ be a simply connected sub-graph of $\Tb_{[\eta]}^*$ such that if $u,v$ are two vertices of $\Omega'_{\eta}$ then the edge $uv$ also belongs to the edge set of $\Omega_{[\eta]}$. Let $a_\eta',b_\eta'$ be two boundary points of $\Omega'_{\eta}$. If for some $\delta>0$,
$$
d(\partial \Omega_{[\eta]}, \partial \Omega'_{\eta}) <\eta^\delta, \quad d(a_\eta,a_\eta')<\eta^\delta, \quad d(b_\eta, b_\eta')<\eta^\delta\,.
$$
then there exists a coupling $\Pb$ of $\gamma_\eta$ and $\gamma_\eta'$ and a constant $c(\delta)>0$ such that
$$
\Pb\left[d\left(\gamma_\eta, \gamma_\eta' \right) > \eta^{c} \right]<O(\eta^{c}).
$$
\end{proposition}
\begin{proof}
Without loss of generality, we assume that $\Omega'_{\eta} \supset \Omega_{[\eta]}$ (otherwise we can compare $\Omega_{[\eta]}$ and $\Omega'_{\eta}$ both with $\Omega_{\eta}''$, where $\Omega_{\eta}'' = \big\{  x \in \Omega_{[\eta]}: d( x , \partial \Omega_{[\eta]}) > \eta^\delta\big\}$. We can couple the critical percolation on $\Omega_{[\eta]}$ and $\Omega_{\eta}'$ such that the two configurations have the same color on $\Omega_{[\eta]}$. Under this coupling, the exploration processes in these two domains can only be different after one of them hits $\partial \Omega_{[\eta]}$ (since they have the same path in the interior of $\Omega_{[\eta]}$). Furthermore, we can find some $c \in (0,\delta)$ such that with more than $(1-\eta^c)$ probability, each time they hit $\Omega_{[\eta]}$, they will remain the same after $\eta^c$ distance (otherwise we have a half plane $3$-arm event from an $\eta^\delta$ ball on the boundary to distance $\eta^c$, or, not a non-smooth point, a $1$-arm event in the whole plane. This happens with probability less than $\eta^c$ when $c$ is close enough to $0$). This completes the proof.
\end{proof}
\begin{remark}\label{rmk:insensitivity} With similar arguments, we can show that the sharp asymptotics we obtain in this work for arm events under the specific discretization, also hold for other discretizations, as they differ by at most a power of the mesh size. In other words, arm events are insensitive to the way we discretize domains. In particular, the good coupling we construct in Proposition \ref{prop:couple-SLE} (as well as Proposition~\ref{prop:couple-SLE-half}) also works for the discretized domains under the convention given at the end of Section \ref{subsec:setting}.
\end{remark}
	 \section{Couplings and conditional arm probabilities}\label{sec:coupling}
In this section, we will give various couplings concerning arm events in the half-plane and the plane; see Propsitions \ref{prop:one-coupling}, \ref{prop:coup-1} and \ref{prop:C-in-coupling-even} and deduce useful estimates on the conditional arm probabilities from these couplings; see Propositions \ref{prop:couple-h-half}, \ref{prop:couple-bh-half}, \ref{prop:couple-y-whole}, \ref{prop:couple-a-whole}. In \cite{MR3073882}, the authors have established the coupling results concerning one-arm (\cite[Proposition 5.2]{MR3073882}) and four-arm  (\cite[Proposition 3.1]{MR3073882}) events in the plane, respectively. Although many techniques developed in \cite{MR3073882} can be adapted for variants of arm events we are considering here (e.g., $\Hc_j(r,R)$, $\Xc_j(r,R)$ and $\Yc_j(r,R)$), certain intricacy arises from our choice of such events, causing extra difficulties. We postpone the proofs to Section \ref{sec:proofcoupling} as they are rather technical and independent of the main story. 
  
\subsection{The half-plane case}
In this subsection, we deal the coupling results in the half-plane. For simplicity of presentation, the $j=1$ case and the $j\ge 2$ cases are stated separately.

\subsubsection{One-arm}
In this subsection, we will state the coupling result for one-arm event in the half-plane. 
We say that a b-path in the half-plane is a \textbf{circuit} if it connects the discretized negative real line to the discretized positive real line. 
If $\eta<r<u$, let $\Gamma_{\operatorname{out}}(r,u)$ (resp. $\Gamma_{\operatorname{in}}(r,u))$ be the outermost (resp. innermost) red circuit in the semi-annulus $A^+(r,u)$, i.e., the red circuit which is closest to $C_u^+$ (resp. $C_r^+$). 
If such circuits do not exist, then set $\Gamma(r,u)=\emptyset$.  
\begin{proposition}\label{prop:one-coupling}
	There exists $\delta>0$ such that for all $100<10r<R$ and $m\in (1.1,10)$, denoting $u:=\sqrt{rR}$, then the following hold:
	\begin{itemize}
		\item {\rm (Inward coupling)} There is a coupling $\Qb$ of the conditional laws $\Pb[\cdot\mid \Hc_1(r,R)]$ and $\Pb[\cdot\mid \Hc_1(r,mR)]$ such that if we sample $(\omega_1,\omega_2)\sim\Qb$, then with probability at least $\big(1-(r/R)^{\delta}\big)$,
		$\Gamma_{\operatorname{out}}(r,u)$ for both $\omega_1$ and $\omega_2$ are non-empty and identical, Furthermore, the percolation configurations in the connected component of $\Hb\setminus\Gamma_{\operatorname{out}}(r,u)$ whose boundary contains $0$ are also identical.
		\item {\rm (Outward coupling)} A similar coupling exists for the conditional laws $\Pb[\cdot\mid \Hc_1(r,R)]$ and $\Pb[\cdot\mid \Bc_1(r,R)]$ with $\Gamma_{\operatorname{out}}$ and $0$ replaced by $\Gamma_{\operatorname{in}}$ and $\infty$ respectively. 
	\end{itemize}
\end{proposition}

\begin{remark}\label{rmk:aboutu}
In the statement above, we take $u=\sqrt{rR}$ merely for convenience. In fact, it can be taken as $r^dR^{1-d}$ for any $0<d<1$ (with the corresponding $\delta(d)$). The same applies to all couplings in this section.
\end{remark}
The proof of the above proposition is essentially the same as that of \cite[Proposition 5.2]{MR3073882}, which deals with one-arm events in the plane in which the key tools are the RSW theory and FKG inequality. Therefore, we omit the proof and refer the readers to the said paper for details.

\subsubsection{$j$-arms with $j\ge2$}\label{subsec:H-j}
In this subsection, we concentrate on the $j$-arm events with $j\ge 2$ in the half-plane. We shall present couplings in both the inward and outward directions simultaneously for they share the same virtue. Note that for $j\ge2$, $\Hc_j(r,R)$ is equivalent to the event that there are at least $(j-1)$ interfaces starting from $C_R^+$ and end at their first hitting of $[-r,r]$; similarly, $\Bc_j(r,R)$ is equivalent to the event that there are at least $(j-1)$ interfaces crossing the annulus $A^+(r,R)$. 
In both cases, we write $\Gamma=\{\gamma_1,\dots,\gamma_n\}$ (where $n\ge j-1$) for the set of interfaces crossing the whole region, and order them counterclockwise. We say two interfaces in $\Gamma$ are \textbf{adjacent} if there is a hexagon on $\Tb^*$ touched by both of them.

\begin{proposition}\label{prop:coup-1}
For any $j\ge 2$, there exists $\delta(j)>0$ such that for any $100j<10r<R$ and $m\in(1.1,10)$, denoting $u=\sqrt{rR}$, the following hold: 
	\begin{itemize}
		\item {\rm (Inward coupling)}  There is a coupling $\Qb$ of $\Pb[\;\cdot\mid \Hc_j(r,R)]$ and $\Pb[\;\cdot\mid \Hc_j(r,mR)]$ such that, if we sample $(\omega_1,\omega_2)\sim\Qb$, then with probability at least $\big(1-(r/R)^{\delta}\big)$, there exists a common configuration of outer faces $\Theta^*$ with $j$ faces around $C_u^+$ in both $\omega_1$ and $\omega_2$, and $\omega_1$ coincides with $\omega_2$ in $\Vc_{\Theta^*}$. Furthermore, when this is the case, $\Theta^*$ is a stopping set (recall the definition in Section \ref{subsec:faces}) and for any $\eta< r'<r$, we have
	    \begin{equation}\label{eq:DMP_1}
	    \Pb[\Hc_j(r',R)\mid \Hc_j(r,R),\Theta^*]=\Pb[\Hc_j(r',mR)\mid \Hc_j(r,mR),\Theta^*]\,.
	    \end{equation}
		\item {\rm (Outward coupling)} A similar coupling exists for $\Pb[\cdot\mid \Hc_j(r,R)]$ and $\Pb[\cdot\mid \Bc_j(r,R)]$ with outer faces replaced by inner faces. 
		For any $R'>R$, when the coupling succeeds (so that the common inner faces $\Theta^*$ exists), we have
		\begin{equation}\label{eq:DMP_2}
		\Pb[\Hc_j(r,R')\mid \Hc_j(r,R),\Theta^*]=\Pb[\Bc_j(r,R')\mid \Bc_j(r,R),\Theta^*]\,.
		\end{equation}
	\end{itemize}
\end{proposition}

As a consequence of these couplings in Proposition \ref{prop:one-coupling} (for the $j=1$ case) and Proposition \ref{prop:coup-1} (for the $j \geq 2$ cases), we know that the law $\Pb[ \cdot | \Hc_j(n^\alpha,n)]$ will coincide with $\Pb[ \cdot | \Hc_j(n^\alpha,mn)]$ at scale $n^\alpha$ with high probability for $\alpha\in(0,1)$\footnote{Not to be confused with the arm exponents $\alpha_j$. Same for $\beta$ in the proof of Proposition \ref{prop:couple-h-half} below.}. From the (slightly stronger) separation lemma Proposition~\ref{prop:half-sep-2}, we can show that the failure event will not contribute much to $\Pb[\Hc_j(r,n) | \Hc_j(n^\alpha,n)]$ much. Thus, we obtain the following comparison of conditional arm probabilities.
	\begin{proposition}
		\label{prop:couple-h-half}
		For any $r\ge r_h(j)$, $m\in (1.1,10)$ and $\alpha\in (0,1)$,
		\begin{equation}\label{eq:couple-h-half}
			\Pb\left[\Hc_j(r,n)\large{|}\Hc_j(n^{\alpha},n)\right] =\Pb\left[\Hc_j(r,mn)\large{|}\Hc_j(n^{\alpha},mn)\right] \Big(1+O(n^{-c})\Big)\,.
		\end{equation}
		Here, $O(n^{-c})$ may depend on $r$ and $\alpha$, but not $m$.
	\end{proposition}
	\begin{proof}
		We only give the proof for $j\ge 2$, for the case $j=1$ is similar (and easier indeed). Denote $\beta=(1+\alpha)/2$, and choose a large integer $K_0=K_0(j,\alpha)$ such that 
		\begin{equation}\label{eq:large-integer-K_0}
			\Pb[\Bc_{K_0}(n^{\alpha},n^{\beta})]\le n^{-\beta_j-1}.
		\end{equation}
		Consider probability measures $\Pb_1=\Pb[\cdot\mid \Hc_j(n^\alpha,n)]$ and $\Pb_2=\Pb[\cdot\mid \Hc_j(n^\alpha,mn)]$ with sample spaces $\Omega_1$ and $\Omega_2$. For a pair $(\omega_1,\omega_2)\in \Omega_1\times \Omega_2$, we say $\omega_1=_{n^\beta}\omega_2$ if there exists a common configuration $\Theta^*$ of outer faces around $C_{n^{\beta}}^+$, such that $\omega_1$ and $\omega_2$ are identical in $\Vc_{\Theta^*}$. Divide $\Omega_1\times \Omega_2$ into the union of  $\Kc_1=\{(\omega_1,\omega_2): \omega_1=_{n^\beta}\omega_2\}$, $\Kc_2=\big\{(\omega_1,\omega_2):\omega_1\in \Bc_{K_0}(n^\alpha,n^\beta)\mbox{ or }\omega_2\in \Bc_{K_0}(n^\alpha,n^\beta)\big\}\setminus\Kc_1$ and $\Kc_3=(\Kc_1\cup \Kc_2)^c$. By Proposition~\ref{prop:coup-1}, we can couple $\Pb_1,\Pb_2$ by $\Qb$ in a way that
		\begin{equation}
			\label{lemma3.10_1}
			\Qb[\Kc_1]\geq 1-n^{-C}\mbox{ for some }C(\alpha,j)>0\,.
		\end{equation}
		
		We define a pair of random stopping sets $(\mathbf{\Theta}_1,\mathbf{\Theta}_2)$ under $\Qb$ as follow. For $(\omega_1,\omega_2)\sim \Qb$, since $\Theta^*$ is a stopping set by Proposition~\ref{prop:coup-1}, we can perform an exploration process outside $C^+_{n^\beta}$ for both configurations to check that whether there exists a common configuration of outer faces $\Theta^*$ with $j$ faces for both $\omega_i,i=1,2$ and without exploring any hexagon in $\Vc_{\Theta^*}$ when this is the case. If so (see the left part of Figure~\ref{fig:succeedfail}), we then further explore the $(j-1)$ interfaces starting from the end-points of $\Theta^*$ and end upon reaching $C^+_{n^\alpha}$ and hence obtain a pair of configurations of outer faces $(\Theta_1(\omega_1),\Theta_2(\omega_2))$ around $C^+_{n^\alpha}$. Otherwise (see the right part of Figure~\ref{fig:succeedfail}), we further explore all the hexagons outside $C_{n^\beta}^+$ together with all the interfaces in $A^+(n^\alpha,n^\beta)$ from $C_{n^\beta}^+$ to $C_{n^\alpha}^+$ for both $\omega_i$, $i=1,2$. This will also yield a pair of configurations of outer faces $(\Theta_1(\omega_1),\Theta_2(\omega_2))$ around $C_{n^\alpha}^+$. Note that in either case, all hexagons in $\Vc_{\Theta_i(\omega_i)},i=1,2$ are left unexplored, hence we get a pair of stopping sets $(\mathbf{\Theta}_1,\mathbf{\Theta}_2)(\omega_1,\omega_2)=(\Theta_1(\omega_1),\Theta_2(\omega_2))$.
		
		\begin{figure}
			\centering
			\includegraphics[width=0.4\linewidth]{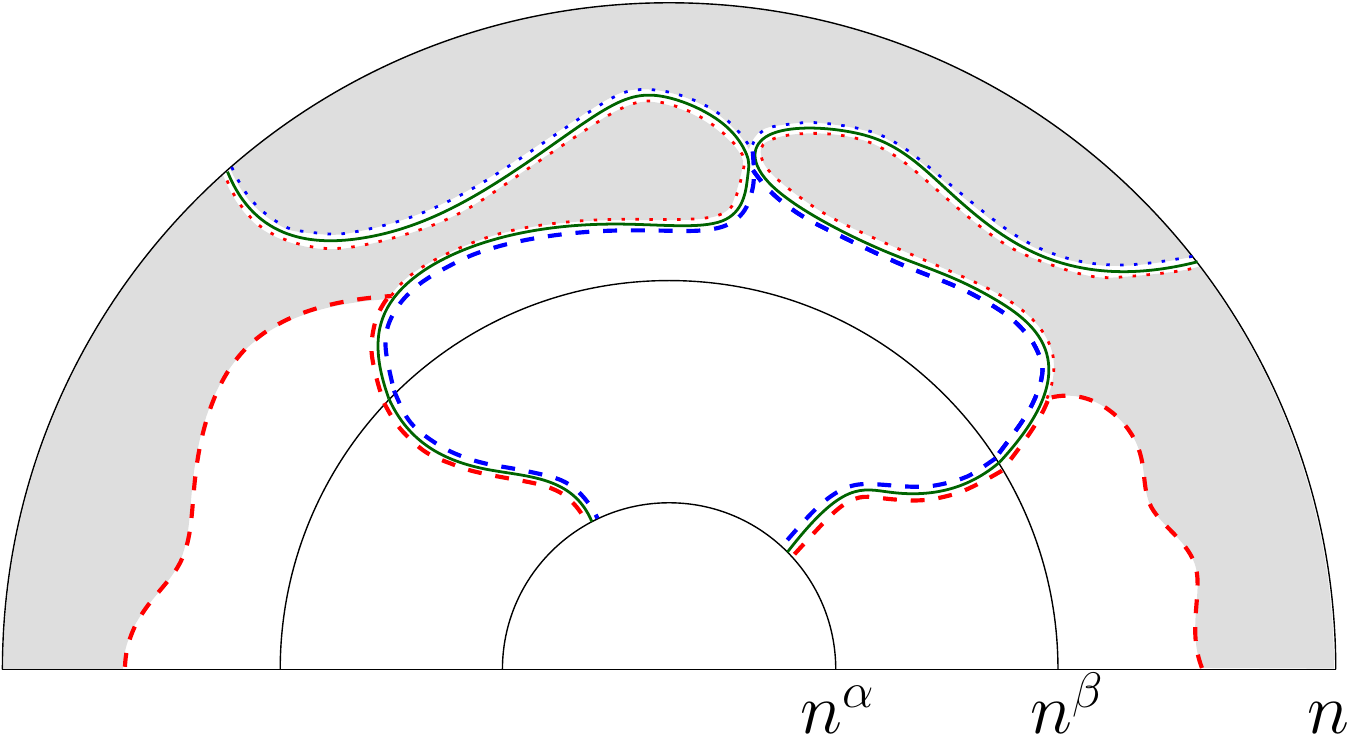} \quad\quad\quad \includegraphics[width=0.4\linewidth]{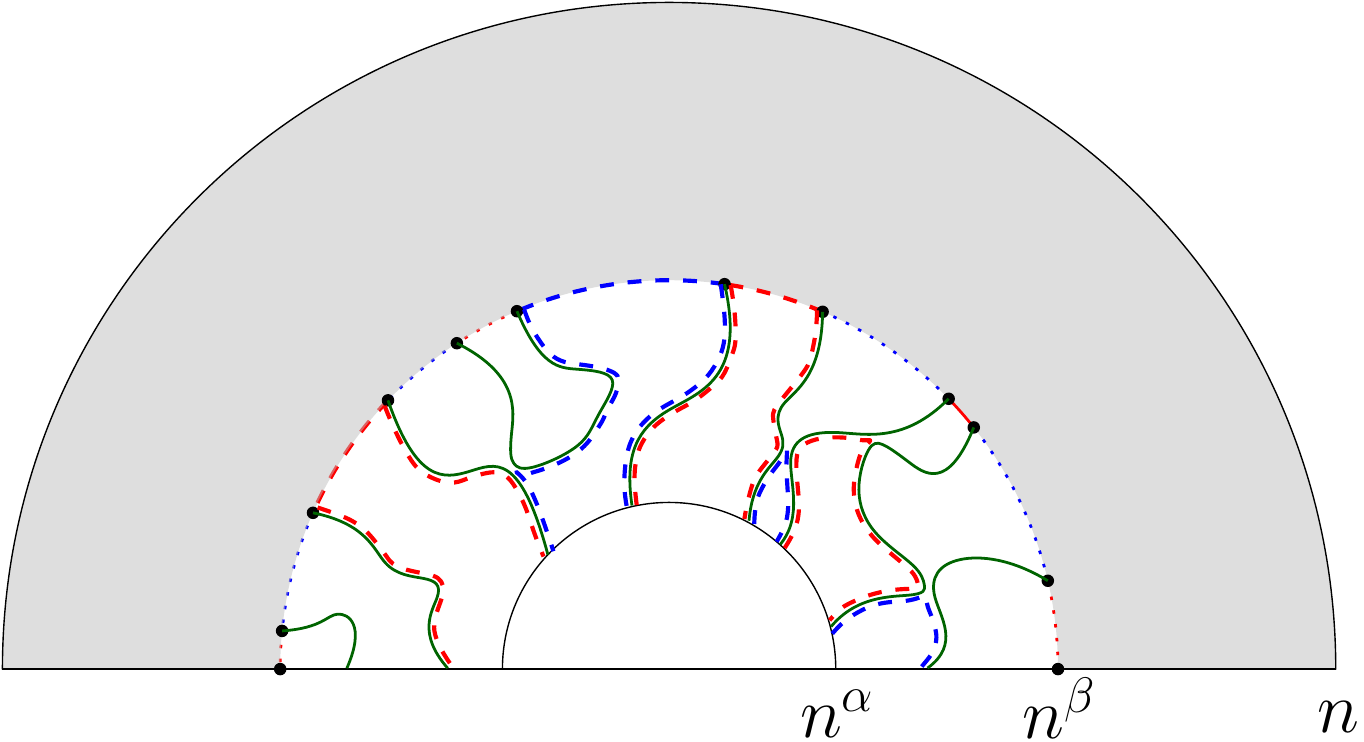}
				\caption{Proof of Proposition~\ref{prop:couple-h-half}. Sketch for $\Pb[\cdot\mid \Hc_3(n^\alpha,n)]$. \textbf{Left:} The case when the coupling succeeds: there is a common configuration of outer faces $\Theta^*$ with $3$ faces around $C^+_{n^{\beta}}$ for both $\omega_1\sim \Pb[\cdot\mid \Hc_j(n^\alpha,n)]$ and $\omega_2\sim\Pb[\cdot\mid \Hc_j(n^\alpha,mn)]$. The two interfaces that cross the annulus $A^+(n^{\alpha,n})$ are in green. Faces $\Theta^*$ together with hexagons that are adjacent to these two interfaces further create a common configuration of outer faces around $C^+_{n^{\alpha}}$, which is sketched in dashed red/blue curves. The region that has been explored ($\Vc_{\Theta^*}$) is in gray. \textbf{Right:} The case when the coupling fails. In this case, we reveal the color of all hexagons in $A^+(n^{\beta},n)$, which is illustrated in gray. Then, the color of all hexagons neighboring $C^+_{n^{\beta}}$ are revealed, which is sketched in dotted/dashed red/blue arcs on $C^+_{n^{\beta}}$. Running all of the exploration processes (in green) from $C^+_{n^{\beta}}$ towards $C^+_{n^{\alpha}}$, we obtain a configuration of outer faces around $C^+_{n^{\alpha}}$, which is sketched in dashed red/blue curves. In this case, the number of faces could be larger than $j$.}
				\label{fig:succeedfail}
		\end{figure}

		For $(\omega_1,\omega_2)\in \Kc_1$, $\Theta_1(\omega_1)$ and $\Theta_2(\omega_2)$ are identical with $j$ faces, and we deduce from \eqref{eq:DMP_1} that
		\begin{equation}\label{lemma3.10_4}
			\Pb[\Hc_j(r,n)\mid\Hc_j(n^\alpha,n), \omega_{\Dc_{\Theta_1(\omega_1)}}=\omega_1]=\Pb[\Hc_j(r,mn)\mid\Hc_j(n^\alpha,mn), \omega_{\Dc_{\Theta_2(\omega_2)}}=\omega_2]\,.
		\end{equation}
		For $(\omega_1,\omega_2)\in K_3$, we claim it holds uniformly that
		\begin{align}
			\Pb[\Hc_j(r,n)\mid \Hc_j(n^\alpha,n), \omega_{\Dc_{\Theta_1(\omega_1)}}=\omega_1]&\le n^{-\beta_j\alpha+o(1)}\,,\label{lemma3.10_2}\\
			\Pb[\Hc_j(r,mn)\mid\Hc_j(n^\alpha,mn), \omega_{\Dc_{\Theta_2(\omega_2)}}=\omega_2]&\le n^{-\beta_j\alpha+o(1)}\,.\label{lemma3.10_3}
		\end{align}
		Since $\Hc_j(r,n)\subset \Bc_j(r,2^{-K_0}n^\alpha)\cap \Bc_j(2^{-K_0}n^\alpha,n)$  and $\Bc_j(r,2^{-K_0}n^\alpha)$ is independent with $\Bc_j(2^{-K_0}n^\alpha,n)$ and $\{\omega_{\Dc_{\Theta_1(\omega_1)}}=\omega_1\}$, then uniformly for $(\omega_1,\omega_2)\in \Kc_3$, the left hand side of \eqref{lemma3.10_2} is bounded by
		\begin{equation*}
			\frac{\Pb[\Bc_j(r,2^{-K_0}n^\alpha)]\Pb[ \Bc_j(2^{-K_0}n^\alpha,n)\mid \omega_{\Dc_{\Theta_1(\omega_1)}}=\omega_1]}{\Pb[\Hc_j(n^\alpha,n)\mid \omega_{\Dc_{\Theta_1(\omega_1)}}=\omega_1]}\le n^{-\beta_j\alpha+o(1)}\cdot \frac{\Pb[\Bc_j(2^{-K_0}n^\alpha,n)\mid\omega_{\Dc_{\Theta_1(\omega_1)}}=\omega_1]}{\Pb[\Hc_j(n^\alpha,n)\mid \omega_{\Dc_{\Theta_1(\omega_1)}}=\omega_1]}\,.
		\end{equation*}
		Noting that $\Theta_1(\omega_1)$ has no more than $K_0$ faces for any $(\omega,\omega')\in \Kc_3$, and letting $\mathcal U$ be the event that the interfaces crossing $A^+(2^{-K_0}n^\alpha,n^\alpha)$ are well-separated on $C^+_{2^{-K_0}n^\alpha}$, it follows that 
		\begin{equation*}
			\begin{aligned}
			 \Pb[\Hc_j(n^\alpha,n)\mid\omega_{\Dc_{\Theta_1(\omega_1)}}=\omega_1]&\,\overset{\rm RSW}{\ge}\  c(j,K_0)\Pb[\mathcal U\cap \Bc_j(2^{-K_0}n^\alpha,n)\mid\omega_{ \Dc_{\Theta_1(\omega_1)}}=\omega_1]\\
			 &\overset{\eqref{eq:Q-in-2}}{\ge}\  c'(j,K_0)\Pb[\Bc_j(2^{-K_0}n^\alpha,n)\mid \omega_{ \Dc_{\Theta_1(\omega_1)}}=\omega_1]\,.
			\end{aligned}
		\end{equation*} 
This proves \eqref{lemma3.10_2}, and \eqref{lemma3.10_3} follows in the same fashion.
		
		Now, by Beyesian formula,
writing 	
$$
F:=\Pb[\Hc_j(r,n)\mid\Hc_j(n^\alpha,n),\omega_{ \Dc_{\Theta_1(\omega_1)}}=\omega_1]\,,\quad G:=\Pb[\Hc_j(r,mn)\mid\Hc_j(n^\alpha,mn), \omega_{\Dc_{\Theta_2(\omega_2)}}=\omega_2]
$$
for short, we have
		\begin{equation}\label{eq:threecases}
\Big|\Pb_1[\Hc_j(r,n)]-\Pb_2[\Hc_j(r,mn)]\Big|=			\Big|\sum_{(\Theta_1,\Theta_2)}(F-G)\Qb[(\mathbf{\Theta_1},\mathbf{\Theta_2})=(\Theta_1,\Theta_2)]\Big|\,.
		\end{equation}
The RHS of \eqref{eq:threecases} can be bounded by 
	\begin{equation*}
\Big|\sum_{(\Theta_1,\Theta_2)}\big(F-G\big)\Qb\big[(\mathbf{\Theta_1},\mathbf{\Theta_2})=(\Theta_1,\Theta_2),\,\Kc_1\big]\Big|
				+\sum_{i=2}^3\sum_{(\Theta_1,\Theta_2)}\big(F+G\big)\Qb\big[(\mathbf{\Theta_1},\mathbf{\Theta_2})=(\Theta_1,\Theta_2),\,\Kc_i\big].
\end{equation*}
		The first term equals  $0$ by \eqref{lemma3.10_4}, the second term (that corresponds to $\Kc_2$ in the sum) is bounded by $2\Qb[\Kc_2]\le 4n^{-\beta_j-1}$ from \eqref{eq:large-integer-K_0}, and the last term is bounded by $n^{-\beta_j\alpha-C+o(1)}$ from \eqref{lemma3.10_1}, \eqref{lemma3.10_2} and \eqref{lemma3.10_3}. As a result we see $|\Pb_1[\Hc_j(r,n)]-\Pb_2[\Hc_j(r,n)]|\le n^{-\beta_j\alpha-C+o(1)}$ for some $C>0$. In addition, by Lemma~\ref{prior_1} and  \eqref{eq:rmk_prior} in Remark~\ref{rmk:prior}, $\Pb\left[\Hc_j(r,n)\large{|}\Hc_j(n^{\alpha},n)\right] = n^{-\beta_j \alpha+o(1)}$. The claim \eqref{eq:couple-h-half} hence follows.	\end{proof}

	\begin{proposition}
		\label{prop:couple-bh-half}
		For any $r\ge r_h(j)(\ge r_b(j))$, $m\in (1.1,10)$ and $\alpha\in (0,1)$,
		\begin{equation}\label{eq:couple-bh-half}
			\Pb\left[\Bc_j(r,n)\large{|}\Bc_j(r,n^\alpha ) \right] =\Pb\left[\Hc_j(r, n)\large{|}\Hc_j(r, n^\alpha )\right] \Big(1+O(n^{-c})\Big).
		\end{equation}
	\end{proposition}
\begin{proof}
	In contrast to Proposition~\ref{prop:couple-h-half}, this is an outward coupling, in which case we do not need to spare an additional scale $n^{\beta}$ to use Proposition~\ref{prop:half-sep-2}, and thus much easier to deal with. Concretely, we first examine whether $\Pb\left[\cdot\mid\Bc_j(r,n^\alpha ) \right]$ and $\Pb\left[\cdot\mid\Hc_j(r,n^\alpha ) \right]$ induce the same configuration of inner faces around $C^+_{n^{\alpha}}$, which will happen with probability larger than $(1-n^{-C})$ by \eqref{eq:DMP_2}. If it fails, for both events $\Bc_j(r,n)$ and $\Hc_j(r,n)$, we further need to fulfill an arm event $\Bc_j(n^{\alpha},n)$ which has probability bounded by $n^{-\beta_j\alpha+o(1)}$ by \eqref{eq:rmk_prior}. The result follows immediately.
\end{proof}

\subsection{The plane cases}
In this subsection, we provide the coupling results in the plane, which are analogous to those in the half plane in the previous subsection. The one-arm and four-arm cases have already been developed in \cite{MR3073882}, but as new difficulties arise for the variants we consider in this work and for general $j$, we decide to include a proof, which we postpone to the last section. Recall the definitions of $\Xc_j,\Yc_j,\Ac_j$ for $j\ge 2$ in \eqref{def:y} and \eqref{def:a} and below \eqref{def:y} resp. We let $$J=2\lfloor j/2\rfloor.$$ 

In the next proposition, we couple together conditional laws in the plane case. For discussions of these couplings, see Proposition \ref{rmk:whole-coup} below.
\begin{proposition}\label{prop:C-in-coupling-even}
	For any $j\ge2$, there exists $\delta(j)>0$ such that for all $100j\eta\le 10r\le R$ and $m\in (1.1,10)$, denoting $u=\sqrt{rR}$, then the following hold.
	\begin{itemize}
			\item {\rm (Inward coupling)} There is a coupling $\Qb$ of the conditional laws $\Pb[\cdot\mid \Xc_j(r,R)]$ and $\Pb[\cdot\mid \Ac_j(r,R)]$ such that if we sample $(\omega_1,\omega_2)\sim \Qb$, then with probability at least $\big(1-({r}/{R})^{\delta}\big)$,  there exists a common configuration of outer faces $\Theta^*$ with $J$ faces around $C_{u}$ in both $\omega_1$ and $\omega_2$, and $\omega_1$ coincides with $\omega_2$  in $\Vc_{\Theta^*}$. Furthermore, when this is the case, $\Theta^*$ is a stopping set and for any $1\le r'\le r$, 
		\begin{equation}\label{eq:DMP_4}
			\Pb[\Xc_j(r',R)\mid \Xc_j(r,R),\omega_{\Dc_{\Theta^*}}=\omega_1]=\Pb[\Ac_j(r',R)\mid \Ac_j(r,R),\omega_{\Dc_{\Theta^*}}=\omega_2]\,.
		\end{equation}
		\item {\rm (Another inward coupling)} A similar coupling exists for the  the conditional laws $\Pb[\cdot\mid \Yc_j(r,R)]$ and $\Pb[\cdot\mid \Yc_j(r,mR)]$. In this case, if the coupling succeeds (i.e., $\Theta^*$ exists), 
		\begin{equation}\label{eq:DMP_3}
			\Pb[\Yc_j(r',R)\mid \Yc_j(r,R),\omega_{\Dc_{\Theta^*}}=\omega_1]=\Pb[\Yc_j(r',mR)\mid \Yc_j(r,mR),\omega_{\Dc_{\Theta^*}}=\omega_2]\,.
		\end{equation}
	\end{itemize}
\end{proposition}
\begin{remark}\label{rmk:whole-coup}
The major difference of these couplings with Proposition~\ref{prop:coup-1} lies in \eqref{eq:DMP_3}. In the half-plane case, the ``domain Markov property'' can be easily applied, and thus \eqref{eq:DMP_1} and \eqref{eq:DMP_2} follows naturally from the existence of common faces for both percolation configurations. In the  plane case, due to the complicacy of the events $\Xc_j(r,R)$ and $\Yc_j(r,R)$, \eqref{eq:DMP_3} does not hold trivially. We tackle this kind of difficulty in Section~\ref{sec:proofcoupling} by coupling together extra structures rather than just configurations of faces.
\end{remark}

As a result of Proposition~\ref{prop:C-in-coupling-even}, we get the following estimates on conditional arm probabilities.

\begin{proposition}
	\label{prop:couple-y-whole}
	For any $r\ge r_y(j)$, $m\in (1.1,10)$ and $\alpha\in (0,1)$, we have
	\begin{equation}
		\label{eq:couple-y-whole}
		\Pb\left[\Yc_j(r,n)\large{|}\Yc_j(n^{\alpha},n)\right] =\Pb\left[\Yc_j(r,mn)\large{|}\Yc_j(n^{\alpha},mn)\right] \Big(1+O(n^{-c})\Big)\,,
	\end{equation}
	where $O(n^{-c})$ may depend on $r$ and $\alpha$, but not $m$.
\end{proposition}
\begin{proposition}
	\label{prop:couple-A-whole-inner}
	For any $r\ge r_x(j)(=r_a(j))$ and $\varepsilon>0$, we have
	\begin{equation}
		\label{eq:couple-A-whole-inner}
		\Pb[\Xc_j(r,n)\mid \Xc_j(\varepsilon n,n)]=\Pb[\Ac_j(r,n)\mid \Ac_j(\varepsilon n,n)]\Big(1+O(\varepsilon^c)\Big)\,.
	\end{equation}
\end{proposition}
 The proof of Proposition~\ref{prop:couple-y-whole} is almost identical to that of Proposition~\ref{prop:couple-h-half}, and by similar arguments and replacing $n^\alpha$ with $\varepsilon n$, we obtain a similar proof of  Proposition~\ref{prop:couple-A-whole-inner}. We omit the proofs.

We now consider outward coupling in the plane. Recall the definition of $\Zc_j(r,R)$. There are similar couplings of $\Pb[\cdot\mid \Ac_j(r,n^\alpha)]$ and $\Pb[\cdot\mid \Zc_j(r,n^\alpha)]$ as in the half-plane case, from which we can deduce a comparison of the probability of $\Xc_j$ conditioned on $\Ac_j$ and $\Yc_j$ on $\Zc_j$, which is summarized in the following proposition. As the key idea is essentially same as before, for brevity we choose not to state the precise coupling results, and only give a sketch of the proof for the comparison result.
\begin{proposition}
	\label{prop:couple-a-whole}
	For any $r\ge r_z(j)(\ge r_a(j))$, $m\in (1.1,10)$ and $\alpha\in (0,1)$, we have
	\begin{equation}
		\label{eq:couple-a-whole}
		\Pb[\Xc_j(r,n)\mid \Ac_j(r,n^\alpha)]=\Pb[\Yc_j(r,n)\mid \Zc_j(r,n^\alpha)]\Big(D(j)+O(n^{-c})\Big)
	\end{equation}
   where $D(j)$ is a constant given by 
   \begin{equation}
   	\label{eq:C(j)}
   D(j)=
   \begin{cases}
   	1,\quad &j\equiv 1\pmod 2\\
   	\frac{j}{2},\quad &j\equiv 2\pmod 4\\
   	\frac{j}{4},\quad &j\equiv 0\pmod 4
   \end{cases}
   \end{equation}
   and $O(n^{-c})$ may depend on $r$, $\alpha$ and $j$, but not $m$.
   
\end{proposition}

\begin{figure}
	\begin{subfigure}{.24\textwidth}
		\centering
		\includegraphics[width=0.9\linewidth]{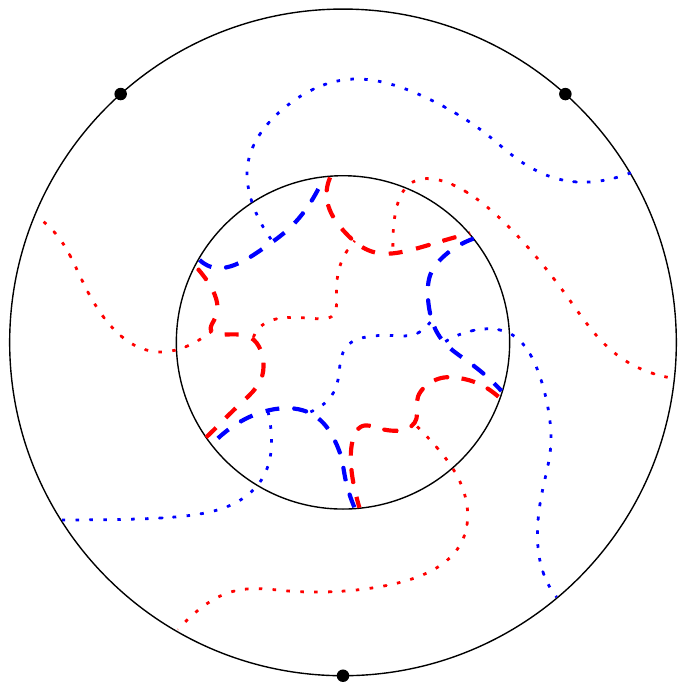}
		\caption{$\Yc_6(r,n)$.}
		\label{fig:Y6}
	\end{subfigure}
	\begin{subfigure}{.24\textwidth}
		\centering
		\includegraphics[width=0.9\linewidth]{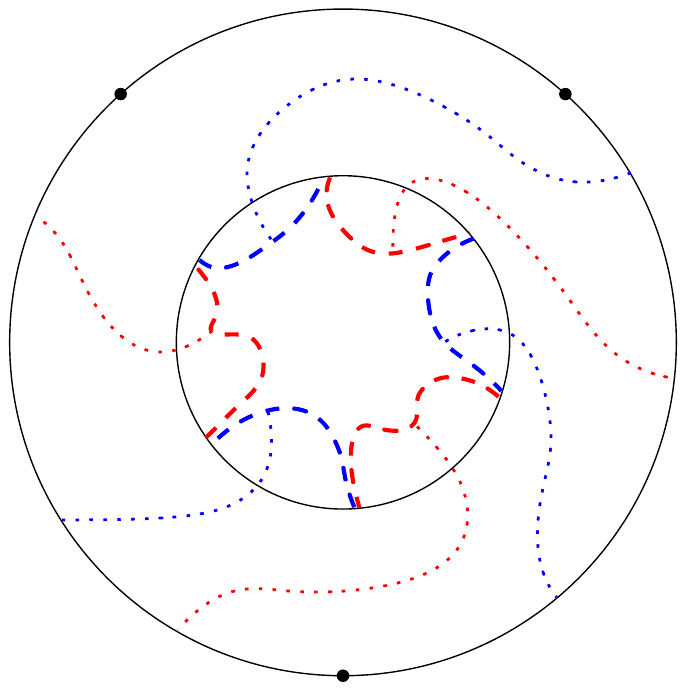}
		\caption{$\Xc^1_6(r,n)$.}
	\end{subfigure}
	\begin{subfigure}{.24\textwidth}
		\centering
		\includegraphics[width=0.9\linewidth]{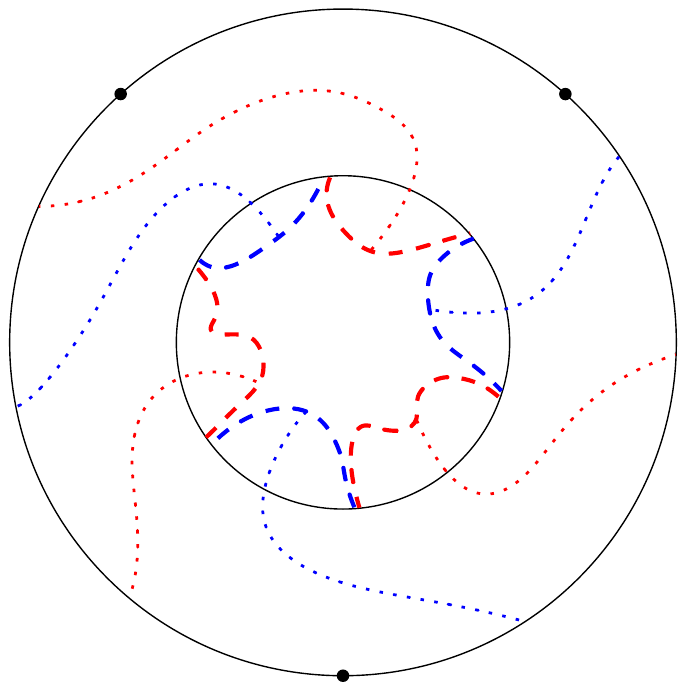}
		\caption{$\Xc^2_6(r,n)$.}
	\end{subfigure}
	\begin{subfigure}{.24\textwidth}
		\centering
		\includegraphics[width=0.9\linewidth]{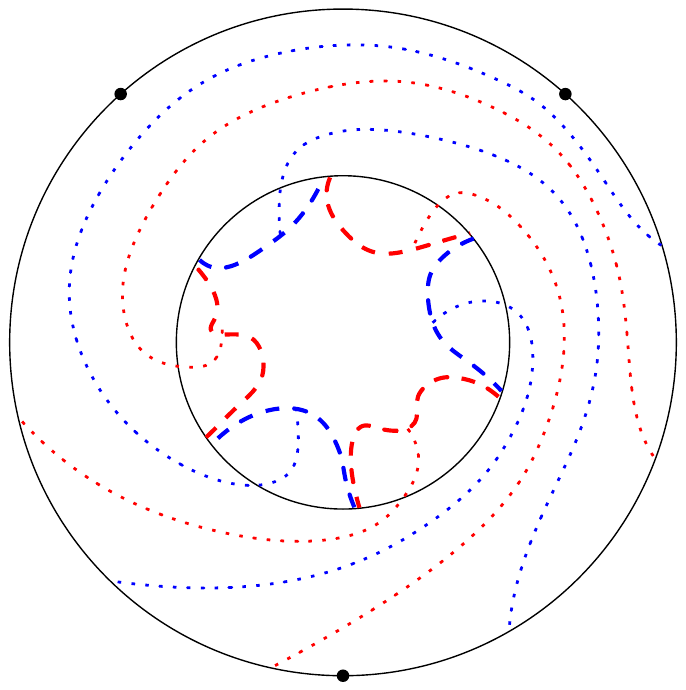}
		\caption{$\Xc^3_6(r,n)$.}
	\end{subfigure}
	\caption{An illustration of the events $\Yc_6(r,n)$ and $\Xc_6(r,n)$ of three different types of connecting patterns. All four annuli are $A(n^{\alpha},n)$ with the same configuration of inner faces $\Theta^*$ in dashed red/blue curves. The three marked points on $C_n$ are $n\cdot a$, $n\cdot b$ and $n\cdot c$, respectively (from bottom, counterclockwise). The outer connecting patterns for $\Theta^*$ are sketched in dotted red/blue curves, and there is an additional inner connecting pattern in the first picture which results in the only possibility for the outer connection for $\Yc_6(r,n)$.}
	\label{fig:C(j)}
\end{figure}

We first explain how the constant $D(j)$ in \eqref{eq:C(j)} appears. As in Proposition~\ref{prop:C-in-coupling-even}, we can couple $\Pb[\cdot\mid \Ac_j(r,n^\alpha)]$ and $\Pb[\cdot\mid \Zc_j(r,n^\alpha)]$ in a way such that if we sample $(\omega_1,\omega_2)$ according to this coupling, $\omega_1$ and $\omega_2$ will coincide on a configuration of inner faces $\Theta^*$ around $C_{n^\alpha}$ with high probability. We now consider the events $\Xc_j(r,n)$ and $\Yc_j(r,n)$  conditioned on the color configuration in $\Dc_{\Theta^*}$. When $j$ is odd, conditioned on $\Dc_{\Theta^*}$, $\Xc_j(r,n)$ is equivalent to $\Yc_j(r,n)$ and the proof is almost identical with the previous cases. However, when $j$ is even, these two events are distinct: there are $j/2$ kinds of connecting patterns for $\Xc_j(r,n) $, since the information inside $\Theta^*$ cannot tell us which faces should be connected to the corresponding segments on $C_n$; but for $\Yc_j(r,n)$, there are only $1$ (if $j\equiv 2\pmod 4$) or $2$ (if $j\equiv 0\pmod 4$) possible connecting patterns, because the information inside $\Theta^*$ already determines the order of faces in $\Theta^*$. See Figure~\ref{fig:C(j)} for an illustration of the case when $j=6$. If we realize that all possible connecting patterns have nearly the same probability thanks to the coupling results, then the definition of $D(j)$ in \eqref{eq:C(j)} makes sense.

\begin{proof}[Sketch of Proof for Proposition~\ref{prop:couple-a-whole} when $j$ is even]
Denote $\Gc$ for the event that there are exactly $j$ interfaces crossing $A(r,n)$, then comparing with $\Xc_j(r,n)$ and $\Yc_j(r,n)$, $\Gc^c$ has negligible probability by BK-Reimer's inequality. We partition $\Xc_j(r,n)\cap \Gc$ into the disjoint union of $D(j)$ events according to the locations  of interfaces. More precisely, let $\hat{\Gamma}$ be the set of the $\frac{j}{2}$ interfaces with red on their left (seen from outside to inside). We label each $\gamma\in \hat{\Gamma}$ in the following two ways: label $\gamma$ by $\gamma_i$ (resp.\ $\gamma^i$)  with $1\le i\le \frac{j}{2}$, if $\gamma$ is the $i$-th element in $\hat{\Gamma}$ when counting around the circle $C_r$ (resp.\ $C_n$) counterclockwise starting from $(0,{-r})$ (resp.\ $(0,{-n})$). For $1\le i\le D(j)$, define $\Xc_j^i(r,n)$ as 
\begin{equation}
	\Xc_j^i(r,n)=\begin{cases}
		\Xc_j(r,n)\cap\Gc\cap\{\gamma_1=\gamma^i\}\quad &j\equiv 2 \pmod 4\,,\\
	\Xc_j(r,n)\cap \Gc\cap	\{\gamma_1=\gamma^i\mbox{ or }\gamma^{i+\frac{j}{4}}\}\quad &j\equiv 0\pmod 4\,.
	\end{cases}
\end{equation}  
See Figure~\ref{fig:C(j)} for illustration of the events $\Xc_j^i(r,n)$'s.
Then $\Xc_j(r,n)\cap \Gc$ is the disjoint union of $\Xc_j^i(r,n),1\le i\le D(j)$. We construct a coupling of $\Pb[\;\cdot\mid \Ac(r,n^\alpha)]$ and $\Pb[\;\cdot\mid \Zc_j(r,n^\alpha)]$ for each $1\le i\le D(j)$.
In such couplings, we keep a record of the relative location of interfaces, so once we successfully couple a face, the conditional version of $\Xc_j^i(r,n)$ and $\Yc_j(r,n)\cap \Gc$ are identical. This allows us to deduce that $$\Pb[\Xc_j^i(r,n)\mid \Ac_j(r,n^\alpha)]=\Pb[\Yc_j(r,n)\cap \Gc\mid \Zc_j(r,n^\alpha)]\Big(1+O(n^{-c})\Big)$$ holds for each $1\le i\le D(j)$, and the desired result follows readily.
\end{proof}
	\section{Proof of main theorems}
\label{sec:compare}
In this section, we will combine the coupling results from the previous section and the power-law convergence of the exploration process discussed in Section \ref{subsec:exploration} to prove  the main results of this work, namely Theorems~\ref{thm:newes-half}, \ref{thm:newes-whole} and \ref{thm:newes2-whole}. In Subsections~\ref{subsec:4.1} and \ref{subsec:4.2}, we will prove Theorem~\ref{thm:newes-half}, which is derived directly from Proposition~\ref{prop:ratioestimate-half}. In order to prove Proposition~\ref{prop:ratioestimate-half}, we use Propositions~\ref{prop:couple-h-half}, \ref{prop:couple-bh-half} and \ref{prop:compare-half} as inputs. In Subsections~\ref{subsec:4.3} and \ref{subsec:4.4}, we will prove Theorems~\ref{thm:newes-whole} and \ref{thm:newes2-whole}. The former is derived directly from Proposition~\ref{prop:ratioestimate-whole}, while the latter theorem is a quick corollary of the former. In order to prove Proposition~\ref{prop:ratioestimate-whole}, we use Propositions~\ref{prop:couple-y-whole} and \ref{prop:compare-whole} (the latter is an analogue of Proposition~\ref{prop:compare-half} in the plane), as inputs. 

Note that in Sections \ref{subsec:4.1} and \ref{subsec:4.3}, we will consider the rescaled lattice and the discretization scheme in Section \ref{subsec:setting} comes into play.

\subsection{Comparison estimates in the half-plane}
\label{subsec:4.1}
In this subsection, our main goal is Proposition~\ref{prop:compare-half} in which we compare the probabilities of $\Hc_j$, the modified half-plane arm events (recall \eqref{def:h} for definitions), at different scales. 

As discussed in Section \ref{subsec:comments}, the key ingredient is the convergence rate of percolation exploration process towards ${\rm SLE}_6$ upon some stopping time as shown in Theorem 4.1.10 in \cite{richards2021convergence} and discussed in Section \ref{subsec:exploration}. Instead of citing Proposition \ref{prop:couple-SLE} directly, we will use a modified version stated as follows. 

Recall the notation from Section \ref{subsec:exploration}. We consider the nice domain (in the sense of Section \ref{subsec:setting}) $$\Omega=B^+_1 \cup B^-\Big((-3/4,0),1/4\Big),\,a=(1,0)\mbox{ and }b=(-3/4,-1/4);
$$ see Figure~\ref{figure:explore-half} for an illustration. We assign red color to the counterclockwise arc $\wideparen{ab}$ and blue color to the counterclockwise arc $\wideparen{ba}$. Given mesh size $\eta>0$, the percolation exploration process starts from $a_\eta$ and end at $b_\eta$. Given $\alpha \in (0,1)$ and $m>0$, let $\widetilde{T}_{1/n}$ (resp. $\widetilde{T}_{1/(mn)}, \widetilde{T}$) denote the first time that $\gamma_{1/n}$ (resp. $\gamma_{1/(mn)}, \gamma$) hits $[-1,-n^{\alpha-1})$.

Recall \eqref{def:h} that $\Hc_j(r,R)$ is the event that in $B^+_R$ with the lattice scale $\eta=1$ there exist $j$ disjoint arms connecting $[-r,r] $ to $C^+_R$ with alternating colors. Rescaling the lattice by $\eta=1/R$,\footnote{In fact we will consider the cases $R=n$ and $R=mn$ in this and the next subsection.} we get that $\Hc_j(r,R)$ is equivalent to the event that in $B^+_1$ with the lattice scale $\eta$, there exist $j$ such arms connecting $ [-r/R , r/R]$ to $C^+_1$.  From here till the end of this subsection, with slight abuse of notation, we will tacitly assume the equivalence between events on the rescaled and those on the unrescaled lattices.
\begin{proposition}
\label{prop:couple-SLE-half}
There exists $u(\Omega) >0$ such that for any $\alpha \in (0,1)$ there exists a coupling of $\gamma_{1/n}$ and $\gamma$ (which we will refer to as the good coupling)
\begin{equation*}
\Pb\left[d\left(\gamma_{1/n}\big|_{[0,\widetilde{T}_{1/n}]}, \gamma\big|_{[0,\widetilde{T}]}\right) >n^{-u}\right]=O(n^{-u}) \,.
\end{equation*}
One can further check that $O(n^{-u})$ is independent of $\alpha \in (0,1)$.
\end{proposition}
\begin{proof}
The proof is the same as that of Proposition \ref{prop:couple-SLE} except that here we let $\Fc$ denote the event that there does not exist $x$ in $[-1,-n^{\alpha-1})$ such that  $\gamma_{1/n}$ enters $B(x, n^{-c})$ ($c(\Omega)$ is the constant associated with the coupling in Theorem 4.1.10 of \cite{richards2021convergence}) but does not hit $[-1,-n^{\alpha-1}) $ between the last entrance and first exit times of $B(x, n^{-c'})$, and the same for $\gamma$. We can still choose $c'$ small such that $\Pb[\Fc^c] \leq O(n^{-c''})$ and show that under good coupling of $\gamma_{1/n}$ and $\gamma$
$$\Fc \bigcap \left\{ \left(\gamma_{1/n}\big|_{[0,T'_{1/n}]}, \gamma\big|_{[0,T']}\right) \leq n^{-c} \right\} \mbox{\LARGE$\subset$}
 \left\{ d\left(\gamma_{1/n}\big|_{[0,\widetilde{T}_{1/n}]}, \gamma\big|_{[0,\widetilde{T}]}\right) \leq n^{-c'} \right\} \,.$$ Hence, 
\begin{equation*}
\Pb\left[d\left(\gamma_{1/n}\big|_{[0,\widetilde{T}_{1/n}]}, \gamma\big|_{[0,\widetilde{T}]}\right) >n^{-c'}\right] \leq \Pb\left[d\left(\gamma_{1/n}\big|_{[0,T'_{1/n}]}, \gamma\big|_{[0,T']}\right) >n^{-c}\right] +\Pb[\Fc^c] \leq O(n^{-c })+O(n^{-c''}) \,.
\end{equation*}
Taking $u= \min \{c,c',c'' \}$, we obtain the proposition.
\end{proof}

Similarly, for the same constant $u$, we can show that one can further couple $\gamma_{1/mn}$ and $\gamma$ (we still denote the law by $\Pb$, view it as the law of $\gamma_{1/n}$, $\gamma_{1/mn}$ and $\gamma$ coupled together and call it the good coupling) such that
\begin{equation}
\label{eq:prop4.1-1}
\Pb\left[d\left(\gamma_{1/(mn)}\big|_{[0,\widetilde{T}_{1/(mn)}]}, \gamma\big|_{[0,\widetilde{T}]}\right) >n^{-u}\right]=O(n^{-u}) \,.
\end{equation}One can further check that $O(n^{-u})$ is independent of $\alpha \in (0,1)$ and $m$ in a finite interval away from zero, e.g., $(1.1,10)$. 

Recall the definition of  $h_j(r,R)$ in \eqref{def:sequence}. 
The following proposition states that $h_j$ are almost the same at different scales within a power-law error at different scales when $r$ is not too small compared to $R$.

\begin{proposition}
\label{prop:compare-half}
There exists $c_8>0$ such that for all $ \alpha \in (1-c_8,1)$ and $m \in (1.1,10)$,
\begin{equation*}
h_j(n^\alpha,n)=h_j(mn^\alpha,mn)\Big(1+O(n^{-c})\Big) \,.
\end{equation*}
Importantly, here $O(n^{-c})$ is independent of $m$ and $\alpha$!
\end{proposition}

To prove it, we need Lemmas~\ref{lemma4.2-1} and \ref{lemma4.2-2}. In Lemma~\ref{lemma4.2-1}, we will show that $\Hc_j(n^\alpha,n)$ can be determined by the exploration process $\gamma_{1/n}$ upon time $\widetilde{T}_{1/n}$ and similarly $\Hc_j(mn^\alpha, mn)$ can be determined by $\gamma_{1/mn}$ upon $\widetilde{T}_{1/mn}$. In Lemma~\ref{lemma4.2-2}, we will show that the indicator functions of corresponding events of the exploration process are identical under three constraints. Finally, we prove that all of these constraints happen with high probability under the good coupling of $\gamma_{1/n}$, $\gamma_{1/mn}$ and $\gamma$. Therefore, the probabilities of $\Hc_j(n^\alpha,n)$ and $\Hc_j(mn^\alpha,mn)$ are close. 
\begin{lemma}

\begin{figure}[H]
	\centering
	\includegraphics[scale=0.8]{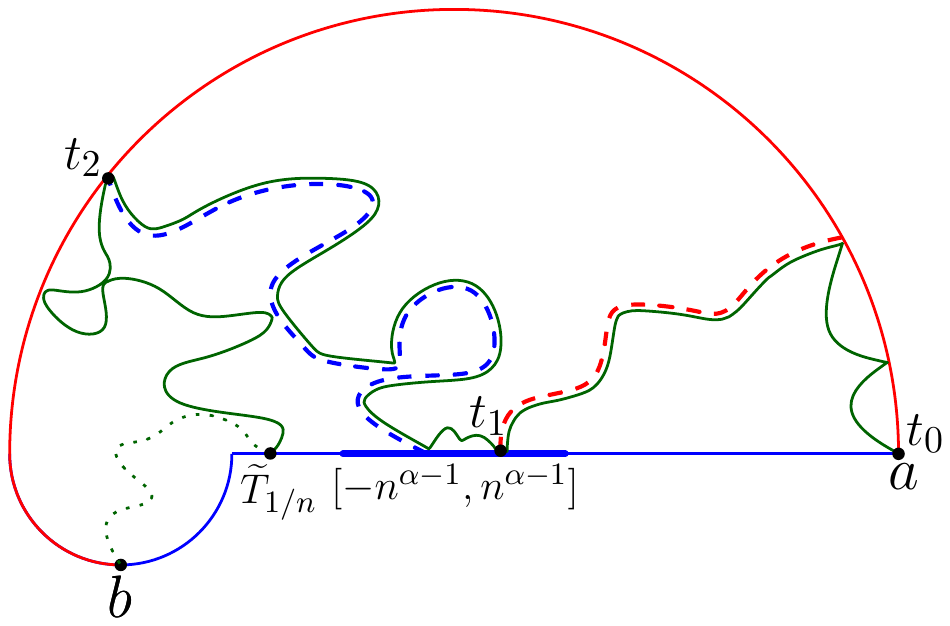}
	\caption{Sketch of the event $\Hc_j(n^\alpha,n)$ when $j=2$. The whole domain is $\Omega_{[1/n]}$. The exploration process $\gamma_{1/n}$ from $a_{1/n}$ to $b_{1/n}$ on $\Omega_{[1/n]}$ is in green. $\wt T_{1/n}$ is the first time that $\gamma_{1/n}$ hits $[-1,-n^{\alpha-1})$, after which time $\gamma_{1/n}$ is sketched in dotted. The segment $[-n^{\alpha-1},n^{\alpha-1}]$ is in bold blue line. Two arms, neighboring the interface $\gamma_{1/n}$, from this segment to the top boundary $C^+_{1}$ are in dashed red and blue respectively.}
	\label{figure:explore-half}
\end{figure}

\label{lemma4.2-1}
One has
\begin{equation}
\label{eq:behavior_gamma_n}
\Hc_j(n^\alpha,n) = \{ \gamma_{1/n} \mbox{ travels }j \mbox{  times between } C^+_1 \mbox{ and } [-n^{\alpha-1},n^{\alpha-1}] \; \mbox{before } \widetilde{T}_{1/n} \} 
\end{equation}
and
\begin{equation}
\label{eq:behavior_gamma_mn}
\Hc_j(mn^\alpha,mn) = \{ \gamma_{1/mn} \mbox{ travels }j\mbox{ times between } C^+_1 \mbox{ and } [-n^{\alpha-1},n^{\alpha-1}] \;  \mbox{before } \widetilde{T}_{1/mn} \} \,.
\end{equation}
\end{lemma}
By the RHS of \eqref{eq:behavior_gamma_n}, we mean that there exist $0=t_0<t_1<t_2<\ldots<t_j \leq \widetilde{T}_{1/n}$ such that $\gamma_{1/n}(t_0) \in C^+_1$, $\gamma_{1/n}(t_1) \in [-n^{\alpha-1},n^{\alpha-1}] $, $\gamma_{1/n}(t_2) \in C^+_1$, etc.; the same for \eqref{eq:behavior_gamma_mn}. The proof follows from simple geometrical arguments, so we omit it. See also Figure~\ref{figure:explore-half} for an illustration. Also, note that if we work with $\Bc_j(n^\alpha,n)$ instead, we need to deal with a varying domain which is a less-than-ideal strategy. This is the main reason that we introduce and work with $\Hc_j(n^\alpha,n)$. We refer readers to Remark \ref{rem:HjnotBj} and the paragraph below Proposition \ref{prop:ratioestimate-half} for more discussions.

Before stating the second lemma, we introduce three events which all happen with high probability under good coupling of $\gamma_{1/n}$, $\gamma_{1/mn}$ and $\gamma$. Let $u = u(\Omega)$ be the constant from Proposition~\ref{prop:couple-SLE-half}. We will call $C^+_1 \cup [-n^{\alpha-1}, n^{\alpha-1}]$ the {\bf designated boundary} and note that $(-1,0)$, $(1,0)$, $(-n^{\alpha-1},0)$ and $(n^{\alpha-1},0)$ are the extremal points of the designated boundary. Define $\Jc_1,\Jc_2,\Jc_3$ as
\begin{align*}
\Jc_1 &= \Big{\{} d\Big(\gamma_{1/n}\big|_{[0,\widetilde{T}_{1/n}]}, \gamma_{1/mn}\big|_{[0,\widetilde{T}_{1/mn}]}\Big) \leq 2n^{-u} \Big{\}};\\
\Jc_2 &= \Big\{ \mbox{There does not exist } x \mbox{ on the boundary such that it is }n^{-v} \mbox{ away from the extremal}\\
&\quad \quad \mbox{points, and }\gamma_{1/n}\mbox{ enters }B(x,2n^{-u}) \mbox{ but does not hit the designated boundary between }\\
&\quad \quad \mbox{the last entrance  and first exit times of } B(x,n^{-v}),\mbox{ and the same for }\gamma_{1/mn} \Big\} \,;\\
\Jc_3 &= \Big\{ \gamma_{1/n} \mbox{ and }\gamma_{1/mn}\mbox{ do not enter }B\big((-1,0), 2n^{-v}\big),\; B\big((-n^{\alpha-1},0), 2n^{-v}\big)\\
&\quad \quad \mbox{ or }  B\big((n^{\alpha-1},0), 2n^{-v}\big)\mbox{ and do not reenter }B\big((1,0),2n^{-v}\big) \mbox{ after leaving }B\big((1,0),1/4\big) \Big\}\,.
\end{align*}
Here, $v \in (0,u)$ is a constant such that \eqref{eq:choose-v-half} holds. If $\Jc_2$ happens, then each time $\gamma_{1/n}$ (or $\gamma_{1/mn}$) enters $B(x, 2n^{-u})$, it will hit the designated boundary between the last entrance and first exit times of $B(x,n^{-v})$  for all points $x$ on the designated boundary that are $n^{-v}$ away from the extremal points. 

We now give the second lemma which states that the indicator functions of $\Hc_j(n^\alpha,n)$ and $\Hc_j(mn^\alpha,mn)$ are identical under three constraints. In the proof, we will use the equivalence established in Lemma~\ref{lemma4.2-1}.
\begin{lemma}
\label{lemma4.2-2} One has
\begin{equation}
\label{eq:indentity-half} 
\mathbbm{1}_{\Hc_j(n^\alpha,n)}=\mathbbm{1}_{\Hc_j(mn^\alpha,mn)} \mbox{ on }\Jc_1 \cap \Jc_2 \cap \Jc_3\,.
\end{equation}
\end{lemma}

\begin{proof}
Assume that $\Jc_1,\Jc_2,\Jc_3$ and $\Hc_j(n^\alpha,n)$ happen so that $\gamma_{1/n}$ fulfills the RHS of \eqref{eq:behavior_gamma_n}. Suppose that $\gamma_{1/n}$ hits $[-n^{\alpha-1},n^{\alpha-1}] $ at the point $x$. From $\Jc_3$ we know that $x$ is $n^{-v}$ away from the extremal points. By $\Jc_1$ and $\Jc_2$, the process $\gamma_{1/mn}$ enters $B(x,2n^{-u})$ and so must hit the designated boundary 
between the last entrance and first exit times of $B(x,n^{-v})$, therefore also before $\widetilde{T}_{1/mn}$. Thus, $\gamma_{1/mn}$ must hit $[-n^{\alpha-1},n^{\alpha-1}]$ at the same place as $\gamma_{1/n}$ before $\widetilde{T}_{1/mn}$. Similarly, when $\gamma_{1/n}$ hits a point $x$ in $C^+_1$ before $\widetilde{T}_{1/n}$ and after hitting $[-n^{\alpha-1},n^{\alpha-1}]$ (this together with $\Jc_3$ ensures that the hitting point is $2n^{-v}$ away from the point $(1,0)$), the process $\gamma_{1/mn}$ must hit $C^+_1$ at the same place as $\gamma_{1/n}$ before $\widetilde{T}_{1/mn}$. Hence, $\gamma_{1/mn}$ fulfilling the RHS of \eqref{eq:behavior_gamma_mn} and so $\Hc_j(mn^\alpha,mn)$ happens. Therefore, $$\Jc_1 \cap \Jc_2 \cap \Jc_3 \cap \Hc_j(n^\alpha,n) \subset \Hc_j(mn^\alpha,mn)\mbox{ and similarly }\Jc_1 \cap \Jc_2 \cap \Jc_3 \cap \Hc_j(mn^\alpha,mn) \subset \Hc_j(n^\alpha,n)\,.$$ This completes the proof of \eqref{eq:indentity-half}.
\end{proof}
\begin{proof}[Proof of Proposition~\ref{prop:compare-half}]
First, we give upper bounds to $\Pb[\Jc_1^c], \Pb[\Jc_2^c]$ and $\Pb[\Jc_3^c]$. By Proposition \ref{prop:couple-SLE-half} and \eqref{eq:prop4.1-1}, under the good coupling of $\gamma_{1/n}, \gamma_{1/mn}$ and $\gamma$
\begin{equation*}
\Pb[\Jc_1^c]= \Pb\left[d\left(\gamma_{1/n}\big|_{[0,\widetilde{T}_{1/n}]}, \gamma_{1/mn}\big|_{[0,\widetilde{T}_{1/mn}]}\right) >2n^{-u}\right] =O(n^{-u}) \,.
\end{equation*}
We claim that if $\Jc_2^c$ happens, we have a half-plane $3$-arm event from a $(2n^{-u)}$-ball on the designated boundary to distance $n^{-v}$. Suppose $\Jc_2^c$ happens for $\gamma_{1/n}$ and some point $x$ on the designated boundary (the case of $\gamma_{1/mn}$ can be handled similarly). We define $\tau^1$ as the last\footnote{Here, last and first are defined with respect to some hitting time of $B(x,2n^{-u})$.} entrance time of $B(x,n^{-v})$, $\sigma^1$ as the first hitting time of $B(x,2n^{-u})$ after $\tau^1$, $\sigma^2$ as the first exit time of $B(x,n^{-v})$ and $\tau^2$ as the last exit time of $B(x,2n^{-u})$ before $\sigma^2$.  Then, the left and right boundaries of $\gamma_{1/n}[\tau^1,\sigma^1]$ and $\gamma_{1/n}[\tau^2, \sigma^2]$ contain three disjoint crossings from $B(x,2n^{-u})$ to distance $n^{-v}$. In other words, there are two crossings of the same color at the two boundaries close to $C^+_1 \cup [-n^{\alpha-1}, n^{\alpha-1}] $ and another crossing of a different color in the area sandwiched by the process. Moreover, we can cover the designated boundary by a $(2n^{-u})$-net with $O(n^u)$ elements. So, for some $c>0$
\begin{equation}
\label{eq:proof-3arm}
\Pb[\Jc_2^c] \leq O(n^u) \times O(n^{(1+c)(v-u)}) \,.
\end{equation}If $\Jc_3^c$ happens, we have a planar $1$-arm event from a $(2n^{-v})$-ball centered at one of the extremal points to distance $1/4$ since we can find an arm along the exploration process that is $2n^{-v}$ close to that extremal point. So, $$\Pb[\Jc_3^c] \leq O(1) \times O(n^{-cv})$$(where we reduce the value of $c$ if necessary).

Therefore, under the good coupling
\begin{equation*}
\begin{split}
|h_j(n^\alpha,n)-h_j(mn^\alpha,mn)| &\overset{\eqref{eq:indentity-half}}{\leq} 1- \Pb[\Jc_1^c \cup \Jc_2^c \cup \Jc_3^c] \;\leq\; \Pb(\Jc_1^c) +\Pb(\Jc_2^c)+\Pb(\Jc_3^c) \\
&\;\,=\;\, O(n^{-u}) +O(n^{u}) \times  O\big(n^{(1+c)(v-u)}\big)+ O(n^{-c v})\,.
\end{split}
\end{equation*}
Take a small $v$ such that 
\begin{equation}
\label{eq:choose-v-half}
u+(1+c)(v-u) <0\,.
\end{equation}By applications of the RSW theory, $h_j(mn^\alpha,mn) \geq n^{C (\alpha-1)}$. We can take $c_8$ such that $$-C c_8 >\max \{ -u, u+(1+c)(v-u),-cv \}\,.$$ Then, for all $\alpha \in (1-c_8,1)$, $|h_j(n^\alpha,n)-h_j(mn^\alpha,mn)| = O(n^{-c}) h_j(mn^\alpha,mn)$ as desired.
\end{proof}
\begin{remark}
When $j=1$, the event $\Hc_j(r,R)$ is simply equivalent to the event that there is a red crossing between $[-r/R,r/R]$ and $C^+_1$ in $B^+_1$ with the lattice scale $\eta=1/R$. This special case of Proposition~\ref{prop:compare-half} has been proved in Proposition 5.6 of \cite{mendelson2014rate} and the Main Theorem of \cite{binder2015rate}.
\end{remark}
 
\subsection{Proof of Theorem~\ref{thm:newes-half}}
\label{subsec:4.2}
In this subsection, we will prove Proposition~\ref{prop:ratioestimate-half} and complete the proof of Theorem~\ref{thm:newes-half}. We will first prove a version of \eqref{eq:ratioestimate-half} for the case of $h_j$ with Propositions \ref{prop:compare-half} and \ref{prop:couple-h-half} as inputs, whose proof is simpler than that of \eqref{eq:ratioestimate-half} for $b_j$ yet already contains the main idea. Then, we will turn back to \eqref{eq:ratioestimate-half} for $b_j$ with Proposition~\ref{prop:couple-bh-half} as additional inputs.
\begin{proof}[Proof of Proposition~\ref{prop:ratioestimate-half}]
As discussed above, we start with the case of $h_j$. Fix $r\ge r_h(j)$. Let $m \in (1.1,10)$ and $\alpha \in (1-c_8,1)$ ($c_8$ is the constant defined in Proposition~\ref{prop:compare-half}). In this proof we write $f(n)\simeq g(n)$ as the shorthand of $f(n)=g(n)\big(1+O(n^{-c})\big)$ where the constants may depend on $j$ and on the choice of $\alpha$ but are independent of the choice of $m$.

We first transform $\frac{h_j(r,mn)}{h_j(r,n)}$ into comparisons of mesoscopic arm probabilities thanks to Proposition~\ref{prop:couple-h-half}. Since $\Hc_j(r,n) \subset \Hc_j(n^{\alpha},n)$ and $\Hc_j(r,mn) \subset \Hc_j(n^{\alpha},mn)$, 
\begin{equation}
\label{eq:thm1.1-1}
\frac{h_j(r,mn)}{h_j(r,n)} = \frac{\Pb\left[\Hc_j(r,mn)\large{|}\Hc_j(n^{\alpha},mn)\right]  \cdot \Pb\left[\Hc_j(n^{\alpha},mn)\right]}{\Pb\left[\Hc_j(r,n)\large{|}\Hc_j(n^{\alpha},n)\right] \cdot \Pb\left[\Hc_j(n^{\alpha},n)\right]}\overset{\eqref{eq:couple-h-half}}{\simeq} \frac{\Pb\left[\Hc_j(n^{\alpha},mn)\right] }{\Pb\left[\Hc_j(n^{\alpha},n)\right]} \,.
\end{equation}
We then use Proposition~\ref{prop:compare-half} to pass from scales $mn,n$ to scales $m^2n,mn$:
\begin{equation}
\label{eq:thm1.1-2}
 \frac{\Pb\left[\Hc_j(n^{\alpha},mn)\right] }{\Pb\left[\Hc_j(n^{\alpha},n)\right]}\simeq \frac{\Pb\left[\Hc_j(mn^{\alpha},m^2n)\right] }{\Pb\left[\Hc_j(mn^{\alpha},mn)\right]} \,.
\end{equation}
Next, similar to the reverse of the first step, we transform from mesoscopic comparison to the ratio $\frac{h_j(r,m^2 n)}{h_j(r,mn)}$ by Proposition~\ref{prop:couple-h-half}. We have
\begin{equation}\label{eq:thm1.1-3}
\frac{\Pb\left[\Hc_j(mn^{\alpha},m^2n)\right] }{\Pb\left[\Hc_j(mn^{\alpha},mn)\right]} \overset{\eqref{eq:couple-h-half}}{\simeq}\frac{\Pb\left[\Hc_j(mn^{\alpha},m^2n)\right] }{\Pb\left[\Hc_j(mn^{\alpha},mn)\right]} \cdot \frac{\Pb\left[\Hc_j(r,m^2n)\large{|}\Hc_j(mn^{\alpha},m^2n)\right] }{\Pb\left[ \Hc_j(r,mn)\large{|}\Hc_j(mn^{\alpha},mn)\right] }  
 = \frac{h_j(r,m^2n)}{h_j(r,mn)} \,.
\end{equation}
The combination of \eqref{eq:thm1.1-1}, \eqref{eq:thm1.1-2} and \eqref{eq:thm1.1-3} gives
\begin{equation}\label{eq:hporpotion}
\frac{h_j(r,mn)}{h_j(r,n)}  \overset{\eqref{eq:thm1.1-1}}{\simeq}  \frac{\Pb\left[\Hc_j(n^{\alpha},mn)\right] }{\Pb\left[\Hc_j(n^{\alpha},n)\right]}\overset{\eqref{eq:thm1.1-2}}{\simeq}\frac{\Pb\left[\Hc_j(mn^{\alpha},m^2n)\right] }{\Pb\left[\Hc_j(mn^{\alpha},mn)\right]} \overset{\eqref{eq:thm1.1-3}}{
\simeq} \frac{h_j(r,m^2n)}{h_j(r,mn)} \,.
\end{equation}

We now turn back to $b_j$. First,
\begin{equation}\label{eq:bhrelation}
\frac{b_j(r,mn)}{b_j(r,n)} = \frac{\Pb\left[\Bc_j(r,mn)\large{|}\Bc_j(r,n^\alpha) \right] }{\Pb\left[\Bc_j(r,n)\large{|}\Bc_j(r,n^\alpha)\right] } \overset{\eqref{eq:couple-bh-half}}{\simeq}  \frac{\Pb\left[\Hc_j(r,mn)\large{|}\Hc_j(r,n^\alpha) \right] }{\Pb\left[\Hc_j(r,n)\large{|}\Hc_j(r,n^\alpha)\right]}  \overset{\eqref{eq:thm1.1-1}}{\simeq}  \frac{\Pb\left[\Hc_j(n^\alpha,mn) \right] }{\Pb\left[\Hc_j(n^\alpha,n) \right]}  \,.
\end{equation}
The next step is the same as \eqref{eq:thm1.1-2}. The third step is similar to the reverse of the first step, but we use \eqref{eq:thm1.1-3} instead of \eqref{eq:thm1.1-1}. As a result, we obtain
\begin{equation}\label{eq:bhrelation2}
\frac{\Pb\left[\Hc_j(mn^{\alpha},m^2n)\right] }{\Pb\left[\Hc_j(mn^{\alpha},mn)\right]} \simeq \frac{b_j(r,m^2n)}{b_j(r,mn)}\,.
\end{equation}
The combination of the above two estimates on proportions as well as \eqref{eq:thm1.1-2} gives
\begin{equation}
\frac{b_j(r,mn)}{b_j(r,n)}  \overset{\eqref{eq:bhrelation}}{\simeq}  \frac{\Pb\left[\Hc_j(n^{\alpha},mn)\right] }{\Pb\left[\Hc_j(n^{\alpha},n)\right]} 
\overset{\eqref{eq:thm1.1-2}}{\simeq}\frac{\Pb\left[\Hc_j(mn^{\alpha},m^2n)\right] }{\Pb\left[\Hc_j(mn^{\alpha},mn)\right]} \overset{\eqref{eq:bhrelation2}}{\simeq} \frac{b_j(r,m^2n)}{b_j(r,mn)}  \,.
\end{equation}
This finishes the proof of Proposition \ref{prop:ratioestimate-half}.
\end{proof}
\begin{proof}[Proof of Theorem~\ref{thm:newes-half}]
For a given $j$, it suffices to consider only $r\geq r_b(j)$.  It follows from Proposition~\ref{prop:ratioestimate-half} and Lemmas~\ref{prior_3} and \ref{lem:sequence} that there exist $0<C<\infty$ and $-\infty< \alpha<\infty$ such that $b_j(r,n)=Cn^{\alpha}\Big(1+O(n^{-c})\Big)$. By Lemma~\ref{prior_1}, $\alpha=-\beta_j$ and so
\begin{equation*}
b_j(r,n)=Cn^{-\beta_j}\big(1+O(n^{-c})\big)\,.
\end{equation*}
And the same holds for $h_j$. With  This finishes the proof of Theorem~\ref{thm:newes-half}.
\end{proof}

\subsection{Comparison estimates in the plane}
\label{subsec:4.3}
Our goal of this subsection is Proposition~\ref{prop:compare-whole} in which we compare the probabilities of $\Yc_j$ at different scales. The proof is similar to that of the half-plane case in spirit, in which we relate $\Yc_j$ to the exploration process and then make use of the power-law rate of convergence from \cite{richards2021convergence}. However when  the number of arms is odd, there are inevitably two neighboring ones of the same color that are not  separated by an interface. This issue poses extra difficulties. See Lemma \ref{lem:msrb-whole} and the paragraph right above it for more details.

\begin{figure}
	\centering
	\includegraphics[height=.4\textwidth]{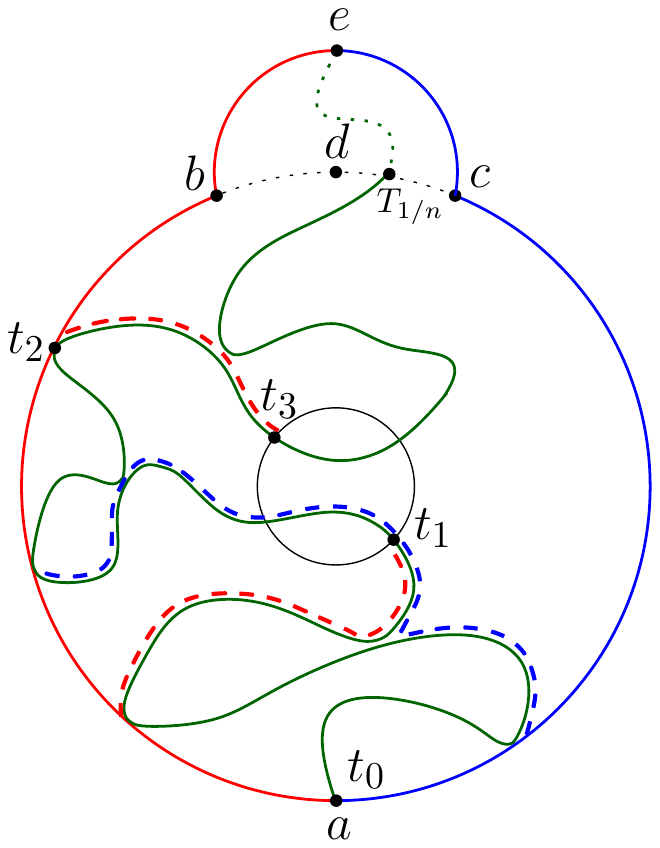}\quad\quad\quad\quad
	\includegraphics[height=.4\textwidth]{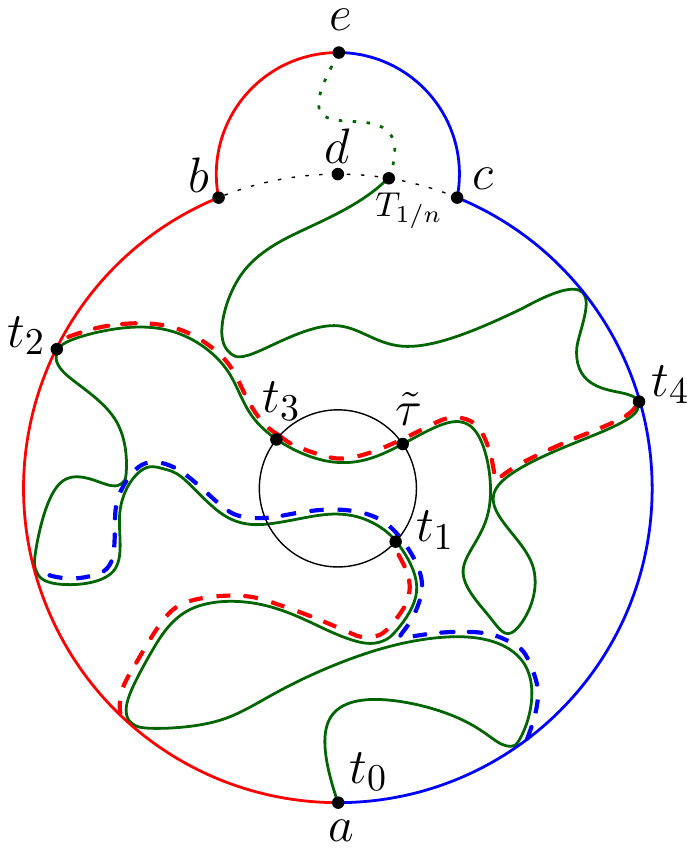}
	\caption{Sketch of events $\Yc_4(n^\alpha,n)$ on the left and  $\Yc_5(n^\alpha,n)$ on the right. The whole domain is $\Omega_{[1/n]}$. 
		The small circle inside this domain is $C_{n^{\alpha-1}}$.
		The exploration process $\gamma_{1/n}$ from $a_{1/n}$ to $e_{1/n}$ on $\Omega_{[1/n]}$ is in green. $T_{1/n}$ is the first time that $\gamma_{1/n}$ enters $U_{[1/n]}$, after which time $\gamma_{1/n}$ is sketched in dotted curves. The arms, neighboring the interface $\gamma_{1/n}$, to enforce the occurrence of arm events are in dashed red and blue. Note that on the right the ``disjointedness condition'' holds: $\gamma_{1/n}(t_2,t_4)$ does not hit $\wideparen{cdb}$, $\tilde\sigma=t_3$, and the left boundaries of $\gamma_{1/n}(t_2,\tilde \sigma)$ and $\gamma_{1/n}(\tilde\tau,t_4)$ are disjoint.}
	\label{figure:explore-whole}
\end{figure}

Fix in this subsection the nice domain (in the sense of Section \ref{subsec:setting}) and four points on its boundary 
$$\Omega= B_1 \cup B\Big((0,1), 2\sin (\pi/12)\Big),\; a=(0,-1),\;b=\left(-\frac{1}{2},\frac{\sqrt{3}}{2}\right),\;c=\left(\frac{1}{2},\frac{\sqrt{3}}{2}\right),\mbox{ and  }d = (0,1) .$$ See Figure \ref{figure:explore-whole} for an illustration. We also let
$$
e=\big(0,1+2\sin (\pi/12)\big) \mbox{ and $U=B\big(d,2 \sin (\pi /2)\big) \backslash B_1 \subset \Omega$ which contains a neighborhood of $e$}.
$$
Note that the value $2 \sin (\pi/12)$ equals to the length of $\overline{bc}$ (and $\overline{cd}$).  Given a mesh size $\eta$, the percolation exploration process starts from $a_\eta$ and ends at $e_\eta$. We assign red color to the counterclockwise arc $\wideparen{eba}$ and blue color to the counterclockwise arc $\wideparen{ace}$. See Figure~\ref{figure:explore-whole} for an illustration of the setup. Recall that $T_{1/n}$ is the first time that $\gamma_{1/n}$ enters $U_n$.  The stopping times $T_{1/mn}$ and $T$ are defined similarly.

Recall \eqref{def:y} that $\Yc_j(r,R)$ is the event that in $C_R$ with mesh size $\eta = 1$ there exist $j$ disjoint arms connecting $C_r$ with $C_R$ with the prescribed color pattern and some additional constraints. Rescaling the lattice by $1/R$, we get that $\Yc_j(r, R)$ is equivalent to the event there exist $j$ such arms connecting $C_{r/R}$ to $C_1$ in $B_1$ with mesh size $\eta = 1/R$. Similar to the half-plane case in Section \ref{subsec:4.1}, we will consider two cases $R=n$ and $mn$ in this subsection and we will tacitly assume the equivalence between events on the rescaled and those on the unrescaled lattices. 

Recall the definition of $y_j(r,R)$ in  \eqref{def:sequence}. The following proposition plays a similar role to that of Proposition~\ref{prop:compare-half} in the half-plane case.  
It states that $y_j$ are almost the same at different scales within a power-law error when $r$ is not too small compared to $R$. 

\begin{proposition}
\label{prop:compare-whole}
There exists $c_9>0$ such that for all $\alpha \in (1-c_9,1)$ and $m \in (1.1,10)$,
$$
y_j(n^{\alpha},n)=y_j(mn^{\alpha},mn)\big(1+O(n^{-c})\big)\,.
$$
Here, $O(n^{-c})$ is independent of $m$ and $\alpha$.
\end{proposition}
We postpone the proof till the end of this subsection and turn to Lemmas~\ref{lem:msrb-whole} and \ref{lemma4.7-1} which play same roles as Lemmas \ref{lemma4.2-1} and \ref{lemma4.2-2} in the half-plane case.  However, the situation here is more complicated, since $\Yc_j(n^\alpha,n)$ corresponds to different events according to whether $j$ is even or odd. When $j$ is odd, in the event $\Yc_j(n^\alpha,n)$ there are two neighboring arms of the same color which lead to the introduction of the ``disjointedness condition'' (see \eqref{eq:behavior-gamma-whole-odd} and \eqref{eq:def-disjoint-cond}).

\begin{lemma}
\label{lem:msrb-whole}
For even $j \geq 2$, one has
\begin{equation}
\label{eq:behavior-gamman-whole}
\Yc_j(n^{\alpha},n)=\{ \gamma_{1/n} \mbox{ hits } C_{n^{\alpha-1}}, \wideparen{ba}, C_{n^{\alpha-1}},  \wideparen{ac},\cdots,\;(j-1) \mbox{ times before } T_{1/n} \} \,,
\end{equation}
and for odd $j \geq 3$, one has
\begin{equation}
\begin{split}
\label{eq:behavior-gamma-whole-odd}
&\Yc_j(n^{\alpha},n)=\Big\{ \gamma_{1/n} \mbox{ hits } C_{n^{\alpha-1}}, \wideparen{ba}, C_{n^{\alpha-1}}, \wideparen{ac},\cdots, \;(j-1) \mbox{ times before } T_{1/n} \mbox{and the last}\\ 
& \quad \quad\quad\quad\quad\quad\mbox{  two crossings satisfy the ``disjointedness condition'' defined in \eqref{eq:def-disjoint-cond}}\Big\} \,.
\end{split}
\end{equation}
\end{lemma}

For $\gamma_{1/n}$ satisfying the requirements on the RHS of \eqref{eq:behavior-gamman-whole} and \eqref{eq:behavior-gamma-whole-odd},
there exists a sequence of times $0=t_0<t_1<t_2<...<t_{j-1} \leq T_{1/n}$ such that 
\begin{equation}
\label{eq:def-condition-whole}
\gamma_{1/n}(t_1) \in C_{n^{\alpha-1}}\,,\quad \gamma_{1/n}(t_2) \in \wideparen{ba}\,,\quad \gamma_{1/n}(t_3) \in C_{n^{\alpha-1}}\,,\quad \gamma_{1/n}(t_4) \in  \wideparen{ac} \,, \quad\mbox{etc.}
\end{equation}In particular, $\gamma_{1/n}$ is allowed to hit other arcs in the middle of two consecutive hitting times. For instance, it can hit $\wideparen{ac}$ between the times hitting $C_{n^{\alpha-1}}$ and $\wideparen{ba}$. We say that this time sequence satisfies the ``disjointedness condition'' if
\begin{equation}
\label{eq:def-disjoint-cond}
\begin{aligned}
&\gamma_{1/n}(t_{j-3}, t_{j-1}) \mbox{ does not hit }\wideparen{cdb}\mbox{, and the left boundaries of }\gamma_{1/n}[t_{j-3}, \tilde \sigma]\mbox{ and }\gamma_{1/n}[\tilde \tau,t_{j-1}] \\
&\mbox{(which are the last two red arms) are disjoint} \,,
\end{aligned}
\end{equation}
where $\tilde \sigma$ denotes the first hitting time of $C_{n^{\alpha-1}}$ after $t_{j-3}$ and $\tilde \tau$ denotes the last exit time from $C_{n^{\alpha-1}}$ before $t_{j-1}$. The proof follows from simple geometrical arguments, so we omit it. One can see the following figure for an illustration.

Similarly, we can show that for even $j \geq 2$, 
\begin{equation}
\label{eq:behavior-gammamn-whole}
\Yc_j(mn^{\alpha},mn)=\{ \gamma_{1/mn} \mbox{ hits } C_{n^{\alpha-1}},  \wideparen{ba}, C_{n^{\alpha-1}},  \wideparen{ac},\cdots,\; (j-1) \mbox{ times before } T_{1/mn} \} \,,
\end{equation}
and the same for odd $j \geq 3$ with the extra ``disjointedness condition''.

Before stating the second lemma, we introduce several events $\Kc_1,\Kc_2, \Kc_3$ (and $\Kc_4$ when $j$ is odd) which all happen with high probability under the good coupling from Proposition \ref{prop:couple-SLE} of $\gamma_{1/n}, \gamma_{1/mn}$ and $\gamma$, whose joint law\footnote{Similar to the half-plane case, here we also consider a probability space in which all three objects are coupled together.} we also denote by $\Pb$. Let $u$ be the constant from Proposition~\ref{prop:couple-SLE}. We will call $\wideparen{ba} \cup \wideparen{ac} \cup\wideparen C_{n^{\alpha-1}}$ the {\bf designated boundary} and note that $b,a,d$ are the extremal points of the designated boundary. Define $\Kc_1$, $\Kc_2$ and $\Kc_3$ as
\begin{align*}
\Kc_1 &= \Big{\{}d(\gamma_{1/n}|_{[0,T_{1/n}]}, \gamma_{1/mn}|_{[0,T_{1/mn}]}) >2n^{-u}\Big{\}}\,;\\
\Kc_2 &= \Big\{ \mbox{there does not exist any point }x\mbox{ on the designated boundary such that it is }n^{-v} \mbox{ away} \\
&\quad \quad \mbox{from the extremal points, and }\gamma_{1/n}\mbox{ enters }B(x,2n^{-u})\mbox{ but does not hit the designated}\\
&\quad \quad \mbox{hit the designated boundary between the last  entrance and first exit times of }B(x,n^{-v}),\\
&\quad\quad \mbox{and the same for }\gamma_{1/mn} \Big\}\,;\\
\Kc_3 &= \Big\{ \gamma_{1/n} \mbox{ and }\gamma_{1/mn}\mbox{ do not enter } B(b, 2n^{-v}),\mbox{ or } B(d, 2n^{-v}) \mbox{, and do not reenter } B(a,2n^{-v})\\
&\quad \quad \mbox{ after leaving }B(a,1/4) \Big\} \,,
\end{align*}
and for odd $j$'s, define $\Kc_4$ as
\begin{align*}
\Kc_4 &= \Big\{\mbox{there does not exist a time sequence that satisfies }\eqref{eq:def-condition-whole} \mbox{ and the} \\
&\quad \quad \mbox{``disjointedness condition'' }\eqref{eq:def-disjoint-cond}\mbox{ but the last two arms are }(2n^{-u})\mbox{-close} \Big\}\,.
\end{align*}
Here, $v \in (0,u)$ is a constant to be chosen later such that \eqref{eq:choose-v-whole} holds. If $\Kc_2$ happens, then each time $\gamma_{1/n}$ (or $\gamma_{1/mn}$) enters $B(x, 2n^{-u})$, it will hit the designated boundary between the last entrance and first exit times of $B(x,n^{-v})$  for all points $x$ on the designated boundary that are $n^{-v}$ away from the extremal points. 

We now give the second lemma which states that the indicator functions of $\Yc_j(n^\alpha,n)$ and $\Yc_j(mn^\alpha,mn)$ are identical under $\cap_{i} \Kc_i$. 
In the proof, we will use the equivalence between arm events and behavior of exploration processes. established in  Lemma~\ref{lem:msrb-whole} and  \eqref{eq:behavior-gammamn-whole}.
\begin{lemma}
\label{lemma4.7-1} One has
\begin{equation}
\label{eq:identity-whole-1}
\mathbbm{1}_{\Yc_j(n^{\alpha},n)} =\mathbbm{1}_{\Yc_j(mn^{\alpha},mn)} \mbox{ on } \cap_{i=1}^3\Kc_i \mbox{ (for $j$ even) or } \cap_{i=1}^4\Kc_i \mbox{ (for $j$ odd)}\,.
\end{equation}
\end{lemma}
\begin{proof}
We begin with the case where $j$ is even. Assume that $\Kc_1,\Kc_2, \Kc_3$ and $\Yc_j(n^\alpha,n)$ happen and so $\gamma_{1/n}$ fulfills the RHS of \eqref{eq:behavior-gamman-whole}. Suppose that $\gamma_{1/n}$ hits $C_{n^{\alpha-1}}$ at the point $x$. By $\Kc_1$ and $\Kc_2$, the process $\gamma_{1/mn}$ enters $B(x,2n^{-u})$ and so must hit $C_{n^{\alpha-1}}$ between the last entrance and first exit times of $B(x,n^{-v})$. Thus, $\gamma_{1/mn}$ hits $C_{n^{\alpha-1}}$ at the same place. Suppose that $\gamma_{1/n}$ hits $\wideparen{ba}$ at the point $x$. From $\Kc_3$ we know that $x$ is $n^{-v}$ away from the extremal points. By $\Kc_1$ and $\Kc_2$, the process $\gamma_{1/mn}$ enters $B(x,2n^{-u})$ and so hits $\wideparen{ba}$ between the last entrance and first exit times of $B(x,n^{-v})$, therefore also before $T_{1/mn}$. The same holds for $\wideparen{ac}$. Therefore, $\gamma_{1/mn}$ hits each side at the same place as $\gamma_{1/n}$ before $T_{1/mn}$ and so $\gamma_{1/mn}$ fulfills the RHS of \eqref{eq:behavior-gammamn-whole} which implies that $\Yc_j(mn^{\alpha},mn)$ happens. Therefore,
$$
\Kc_1 \cap \Kc_2 \cap \Kc_3 \cap \Yc_j(n^{\alpha},n) \subset \Yc_j(mn^{\alpha},mn)\,.
$$In the same way, we can show that $$\Kc_1 \cap \Kc_2 \cap \Kc_3 \cap \Yc_j(mn^{\alpha},mn) \subset \Yc_j(n^{\alpha},n)\,.$$This completes the proof of \eqref{eq:identity-whole-1}.

For odd $j$, we add one more constraint $\Kc_4$ to ensure that we can still find a time sequence that satisfies the ``disjointedness condition'' even after minor perturbation of the process.  Then, it is easy to verify that $\mathbbm{1}_{\Yc_j(n^\alpha,n)}=\mathbbm{1}_{\Yc_j(mn^\alpha,mn)}$ on $\Kc_1 \cap \Kc_2 \cap \Kc_3 \cap \Kc_4 $, since after a minor perturbation of $2n^{-u}$, the same time sequence (which may differ by a distance of $n^{-v}$) still satisfies \eqref{eq:def-condition-whole} and the ``disjointedness condition''. 
\end{proof}
We are now ready to prove the main result of this subsection.
\begin{proof}[Proof of Proposition~\ref{prop:compare-whole}]
First, we give upper bounds to $\Pb[\Kc_1^c]$, $\Pb[\Kc_2^c]$, $\Pb[\Kc_3^c]$ and $\Pb[\Kc_4^c]$. By Proposition~\ref{prop:couple-SLE}, there exists $u>0$ such that under the good coupling of $\gamma_{1/n}, \gamma_{1/mn}$ and $\gamma$
\begin{equation}
\label{eq:couple-SLE-whole}
\Pb[\Kc_1^c] = \Pb\left[d\left(\gamma_{1/n}|_{[0,T_{1/n}]}, \gamma_{1/mn}|_{[0,T_{1/mn}]}\right) >2n^{-u}\right] =O(n^{-u}) \,.
\end{equation}
If $\Kc_2^c$ happens, we have a half-plane $3$-arm event from a $(2n^{-u})$-ball on the designated boundary to distance $n^{-v}$, see the proof of \eqref{eq:proof-3arm} for more details. Thus, $\Pb[\Kc_2^c] < O(n^{u}) \times O(n^{(1+c)(v-u)})$ for some $c>0$. If $\Kc_3^c$ happens, we have a planar $1$-arm event from a $(2n^{-v})$-ball centered at one of the extremal points to distance $1/4$. So, $\Pb[\Kc_3^c] < O(1) \times O(n^{-cv})$ (where we reduce $c$ if necessary). 

Therefore, for even $j$, under the good coupling
\begin{equation}
\label{eq:prop4.5-11}
\begin{split}
|y_j(n^{\alpha},n)-y_j(mn^{\alpha},mn)| &\overset{\eqref{eq:identity-whole-1}}{\leq} \Pb(\Kc_1^c) +\Pb(\Kc_2^c) +\Pb(\Kc_3^c) \\
&\;\;=\;\; O(n^{-u}) +O(n^{u}) \times O\big(n^{(1+c)(v-u)}\big) + O(n^{-c v}) \,.
\end{split}
\end{equation}

Now, we give upper bound to $\Pb[\Kc_4^c]$. We will call $C_1 \cup C_{n^{\alpha-1}}$ the extended boundary (which is different from the designated boundary we defined before). If $\Kc_4^c$ happens, we have a whole-plane $6$-arm event from a $(2n^{-u})$-ball in $B_1$ to distance $n^{-v}$, or on the extended boundary a half-plane $4$-arm event from distance $3n^{-v}$ to $1/4$. We assume that $\Kc_4^c$ happens for $\gamma_{1/n}$ and $j \equiv 1 \pmod{4}$. The cases for $\gamma_{1/mn}$ or $j \equiv 3 \pmod{4}$ can be treated similarly. 

If $\Kc_4^c$ happens, then the last two red arms, i.e., the left boundaries of $\gamma_{1/n}[t_{j-3}, \tilde \sigma]$ and $\gamma_{1/n}[\tilde \tau, t_{j-1}]$ are both $(2n^{-u})$-close to a point $x$. If $x$ is $n^{-v}$ away from the extended boundary, we write $\sigma^1$ (resp.\ $\sigma^3$) for the first hitting time of $C(x,2n^{-u})$ after $t_{j-3}$ (resp.\ $\tilde \tau$) and $\tau^1$ (resp.\ $\tau^3$) for the last exit time of $C(x,n^{-v})$ before $\sigma^1$ (resp.\ $\sigma^3$). Write $\sigma^2$ (resp.\ $\sigma^4$) for the first exit time of $C(x,n^{-v})$ after $\sigma^1$ (resp.\ $\sigma^3$) and $\tau^2$ (resp.\ $\tau^4$) for the last exit time of $C(x,2n^{-u})$ before $\sigma^2$ (resp.\ $\sigma^4$). Then, $$\tau^1 < \sigma^1 \leq \tau^2 <\sigma^2 <\tau^3 <\sigma^3 \leq \tau^4 <\sigma^4.$$ Furthermore, the left boundaries of $\gamma_{1/n}[\tau^1,\sigma^1], \gamma_{1/n}[\tau^2, \sigma^2],  \gamma_{1/n}[\tau^3, \sigma^3],  \gamma_{1/n}[\tau^4, \sigma^4]$ are four disjoint red arms from $B(x,2n^{-u})$ to distance $n^{-v}$, and the right boundaries of $\gamma_{1/n}[\tau^1,\sigma^1], \gamma_{1/n}[\tau^3, \sigma^3]$ are two disjoint blue arms from $B(x,2n^{-u})$ to distance $n^{-v}$. So, there is a whole-plane $6$-arm event from a $(2n^{-u})$-ball in $B_1$ to distance $n^{-v}$. If $x$ is $n^{-v}$-close to a point $y$ on the extended boundary, then the left and right boundaries of $ \gamma_{1/n}[\tau^2, \sigma^2], \gamma_{1/n}[\tau^3, \sigma^3]$ (where we replace $2n^{-u}$ in the definition to $2n^{-v}$ and $n^{-v}$ to $1/4$) are four disjoint crossings from $B(y,3n^{-v})$ to distance $1/4$. We can cover $B_1$ by a $(2n^{-u})$-net with $O(n^{2u})$ elements and cover the extended boundary by a $(2n^{-v})$-net with $O(n^v)$ elements. Thus, 
$$\Pb[\Kc_4^c] =O(n^{-cv}) + O\Big(n^{2u+(2+c)(v-u)}\Big),\mbox{ where we reduce $c$ if necessary.}$$

Therefore, for odd $j \geq 3$, under the good coupling
\begin{equation*}
\begin{split}
&|y_j(n^{\alpha},n)-y_j(mn^{\alpha},mn)| \leq \Pb[\Kc_1^c] +\Pb[\Kc_2^c] +\Pb[\Kc_3^c]+ \Pb[\Kc_4^c] \\
\leq\;& O(n^{-u}) +O(n^{u}) \times O\big(n^{(1+c)(v-u)}\big) + O(n^{-c v})  + O\Big(n^{2u+(2+c)(v-u)}\Big) \,.
\end{split}
\end{equation*}

Take a small $v$ such that 
\begin{equation}
\label{eq:choose-v-whole}
u+(1+c)(v-u) <0\quad\mbox{ and } \quad 2u+(2+c)(v-u)<0\,.
\end{equation}By applications of RSW theory, $y_j(mn^\alpha,mn) > n^{C (\alpha-1)}$. We can take $c_9$ such that $$-C  c_9>\max \{ -u,-c v, u+(1+c)(v-u), 2u+(2+c)(v-u) \}\,.$$ Then, for all $\alpha \in (1-c_9,1)$, $|y_j(n^\alpha,n) -y_j(mn^\alpha,mn)|=O(n^{-c}) y_j(mn^\alpha,mn)$ as desired.
\end{proof}

\subsection{Proof of Theorems~\ref{thm:newes-whole} and \ref{thm:newes2-whole}}
\label{subsec:4.4}
In this subsection, we will complete the proof of Theorems~\ref{thm:newes-whole} and \ref{thm:newes2-whole}.

The proof of Proposition~\ref{prop:ratioestimate-whole} is the same as that of Proposition~\ref{prop:ratioestimate-half} except that we will replace Proposition~\ref{prop:couple-h-half} with Proposition~\ref{prop:couple-y-whole} and Proposition~\ref{prop:compare-half} with Proposition~\ref{prop:compare-whole}. We are now ready to prove Theorem~\ref{thm:newes-whole}.

\begin{proof}[Proof of Theorem~\ref{thm:newes-whole}] We start with $y_j$. For a given $j$, it suffices to consider only $r\geq r_y(j)$.
Combining with Proposition~\ref{prop:ratioestimate-whole} and Lemmas~\ref{prior_3}, \ref{lem:sequence}, there exist $0<C<\infty$ and $\alpha$ such that $y_j(r,n)=Cn^{\alpha}\big(1+O(n^{-c})\big)$. By Lemma~\ref{prior_1}, $\alpha=-\alpha_j$.

We now turn to $x_j$. By applying Proposition~\ref{prop:couple-a-whole} we obtain that
\begin{equation}\label{eq:xycomparison}
\frac{x_j(r,mn)}{x_j(r,n)}=\frac{y_j(r,mn)}{y_j(r,n)}\Big(1+O(n^{-c})\Big)
\end{equation}
for $r,m$ satisfying the requirement of Proposition~\ref{prop:ratioestimate-whole}. We thus establish a proportion estimate similar to \eqref{eq:ratioestimate-whole} for $x_j$, which yields \eqref{eq:newes-whole} for $x_j$.
\end{proof}

We now turn to Theorem~\ref{thm:newes2-whole}.
\begin{proof}[Proof of Theorem~\ref{thm:newes2-whole}]
By Proposition~\ref{prop:couple-A-whole-inner},
\begin{equation}\label{eq:norateprop}
\frac{a_j(1,n)}{x_j(1,n)}\overset{\eqref{eq:couple-A-whole-inner}}{=}\frac{a_j(\epsilon n,n)}{x_j(\epsilon n,n)}\Big(1+O(\epsilon^c)\Big).
\end{equation}
Thanks to Claim (4) of Lemma~\ref{prior_1}, the first term of the RHS of \eqref{eq:norateprop} converges to $g_j(\epsilon)/f_j(\epsilon)$ as $n\to\infty$. Hence the LHS of \eqref{eq:norateprop} must also converge to a constant as $n\to\infty$. Combined with \eqref{eq:newes-whole} for $x_j$, this finishes the proof.
\end{proof}
\begin{remark}\label{rmk:a6}  Note that although $a_6(1,n)>0$, we did not include the case $j=6$ in the statement of Theorem \ref{thm:newes2-whole}, because $y_6(1,n)=0$ by definition and hence \eqref{eq:xycomparison} no longer holds. In this case (plus some other inner initial configurations), obtaining sharp asymptotics for $x_6(1,n)$ as well as for $a_6(1,n)$ requires some extra coupling argument, which we omit for brevity.
\end{remark}

	\section{Proof of coupling results}\label{sec:proofcoupling}
In this section, we give the proof of the couplings in the half-plane and plane cases stated in Propositions~\ref{prop:coup-1} and \ref{prop:C-in-coupling-even} respectively. As the idea for all these couplings are quite similar, we will give a detailed proof for Proposition \ref{prop:coup-1} in Section \ref{subsec:5.1} but only sketch the differences in the details of the proof for Proposition \ref{prop:C-in-coupling-even} in Section \ref{subsec:5.2}.

Before going into details for each specific setup, we briefly explain here the core idea of these couplings. In essence, one divides the domain into exponential scales (``layers'' in the text) and couple configurations of each layer step by step. In each step, when passing from one scale to the next, one can show by the separation lemma and RSW-FKG gluing technique that (although under different types of conditioning) the laws of nice configurations in previous scales are absolutely continuous with respect to each other with bounded Radon-Nykodym derivative, and hence different conditional laws can be coupled with positive probability at each scale. Thus the coupling will succeed in the end with high probability. It can also be proved that with high probability, some good events will happen which encompass all the dependence of the past and future configurations under the conditional laws. The result we want follows from the combination of these two facts.

\subsection{$j\ge2$ arms in the half-plane}\label{subsec:5.1}
In this subsection, we give the proof of Proposition~\ref{prop:coup-1}. We divide the (discretized) half-plane $\mathbb{H}$ into the disjoint union of dyadic semi-annuli\footnote{For simplicity in this subsection we do not add the $+$ in the superscript.} $A_i=B_{r_{i+1}}^+\setminus B_{r_i}^+$, where $r_{i}=2^i$ for $i=0,1,2,\cdots$. Call the $A_i$'s the \textbf{layers} in $\Hb$. We will do the coupling layer by layer.

In the following, we fix some annulus region $A_i\cup A_{i+1}$ between $C_r^+$ and $C_R^+$ and consider the percolation configuration inside this annulus. We start by defining the \textbf{good event} in $A_i\cup A_{i+1}$. 
\begin{definition}[Good event]\label{def:good-event-i}
	Define the good event $\Gc_j^{(i)}$ associated with the percolation configuration inside $A_i\cup A_{i+1}$ as follows: 
	\begin{itemize}
		\item There are exactly $(j-1)$ interfaces crossing the annulus $A_i\cup A_{i+1}$. Denote them by $\gamma_1,\dots,\gamma_{j-1}$ in counterclockwise order, then $\gamma_k$ is adjacent to $\gamma_{k+1}$ for any $1\le k\le j-2$. In addition, $\gamma_1,\dots,\gamma_{j-1}$ are well-separated on both ends.
		\item There is a red path connecting $[r_i,2r_i]$ to $\gamma_1$ in $A_i$ and a red path connecting $[2r_i,4r_i]$ to $\gamma_1$ in $A_{i+1}$.
		\item There are two red (if $j$ is odd) or blue (if $j$ is even) paths which connect $[-2r_i,-r_i]$ and $[-4r_i,-2r_i]$ to $\gamma_{j-1}$ in $A_i$ and $A_{i+1}$, respectively.
	\end{itemize}
\end{definition}

We will need that the good event $\Gc_j^{(i)}$ in $A_i\cup A_{i+1}$ happens with at least some positive probability depending only on $j$. 
\begin{lemma}\label{lem:positive_probability2}
	There is a constant $c(j)>0$ such that for any $r_i\ge 10j$, we have $
	\Pb[\Gc_j^{(i)}]\ge c.
$
\end{lemma}
This lemma could be proved in a way similar to the second proof of \cite[Lemma 3.4]{MR3073882} by constructing pivotal points through the exploring process (to ensure the interface conditions) together with FKG-RSW gluing (to ensure the well-separateness conditions). Here we present the following alternative proof which relies on a combination of FKG-RSW gluing and BK-Reimer's inequality.
\begin{proof}
 Define a {\bf quasi-good} event $\Uc$ in $A_i\cup A_{i+1}$ as follows:  
	\begin{itemize}
		\item There are two configurations of well-separated inner and outer faces, say $\Theta=\{\theta_1,\dots,\theta_j\}$ and $\Theta'=\{\theta_1',\dots,\theta_j'\}$), which are around $C^+_{4r_i}$ and $C_{r_i}^+$ resp., and lying completely in $\Hb\setminus B_{3r_i}^+$ and in $B_{2r_i}^+$ resp.
		\item For each $1\le k\le j$, both $\theta_k$ and $\theta_{k}'$ are red (if $k$ is odd) or blue (if $k$ is even). Furthermore, there is a path connecting $\theta_k,\theta_{k}'$ together with the aforementioned color.
	\end{itemize}
From FKG-RSW gluing technique we see $\Pb[\Uc]\ge c_0$ for some positive constant $c_0(j)$. If $\Uc$ happens, there will be exactly $(j-1)$ interfaces crossing $A_i\cup A_{i+1}$ which end at the $(2j-2)$ endpoints of $\Theta,\Theta'$, respectively. Write
$$\Wc:= \Big\{\mbox{each interface is adjacent to both of its neighbors}\Big\}.
$$ By the observation that $\Uc\cap \Wc\subset \Gc_j^{(i)}$ and $\Uc\cap \Wc^c\subset\Uc\square\Bc_1(2r_i,3r_i)$ (Recall the notation $\square$ from Section \ref{subsec:basictools}), we conclude that (by applying BK-Reimer's inequality (Lemma \ref{lem:Reimer}) in the last step)
\begin{equation}
	\Pb[\Gc_j^{(i)}]\ge \Pb[\Uc\cap \Wc]\ge \Pb[\Uc]-\Pb[\Uc\square \Bc_1(2r_i,3r_i)]\ge \Pb[\Uc](1-\Pb[\Bc_1(2r_i,3r_i)])
\end{equation}
which is bounded below by some $c(j)>0$.
\end{proof}

\begin{figure}[h!]
	\centering
	\includegraphics[width=.5\textwidth]{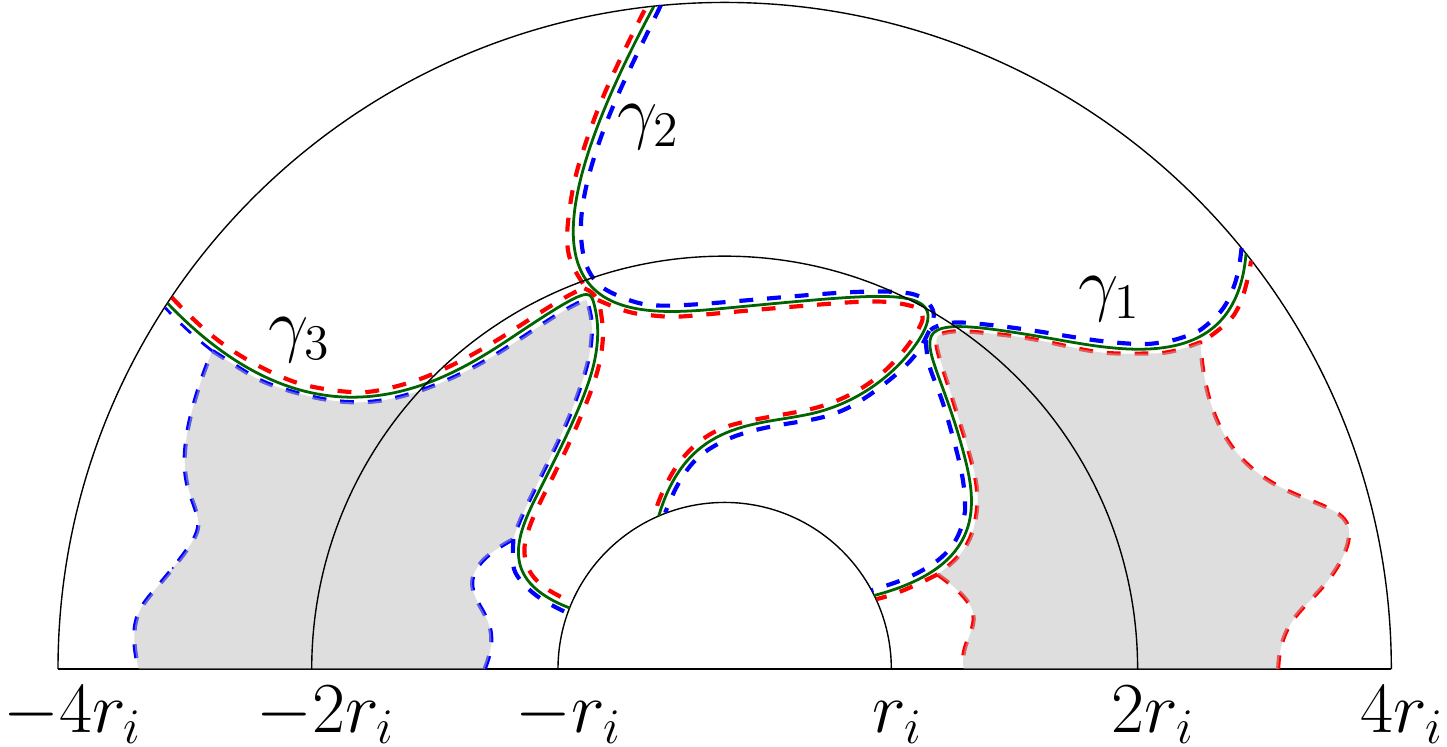}
	\caption{The good set $\Sc$ when $\Gc_j^{(i)}$ holds with $j=4$. The three interfaces $\gamma_i,i=1,2,3$ crossing $A^+(r_i,4r_i)$ are in green. All the hexagons that are adjacent to these interfaces are sketched in dashed red or blue curves, and they make up part of the good set $\Sc$. The gray regions on both sides are bounded by the interfaces $\gamma_1$ and $\gamma_3$ with some innermost paths (crossings) in $A_{i+1}$ and outermost ones in $A_{i}$. The (colored) hexagons in these gray regions make up the rest part of $\Sc$.}
	\label{fig:The_Event_Uc}
\end{figure}

\begin{definition}[Good set]
	If $\Gc_j^{(i)}$ holds, let $\Sc$ be the union of the following colored hexagons (we keep the notation used in the definition of $\Gc_j^{(i)}$):
	\begin{itemize}
		\item The hexagons that touch at least one of the $(j-1)$ interfaces;
		\item The hexagons in the region enclosed by four curves, namely the innermost red path in $A_{i+1}$ from $[2r_i,4r_i]$ to $\gamma_1$, the outermost red path in $A_i$ from $[r_i,2r_i]$ to $\gamma_1$, $\gamma_1$, and the real axis;
		\item The hexagons in the region enclosed by four curves, namely the innermost red (if $j$ is odd) or blue (if $j$ is even) path in $A_{i+1}$ from $[-4r_i,-2r_i]$ to $\gamma_j$, the outermost red (if $j$ is odd) or blue (if $j$ is even) path in $A_i$ from $[-2r_i,-r_i]$ to $\gamma_{j-1}$, $\gamma_{j-1}$, and the real axis.
	\end{itemize}
	 If $\Gc_j^{(i)}$ fails, set $\Sc=\emptyset$. For a nonempty set $S$ of hexagons in $A_i\cup A_{i+1}$, we say that $S$ is a \textbf{good set} if $S$ is a possible value for $\Sc$ such that $\Gc_j^{(i)}$ holds.
	See Figure~\ref{fig:The_Event_Uc} for an illustration.
\end{definition}
Fix a large integer $K_1(j)$ such that\footnote{Note that this is possible since $\beta_j \asymp j^2$; same for the definition of $K_2$ and $K_4$ later in this section, noting that $a_j\asymp j^2$ too.}
\begin{equation}
	\label{eq:choice-of-K1}
p(K_1):=\sup_{i\in \mathbb{N}}\Pb[\Bc_{K_1}(r_i,r_{i+1})]<2^{-10-8\beta_j}\,.
\end{equation}
The following proposition is a generalization of (3.1) and (3.2) of \cite{MR3073882}.
\begin{lemma}\label{lem:estimation2}
	There are constants $C,C'>0$ depending only on $j$, such that for any $2^{K_1+1}r\le r_i<4r_i\le2^{-K_1-1}R$, the following holds: for any good set $S\subset A_i\cup A_{i+1}$, any configuration of outer faces $\Theta$ around $C_{r_{i+K_1+3}}^+$ with no more than $K_1$ faces, and any color configuration $\omega_0$ coincides with $\Theta$ and satisfies $\Pb[\Hc_j(r,R)\mid \omega_{\Dc_{\Theta}}=\omega_0]>0$, 
		\begin{equation}\label{eq:S=S-1}
			C^{-1} 2^{-|S|} \leq \mathbb{P}\left[\mathcal{S}=S\mid \Hc_j(r,R),\omega_{\Dc_\Theta}=\omega_0\right] \leq C 2^{-|S|}\,.
		\end{equation}
		Furthermore,
		\begin{equation}\label{eq:G-i-j-1}
			\mathbb{P}\left[\Gc_j^{(i)}\mid \Hc_j(r,R),\omega_{\Dc_{\Theta}}=\omega_0\right] \geq C'\,.
		\end{equation}
		And similar bounds also hold if we replace $\Hc_j(r,R)$ by $\Bc_j(r,R)$. In a similar fashion, same results also hold for inner faces $\Theta$ around $C_{r_{i-K_1-1}}^+$ satisfying the same requirements.
\end{lemma}
\begin{proof}
	We give the proof of \eqref{eq:S=S-1} and \eqref{eq:G-i-j-1}, and the rest are similar. Rewrite the probability in \eqref{eq:S=S-1} by Bayesian formula as
	\begin{equation}\label{eq:S-1}
		\Pb[\Sc=S\mid \Hc_j(r,R),\omega_{\Dc_{\Theta}}=\omega_0]=2^{-|S|}\times \frac{\Pb[\Hc_j(r,R)\mid \mathcal{S}=S,\omega_{\Dc_{\Theta}}=\omega_0]}{\Pb[\Hc_j(r,R)\mid\omega_{ \Dc_{\Theta}}=\omega_0]}.
	\end{equation} 
	Let $\Theta_1$ (resp. $\Theta_2$) be the configuration of inner faces around $C_{4r_i}^+$ (resp.\ the configuration of outer faces around $C_{r_i}^+$) induced by $S$. Denote $\Rc$ for the event that there are $j$ arms connecting $C_R^+$ to each face of $\Theta_1$, and $\Sc$ for the event that there are $j$ arms connecting $[-r,r]$ to each face of $\Theta_2$ (each arm has the same color with the face it connects to). Then we see that the event in numerator of \eqref{eq:S-1} is equivalent to $\Rc \cap \Sc $ conditioned on $\{\omega_{\Dc_{\Theta_1}}=\omega_0\}$. By independence it has probability $\Pb[\Rc\mid \omega_{\Dc_{\Theta_1}}=\omega_0]\cdot \Pb[\Sc]$.
	
	Since $\Theta_2$ is well-separated, by the separation lemma (Lemma \ref{lem:H-sep}) and FKG-RSW gluing, 
	\begin{equation}\label{eq:S-2}
		\Pb[\Sc]\asymp \Pb[\Hc_j(r,r_{i-1})]\,.
	\end{equation} 
In addition, by \eqref{eq:Q-in-2} in Proposition~\ref{prop:half-sep-2} and FKG-RSW gluing, we also have 
	\begin{equation}\label{eq:S-3}
		\Pb[\Rc\mid \omega_{\Dc_{\Theta_1}}=\omega_0]\asymp \Pb[\Bc_j(r_{i+3},R)\mid \omega_{\Dc_{\Theta_1}}=\omega_0]\,.
	\end{equation}
	Finally, similar to the proof of Proposition~\ref{prop:half-sep-2}, we can show the following quasi-multiplicativity
	\begin{equation}\label{eq:S-4}
		\Pb[\Bc_j(r_{i+3},R)\mid\omega_{ \Dc_{\Theta_1}}=\omega_0]\cdot \Pb[\Hc_j(r,r_{i-1})]\asymp \Pb[\Hc_j(r,R)\mid\omega_{ \Dc_{\Theta_1}}=\omega_0]\,.
	\end{equation} 
	Combining \eqref{eq:S-1} and \eqref{eq:S-2}, \eqref{eq:S-3} and \eqref{eq:S-4} we see that  \eqref{eq:S=S-1}. Summing over all good set $S\subset A_i\cup A_{i+1}$, and using Lemma~\ref{lem:positive_probability2}, we get \eqref{eq:G-i-j-1} as desired.
\end{proof}

\begin{proof}[Proof of Proposition~\ref{prop:coup-1}]
	We focus on the first claim in Proposition~\ref{prop:coup-1}. Recall that $u=\sqrt{rR}$.
	Let $i_0<i_N$ be the integers such that $r_{i_0-1}<u\leq r_{i_0}<r_{i_N}\leq R<r_{i_N+1}$. We inductively couple the conditional laws $\Pb_1=\Pb[\;\cdot \mid \Hc_j(r,R)]$ and $\Pb_2=\Pb[\;\cdot\mid \Hc_j(r,mR)]$ from outside to inside layer by layer. The virtue of the coupling is same as the maximal coupling for a Markov chain on a finite state space.
	
	For initialization, sample the color configuration $\omega_1,\omega_2$ outside $C_{i_N}^+$ according to $\Pb_1$ and $\Pb_2$ independently, and set an index $I=0$. Assume now we have sampled two color configurations $\omega_1$ and $\omega_2$ outside $C^+_{r_{i+K_1+1}}$ for some $i\ge i_0$ according to some coupling of these two conditional distribution, we proceed as follow:
	
	If $I\geq 0$, then for both configurations, we independently explore all the interfaces crossing $A^+(r_{i+K_1},r_{i+K_1+1})$ from outside to inside and stop exploring when reaching $\Rb\cup C^+_{r_{i+K_1}}$. This exploring process would induce two configurations of outer faces around $C^+_{r_{i+K_1}}$ for both $\omega_i,i=1,2$, denoted by $\hat{\Theta}_i,i=1,2$. Also note that $\Vc_{\hat{\Theta}_1}$ and $\Vc_{\hat{\Theta}_2}$ are left unexplored. If either $\hat{\Theta}_1$ or $\hat{\Theta}_2$ has more than $K_1$ faces, we just keep $I$ unchanged, explore all the hexagons outside $C_{r_{i+K_1}}^+$ and proceed the previous procedure on the scale $C_{r_{i+K_1}}$. Otherwise, we add $I$ by $1$, and by applying Lemma~\ref{lem:estimation2} for both $\Hc_j(r,R)$ and $\Hc_j(r,mR)$, we can construct a coupling $\Qb_i$ of the laws $\Pb_1[\;\cdot \mid\omega_{\Dc_{\hat{\Theta}_1}}=\omega_1]$ and $\Pb_2[\;\cdot \mid \omega_{\Dc_{\hat{\Theta}_2}}=\omega_2]$ with the following property: 
	\begin{equation}
		\Qb_i\left[\Sc(\omega_1)=\Sc(\omega_2)\neq \emptyset\right]>C'',\ \mbox{For some }C''(j)>0\,,
	\end{equation}
	where $\Sc(\omega_1),\Sc(\omega_2)$ are the sets induced by $\omega_1,\omega_2$ in $A_i\cup A_{i+1}$, respectively.

	We sample the pair of color configurations $(\omega_1,\omega_2)$ outside $C_{r_i}^+$ from $\Qb_i$, and perform the following exploring process to detect whether $\{\Sc(\omega_1)=\Sc(\omega_2)\neq \emptyset\}$ holds: 
	\begin{itemize}
		\item [\bf{Step 1}] For each of the endpoints $x$ of $\hat{\Theta}_1$, we explore the interface of $\omega_1$ inside $B^+_{r_{i+K_1}}$ starting from $x$, and stop exploring until it hits $\Rb\cup C^+_{r_i}$. Denote $\Gamma$ for the subset of the aforementioned interfaces reaching $C_{r_i}^+$. Then we can check whether $|\Gamma|=j-1$ and each interface in $\Gamma$ is adjacent to its neighbors in $A^+(r_i,4r_i)$. 
		\item [\bf{Step 2}] If the condition in \textbf{Step 1} holds, denote $\gamma_1$ and $\gamma_{j-1}$ for the leftmost and rightmost interface in $\Gamma$. In the quad $Q_1$ enclosed by $C^+_{4r_i},\gamma_1,C^+_{r_i}$ and $\Rb$, we start from $Q_1\cap C^+_{2r_i}$ and explore in $\omega_1$ to find the innermost red path in $A_{i+1}$ connecting $[2r_i,4r_i]$ to $\gamma_1$, and the outermost red path in $A_i$ connecting $[r_i,2r_i]$ to $\gamma_1$. Our exploring process stops whenever we find such paths, or if we cannot find them, we stop exploring until reaching $C_{r_i}^+$ or $C_{4r_i}^+$. Perform the same exploration on the left side of $A_i\cup A_{i+1}$. Up to now, we already have enough information to fix $\Sc(\omega_1)$ and to assert whether the good event $\Gc_j^{(i)}$ happens or not for $\omega_1$.
		\item [\bf{Step 3}] Run the same exploring process for $\omega_2$.  Now we can check whether $\Gc_j^{(i)}$ happens for both $\omega_1$ and $\omega_2$  and $\Sc(\omega_1)=\Sc(\omega_2)$. If so, we set $I=-1$; otherwise we keep $I$ unchanged.
	\end{itemize} 
	
	If $I$ is still non-negative, we explore the entire color configuration $\omega_1,\omega_2$ outside $C_{r_i}^+$, and proceed the procedure above on the scale $r_i$. Otherwise $I=-1$, then the identical nonempty good set $\Sc(\omega_1)=\Sc(\omega_2)$ induces a configuration of outer faces $\Theta^*$ with $j$ faces around $C^+_{r_i}$, which is  common for both $\omega_i,i=1,2$. Note that our exploring process ensures that none of the hexagons in $\Vc_{\Theta^*}$ are explored, i.e.\ the law in $\Vc_{\Theta^*}$ is still the critical Bernoulli percolation, and thus by the domain Markov property, $\Pb[\;\cdot\mid \Hc_j(r,R),\omega_{\Dc_{\Theta^*}}=\omega_1]$ and $\Pb[\;\cdot\mid \Yc_j(r,mR),\omega_{\Dc_{\Theta^*}}=\omega_2]$ are precisely equal. Hence, we can couple the color configuration in $\Vc_{\Theta^*}$ of these two distributions identically. In this case we have $(\omega_1,\omega_2)$ satisfies the requirement in proposition~\ref{prop:coup-1}, and $\Theta^*$ is a stopping set as desired.
	
	Finally, we control the probability that our coupling remains unsuccessful until scale $r_{i_0}$. Write $N_1=\frac{N-2K_1}{2K_1+2}$, then $N_1\ge \frac{1}{4K_1+4}\log_2\frac{R}{r}-2$ since $N\ge \frac{1}{2}\log_2\frac{R}{r}-1$. On the one hand, note that in each time $I$ is plus by $1$, we always have a positive probability $C''$ to couple successfully, no matter how things went in previous layers. Thus $$\Pb\left[I\ge N_1\right]\le(1-C'')^{N_1}\le\left({r}/{R}\right)^\delta\mbox{ for some  }\delta(j)>0\,.$$
	On the other hand, $I<N_1$ implies there are more than $N-(K_1+1)N_1-K_1=N/2$ integers $i\in[i_0.i_N-1]$ satisfying $\mathbf{1}_{\Bc_{K_1}(r_i,r_{i+1})}(\omega_1)+\mathbf{1}_{\Bc_{K_1}(r_i,r_{i+1})}(\omega_2)\ge 1$, so 
	$$
	\Tc:=\Big\{\sum_{i=i_0}^{i_N-1}\mathbf{1}_{\Bc_{K_1}(r_i,r_{i+1})}\ge N/4\Big\}$$ happens for either both $\omega_i,i=1,2$. Combining
	$$\Pb_1[\Tc]\le 2^Np(K_1)^{N/4}/\Pb[\Hc_j(r,R)], \quad\quad\Pb_2[\Tc]\le 2^Np(K_1)^{N/4}/\Pb[\Hc_j(r,mR)],$$ with the a priori estimates of arm probabilities, we obtain
	\begin{align*}
		\Pb\left[0\le I\le N_1\right] 
		\le&\ 2^{N}p(K_1)^{N/4}\left[\left({R}/{r}\right)^{-\beta_j+o(1)}+\left({mR}/{r}\right)^{-\beta_j+o(1)}\right],
	\end{align*}
	which is bounded by $(r/R)^{\delta}$ for some $\delta(j)>0$ from the choice of $K_1$ in \eqref{eq:choice-of-K1}. Altogether we conclude $\Pb[I\ge 0]\le \big(r/R\big)^\delta$ for some $\delta(j)>0$, as desired.
	
	With these tools in hand, we can rerun the coupling from inside to outside to conclude the proof of the second claim of Proposition~\ref{prop:coup-1}.
\end{proof}

\subsection{Coupling arm events in the plane}\label{subsec:5.2}
In this subsection, we prove  Proposition~\ref{prop:C-in-coupling-even}. Recall the definitions of events $\Xc_j(r,R),\Yc_j(r,R)$ and $\Ac_j(r,R)$. We shall couple $\Pb[\;\cdot\mid \Xc_j(r,R)]$ with $\Pb[\;\cdot\mid \Ac_j(r,R)]$ and $\Pb[\;\cdot\mid \Yc_j(r,R)]$ with $\Pb[\;\cdot\mid \Yc_j(r,mR)]$. As discussed at the beginning of this section, we will just sketch the proofs and note necessary modifications of arguments from Section \ref{subsec:5.1}. For clearance, the cases of odd $j$'s and even $j$'s are dealt with separately. 

As before, we divide the plane into the union of dyadic annuli $A_i=B_{r_{i+1}}\setminus B_{r_i}$, where $r_i=2^i$ for $i=0,1,2,\cdots$. Call them the {\bf layers} in the plane. 

\subsubsection{Inward coupling: the odd $j$ case}
We begin with the case for odd $j$. As mentioned in Remark~\ref{rmk:whole-coup}, the major difference between the plane and the half-plane cases lies in \eqref{eq:DMP_3}, which is not obvious from the existence of a single configuration of face $\Theta^*$. But as we have shown in the proof of Proposition~\ref{prop:coup-1}, we can indeed construct a coupling $\Qb$ such that $\omega_1,\omega_2$ are identical on some ``good set'' with high probability under $\Qb$. In this subsection, we define the plane version of good events and good sets, and obtain similar result as before; we shall see the existence of a common good set is enough to give \eqref{eq:DMP_3}.

Now for each $r_i>10j$, define good event in $A_i\cup A_{i+1}$ as follows:
\begin{itemize}
	\item There are exactly $(j-1)$ interfaces crossing the annulus, say $\gamma_1,\dots,\gamma_{j-1}$, in counterclockwise order and they are well-separated on both sides of the annulus. 
	\item for $1\leq k\leq j-2$, $\gamma_k$ is adjacent to $\gamma_{k+1}$, while $\gamma_1$ is \emph{not} adjacent to $\gamma_{j-1}$. 
	\item The hexagons between $\gamma_1$ and $\gamma_{j-1}$ that touch $\gamma_1\cup \gamma_{j-1}$ are all red (if $j\equiv 1\pmod 4$) or blue (if $j\equiv 3\pmod 4$). Further, there are two paths with the aforementioned color connecting $\gamma_1$ and $\gamma_{j-1}$, which lie in $A_i$ and $A_{i+1}$, respectively.
\end{itemize}
And if $\Gc_j^{(i)}$ holds, let $\Sc$ be the union of all hexagons which touch at least one of the interfaces $\gamma_k,1\leq k\leq j-1$, or lie in the quad enclosed by $\gamma_1,\gamma_{j-1}$, the outermost and innermost paths connecting $\gamma_1$ and $\gamma_{j-1}$ in $A_i$ and $A_{i+1}$, respectively. We set $\Sc=\emptyset$ if $\Gc_j^{(i)}$ fails, and for a nonempty set of colored hexagons in $A_i\cup A_{i+1}$, we say that $S$ is a {\bf good set} if $S$ is a possible value of $\Sc$. We have the following estimates similarly as before:
\begin{itemize}
	\item There exists $c(j)>0$, such that for any $r_i>10j$, we have $\Pb[\Gc_j^{(i)}]>c.$
	\item For any fixed large integer $K_2$, there exists $C=C(j,K_2)>0$, such that for any good set $S$ in $A_i\cup A_{i+1}$ with $2^{K_2+1}r\le r_i<r_{i+2}\le 2^{-K_2-1}R$, any configuration of outer faces $\Theta$ around $C_{r_{i+K_2+3}}^+$ with no more than $K_2$ faces and any color configuration $\omega_0$ which satisfies $\Pb[\Ac_j(r,R)\mid \omega_{\Dc_{\Theta}}=\omega_0]>0$, 
	\begin{equation}\label{eq:whole-odd-2}
		C^{-1}2^{-|S|}\le \Pb[\Sc=S\mid \Ac_j(r,R),\omega_{\Dc_{\Theta}}=\omega_0]\le C2^{-|S|}\,.
	\end{equation}
	As a result,
$
		\Pb[\Gc_j^{(i)}\mid \Ac_j(r,R),\omega_{\Dc_{\Theta}}=\omega_0]\ge cC^{-1}.
$
\item In addition, similar results also hold for $\Ac_j(r,R)$ replaced by $\Xc_j(r,R)$ or $\Yc_j(r,R)$.
\end{itemize}

\begin{proof}[Sketch of Proof for Proposition~\ref{prop:C-in-coupling-even}, odd $j$ case]
	Fix a large integer $K_2=K_2(j)$ such that 
$$
	p({K_2}):=\sup_{i\in \mathbb N}\Pb[\Ac_{K_2}(r_i,r_{i+1})]<2^{-10-8\alpha_j}.
$$
The first two items can be established similarly as in Lemma \ref{lem:positive_probability2} and \ref{lem:estimation2}, results for $\Xc_j(r,R)$ and $\Yc_j(r,R)$ can also be obtained from the same manner. From these estimates, following the framework of the proof of Proposition~\ref{prop:coup-1}, we can construct a coupling of $\Pb[\;\cdot\mid \Ac_j(r,R)]$ (resp. $\Pb[\;\cdot\mid \Yc_j(r,R)]$) and $\Pb[\;\cdot\mid \Xc_j(r,R)]$ (resp. $\Pb[\;\cdot\mid \Yc_j(r,mR)]$) such that if we sample $(\omega_1,\omega_2)$ according to such coupling, then with probability at least $\big(1-(r/R)^\delta\big)$ for some $\delta(j)>0$, $\omega_1$ and $\omega_2$ are identical on some good set outside $C_{\sqrt{rR}}$. This good set induces a common configuration of faces in both $\omega_i$, $i=1,2$, which is a stopping set and it also ensures \eqref{eq:DMP_4} (resp. \eqref{eq:DMP_3}), thus completing the proof.
\end{proof}

\subsubsection{Inward coupling: even $j$ case}
Finally, we deal with the case for even $j$. This is essentially the case studied in \cite{MR3073882}; we follow the idea in this paper and define good events $\Gc_j^{(i)}$ in $A_i\cup A_{i+1}$ as 
\begin{itemize}
	\item There are exactly $j$ interfaces crossing the annulus $A_i\cup A_{i+1}$, and they are well-separated on both ends. In addition, each of the interfaces is adjacent with its two neighbors.
\end{itemize}
If $\Gc_j^{(i)}$ holds, let $\Sc$ be the union of all colored hexagons in $A_{i}\cup A_{i+1}$ which touch at least one of the $j$ interfaces, otherwise we set $\Sc=\emptyset$. For a nonempty set $S$ of colored hexagons in $A_i\cup A_{i+1}$, we say that $S$ is a \textbf{good set} if $S$ is a possible value of $\Sc$.

The proof of the first item in Proposition~\ref{prop:C-in-coupling-even} is essentially same as before, since coinciding on some good set implies \eqref{eq:DMP_4}. However, for the second item, in order to guarantee \eqref{eq:DMP_3}, it is not enough to make $\omega_1$ and $\omega_2$ coincide on some good set. Indeed, we need to further specify the connecting pattern between this good set and the outer boundary $C_R$ (or $C_{mR}$). For a color configuration $\omega_0$ sampled from $\Pb[\;\cdot\mid \Yc_j(r,R)]$ (resp. $\Pb[\;\cdot\mid \Yc_j(r,mR)]$), assume that under $\omega_0$, for some $u\in [r,R]$ there is a configuration of outer faces $\Theta$ around $C_u$ with $j$ faces. Denote $\tilde\theta$ for the first red face in $\Theta$ when counting from $(0,-u)$ counterclockwise, and denote $\hat{\theta}$ for the red face in $\Theta$ which is connected to $C_R$ (resp. $C_{mR}$) by the first red arm counting from $(0,-R)$ (resp. $(0,-mR)$) counterclockwise.
Then we have the following estimates:
\begin{itemize}
	\item There exists $c(j)>0$, such that for any $r_i>10j$, we have $\Pb[\Gc_j^{(i)}]>c.$
	\item Fix a large integer $K_3$, there exists $C=C(j)>0$, such that for any good set $S$ in $A_i\cup A_{i+1}$ with $2^{K_3+1}r\le r_i<r_{i+2}<2^{-K_3-1}R$, any configuration of outer faces $\Theta$ around $C^+_{r_{i+K_3+3}}$ with no more than $K_3$ faces and any color configuration $\omega_0$ satisfies $\Pb[ \Yc_j(r,R)\mid \omega_{\Dc_\Theta}=\omega_0]>0$, denoting $\Theta_S$ for the configuration of outer faces around $C_{r_i}$ induced by $S$, then it holds that
	\begin{equation}\label{eq:whole-even-2}
		C^{-1}2^{-|S|}\le \Pb[\Sc=S,\tilde{\theta}_S=\hat{\theta}_S\mid \Yc_j(r,R),\omega_{\Dc_{\Theta}}=\omega_0]\le C2^{-|S|}\,.
	\end{equation}
	As a result, 
$
		\Pb[\Gc_j^{(i)}, \tilde{\theta}_S=\hat{\theta}_S\mid \Yc_j(r,R),\omega_{\Dc_{\Theta}}=\omega_0]\ge cC^{-1}.
$
\end{itemize}
\begin{proof}[Sketch of Proof for Proposition~\ref{prop:C-in-coupling-even}, even $j$ case]
	Fix a large integer $K_3(j)$ such that 
$$
p(K_3):=\sup_{i\in \mathbb N}\Pb[\Ac_{K_3}(r_i,r_{i+1})]\le 2^{-10-8\alpha_j}.
$$ The estimates above can be proved similarly, with the caveat that for \eqref{eq:whole-even-2}, we construct some specific gluing to make sure that $\tilde{\theta}_S$ coincides with $\hat{\theta}_S$.
	
	Note that if two color configurations $\omega_1$ and $\omega_2$ are sampled from $\Pb[\;\cdot\mid \Yc_j(r,R)]$ and $\Pb[\;\cdot\mid \Yc_j(r,mR)]$ respectively, such that $\omega_i,i=1,2$ share a common good set $S$, and $\tilde{\theta}_S=\hat{\theta}_S$ holds for both configurations, then we can conclude that \eqref{eq:DMP_3} is true and $\Theta^*=\Theta_S$ is also a stopping set. Following the same framework, this observation together with the estimates above gives us the desired result. 
\end{proof}

	\appendix

\section{Proof of Lemmas~\ref{prior_1} and \ref{prior_3}}
\label{sec:proof-prior}

\begin{proof}[Proof of Lemma~\ref{prior_1}]
	We begin with the proof of $h_j(n) \asymp b_j(n)$. Since $\Hc_j(r,n) \subset \Bc_j(r,n)$, we have $h_j(n) \leq b_j(n)$. It suffices to prove that $h_j(n) \geq cb_j(n)$ for some constant $c>0$ and all $n$ large enough. Let $r' = 10j \vee r$ and $\Gamma$ be the set of interfaces connecting $C_n^+$ and $C_{r'}^+$ in $A^+(r',n)$. Then,
\begin{equation}
\label{eq:apdx-b-1}
\begin{aligned}
h_j(r,n) & \geq \Pb[\Hc_j(r,n)| Q_{\rm in}(\Gamma) >j^{-1}, \Bc_j(r',n)] \times \Pb[Q_{\rm in}(\Gamma) >j^{-1} |\Bc_j(r',n)] \times \Pb[\Bc_j(r',n)] \\
&\geq c \times c \times  b_j(r',n) \geq cb_j(r,n)\,.
\end{aligned}
\end{equation}
(The second inequality follows from RSW theory and Lemma~\ref{lem:H-sep}. The third inequality follows from the fact that $\Bc_j(r,n) \subset \Bc_j(r',n)$.) Claim (2) is derived in the same way, but using Proposition~\ref{prop:C-sep} instead of Lemma~\ref{lem:H-sep}. 

Now, we consider Claim (3). We can show the up-to-constants equivalence between arm probabilities using similar arguments in \eqref{eq:apdx-b-1}. 
The methods in \cite{MR1879816} for calculating asymptotics of arm probabilities down to microscopic scales (i.e., $b_j(n)$ or $p_j(n)$ in our notation) can also be applied to derive mesoscopic asymptotics.

Claim (4) follows from \cite{smirnov2001critical} and several stability arguments (which need looser estimates than those in the article). One can refer to the proof of Lemma 2.9 in \cite{MR3073882} for more details.
\end{proof}

\begin{proof}[Proof of Lemma~\ref{prior_3}]
	We only prove the case of $h_j$. The cases for other arm probabilities 
can be proved in the same way. Let $\Gamma$ be the set of interfaces connecting $[-r,r]$ and $C_n^+$ in $B_n^+$. For $\epsilon>0$, define $\Hc_j^{\epsilon}(r,n)$ by $$\Hc_j^{\epsilon}(r,n) = \Hc_j(r,n) \cap \{ Q_{\rm ex}(\Gamma) \geq \epsilon \}\,. $$Let $\Rc^{\epsilon}$ denote the event that there exists a point $x$ on $C^+_n$ such that there are three disjoint arms (not all of the same color) connecting $ C(x,\epsilon n)$ and $C(x, n/2)$. By Claim (3) in Lemma~\ref{prior_1}, $\Pb[\Rc^\epsilon] \leq O(1/\epsilon) \times O(\epsilon^{1+c}) = O(\epsilon^c)$. Thanks to the spatial independence of percolation
$$
\Pb[\Hc_j(r,n) \backslash \Hc_j^{\epsilon}(r,n)] \leq \Pb[\Hc_j(r, n/2)] \cdot \Pb[\Rc^{\epsilon}]\leq C \Pb[\Hc_j(r,n)] \times O(\epsilon^{c})=O(\epsilon^c) \Pb[\Hc_j(r,n)]\,.
$$
Therefore,
	\begin{equation}
		\label{lemma4.11_1}
 \Pb[\Hc_j^{\epsilon}(r,n) |\Hc_j(r,n)] \geq 1-O(\epsilon^c)\,.
	\end{equation}
For all $n \leq t \leq (1+\epsilon/K)n$, $\Pb[\Hc_j(r, t)|\Hc_j^{\epsilon}(r,n)] \geq 1-j(1-c)^K$. This is because on the event $\Hc_j^{\epsilon}(r,n)$ each outer face has a length of at least $\epsilon n$ which can be partitioned into at least $K$ pieces with a length of at least $(1+\epsilon/K)n$, and by applications of RSW theory, each piece has a probability at least $c$ to be connected with $C^+_{(1+\epsilon/K)n}$,  independently of others. Therefore, for all $n \leq t \leq (1+\epsilon/K)n$,
\begin{align*}
			h_j(t)&=\Pb[\Hc_j(r,t)] \geq \Pb[\Hc_j(r,t)|\Hc_j^{\epsilon}(r,n)] \times \Pb[\Hc_j^{\epsilon}(r,n) |\Hc_j(r,n)] \times \Pb[\Hc_j(r,n)]\\
			&\geq (1-j(1-c)^K) (1-O(\epsilon^c))h_j(r,n)\,.
\end{align*}
	Picking $K$ large and $\epsilon$ small, we then complete the proof. 
\end{proof}

\section{Proof of Lemma \ref{lem:sequence}}\label{sec:B}
\begin{proof}[Proof of Lemma \ref{lem:sequence}]
For $m \in (1.1, 10)$ and integer $k \geq 1$, let $b_k(m)=a_{m^{k+1}}/{a_{m^k}}$. We note that all constants in $O(\cdot)$'s in this proof are uniform w.r.t.\ $m \in (1.1,10)$. By Assumption (1) in the statement, $b_{k+1}(m)=b_k(m)(1+O(m^{-ck}))$. Then, it is easy to see that $\lim_{k \rightarrow \infty} b_k(m)$ exists (which will be denoted as $C(m)$) and is positive (since $b_k(m)$ is positive when $k$ is sufficiently large). In addition,
\begin{equation*}
b_k(m)=C(m)\Big(1+O(m^{-c k})\Big) \,, \quad \mbox{and equivalently} \quad \frac{a_{m^{k+1}}}{C(m)^{k+1}}=\frac{a_{m^{k}}}{C(m)^{k}}\Big(1+O(m^{-ck})\Big) \,.
\end{equation*}
Therefore, $\lim_{k \rightarrow \infty} \frac{a_{m^k}}{C(m)^k}$ exists (which will be denoted by $B(m)$) and is positive. In addition,
\begin{equation}
\label{lemma4.8_1}
a_{m^k}=B(m)C(m)^k\Big(1+O(m^{-c k})\Big) \,.
\end{equation}
Next, we use Assumption (2) to prove that there exist $\alpha \in (-\infty, \infty)$ and $C \in (0,\infty)$ such that $C(m)=m^{\alpha}$ and $B(m) = C$ for all $m$. It suffices to show that 
$$\mbox{$\log(C(m_1))/\log(m_1)=\log(C(m_2))/\log(m_2)$ and $B(m_1)=B(m_2)$ for all $m_1,m_2 \in (1.1,10)$.}$$ Fix $m_1,m_2 \in (1.1,10)$. WLOG we can assume that $\log_{m_2}(m_1)$ is irrational, because we can always find another $m_3$ such that $\log_{m_3}(m_1)$ and $\log_{m_3}(m_2)$ are both irrational and the case of $m_1,m_2$ follows from those of $m_1,m_3$ and $m_2,m_3$. Let $\delta,\epsilon>0$ denote two small constants to be chosen later. Since $\log_{m_2}(m_1)$ is irrational, we can find a sequence of increasing integers $\{p_j\}_{j \geq 1}$ and $\{q_j=\lfloor p_j\log_{m_2}(m_1) \rfloor \}_{j \geq 1}$ such that for all $j \geq 1$ \begin{equation}\label{eq:secc-1}q_j \leq p_j \log_{m_2}(m_1) \leq q_j+ \epsilon\,.\end{equation} By Assumption (2), there exists a constant $c=c(\delta)>0$ such that for all $\epsilon<c$
\begin{equation*}
\liminf_{n \rightarrow \infty} \inf_{10^{-\epsilon}n \leq s \leq t \leq n} \frac{a_t}{a_s}>1-\delta \,.
\end{equation*}
By \eqref{eq:secc-1} and the fact that $m_2<10$, we have $10^{-\epsilon} m_1^{p_j} \leq m_2^{-\epsilon} m_1^{p_j} \leq m_2^{q_j} \leq m_1^{p_j}$ and so for $j$ large enough $$(1-\delta) a_{m_1^{p_j}} \leq a_{m_2^{q_j}} \leq \frac{1}{1-\delta} a_{m_1^{p_j}} \,.$$ Together with \eqref{lemma4.8_1}
\begin{equation*}
(1-\delta) \Big(1+O(m_1^{-cp_j})+O(m_2^{-cq_j})\Big) \leq   \frac{B(m_2)C(m_2)^{q_j }}{B(m_1)C(m_1)^{p_j}}\leq \frac{1}{1-\delta }\Big(1+O(m_1^{-cp_j})+O(m_2^{-cq_j})\Big)\,.
\end{equation*}
Let $j$ tend to infinity. Since $\lim_{j \to \infty} \frac{q_j}{p_j} = \log_{m_2}(m_1)$, we have 
$$
\alpha:=\log\big(C(m_1)\big)/\log(m_1)=\log\big(C(m_2)\big)/\log(m_2).
$$
Note that when $\alpha \geq 0$, $10^{-\alpha \epsilon} m_1^{\alpha p_j} \leq m_2^{\alpha q_j } \leq m_1^{\alpha p_j}$; when $\alpha < 0$,  $10^{-\alpha \epsilon} m_1^{\alpha p_j} \geq m_2^{\alpha q_j } \geq m_1^{\alpha p_j}$. Therefore,
\begin{equation*}
(1-\delta)10^{-\epsilon(\alpha \wedge 0)} B(m_1)\Big(1+O(m_1^{-cp_j})\Big) \leq   B(m_2) \leq \frac{1}{1-\delta} B(m_1) 10^{\epsilon(\alpha \vee 0)} \Big(1+O(m_1^{-cp_j})\Big).
\end{equation*}
Let $j$ tend to infinity, then $\epsilon$ to zero and finally $\delta$ to zero. We have $B(m_1)=B(m_2)$. This completes the proof.
\end{proof}

\section{Proof of the strong separation lemma}\label{sec:A}
\begin{proof}[Proof of Proposition~\ref{prop:sepa}]
	The proof we employ here is a combination of techniques used in \cite{MR2630053} and \cite{MR3073882}.
	We only prove the first bullet point \eqref{eq:ex-sep}.
	Write $s=\dist(r, x_2,\cdots,x_{j}, -r)$. Then $Q(\Theta)=s/r$. Let $M=\lfloor\log_2\left( r/s \right)\rfloor\vee 0$ and $L=\lfloor\log_2\left(R/r\right)\rfloor$. Set
	\begin{equation*}
	r_i=\begin{cases}
	r, &i=0;\\
	r+2^{i-1}s, & 1\le i\le M;\\
	2^{i-M}r, & M+1\le i\le M+L.
	\end{cases}
	\end{equation*}
	Let $\Gamma_i$ be the collection of interfaces in $\Gamma$ truncated at their first hitting on $C^+_{r_i}$. Recall that $Q_{\rm ex}(\Gamma_i)$ is the exterior quality of $\Gamma_i$ on $C^+_{r_i}$ defined in \eqref{eq:ex-Q}. Let $d_i:=r_iQ_{\rm ex}(\Gamma_i)$ be the minimal distance between $-r_i,\,r_i$ and the endpoints of $\Gamma_i$. Define the \textbf{relative qualities}
	\begin{equation*}
	Q^*(i):=\begin{cases}
	d_i/(2^{i}s), & 0\le i\le M;\\
	Q_{\rm ex}(\Gamma_i), & M+1\le i\le M+L.
	\end{cases}
	\end{equation*}
	We set $Q^*(i):=0$ if not all $j$ interfaces manage to reach $C^+_{r_i}$.
	Furthermore, set 
	\[
	f_i:=\Pb[Q^*(i)>0],\quad g_i(\rho):=\Pb[Q^*(i)>\rho] \text{ for } \rho>0.
	\]
	With the above definitions, it suffices to show 
	\begin{equation}\label{eq:g-f}
	g_{M+L}(j^{-1})\ge c(j)f_{M+L}
	\end{equation}
	for some $c(j)>0$ that only depends on $j$. To this end, we need 
	the following facts about the quality, which can be obtained by using RSW-FKG gluing techniques, similar to the appendices of \cite{MR2630053,MR3073882} for the plane case. We omit the details.
	\begin{itemize}
		\item By application of the RSW theory, there is $c_{10}(j)>0$ such that
		\begin{equation}\label{eq:Q-1}
		f_1\ge c_{10}.
		\end{equation}
		\item For any $\delta>0$, there exists $\rho_0(\delta,j)$ such that for all $\rho\le \rho_0$ and $i\ge 0$,
		\begin{equation}\label{eq:Q-2}
		f_{i+1}-g_{i+1}(\rho)\le \delta f_i.
		\end{equation}
		\item For any $\rho>0$, there exists $R(\rho,j)>0$ such that for all $i\ge 0$,
		\begin{equation}\label{eq:Q-3}
		g_{i+1}(j^{-1})\ge R g_i(\rho).
		\end{equation}
	\end{itemize}
	
	Let $K=K(\rho)$ be the smallest integer in the range $1\le K\le M+L$ such that $g_i(\rho)\le f_i/2$ for all $K< i\le M+L$, where we set $K=M+L$ if $g_{M+L}(\rho)> f_{M+L}/2$. We claim that there exists $\rho_1>0$ such that for all $\rho\le\rho_1$, 
	\begin{equation}\label{eq:f_K}
	f_K\le 2g_K(\rho).
	\end{equation}
	This follows by definition if $K\ge 2$. 
	If $K=1$, we set $\delta_0:=c_{10}/2$ and $\rho_1:=\rho_0(\delta_0,j)$ with $c_{10}$ and $\rho_0$ from the first and second bullet points above respectively. By \eqref{eq:Q-2}, for all $\rho\le\rho_1$, we have $f_1-g_1(\rho)\le\delta_0$. This combined with \eqref{eq:Q-1} gives $f_1-g_1(\rho)\le f_1/2$ for all $\rho\le\rho_1$, which implies \eqref{eq:f_K}.
	
	In the following, for any $\delta>0$, we let $\rho\le \rho_0(\delta,j)\wedge \rho_1$, then both \eqref{eq:Q-2} and \eqref{eq:f_K} hold. Furthermore, by \eqref{eq:Q-2}, $f_{i}\le 2\delta f_{i-1}$ for all $K< i\le M+L$. Iterating this, we get 
	\begin{equation}\label{eq:f_M+L}
	f_{M+L}\le (2\delta)^{M+L-K}f_K\le (2\delta)^{M+L-K} 2g_K(\rho),
	\end{equation}
	where in the last inequality we have used \eqref{eq:f_K}.
	Repeated application of \eqref{eq:Q-3} gives that 
	\begin{equation}\label{eq:g_K}
	g_K(\rho)\le R(\rho)^{-1}R(j^{-1})^{K+1-M-L}g_{M+L}(j^{-1}).
	\end{equation}
	Combining \eqref{eq:f_M+L} with \eqref{eq:g_K}, we have 
	\[
	f_{M+L}\le 2R(\rho)^{-1}R(j^{-1})(2\delta/R(j^{-1}))^{M+L-K}g_{M+L}(j^{-1}).
	\]
	Letting $\delta$ be sufficiently small such that $2\delta/R(j^{-1})\le 1$, we obtain \eqref{eq:g-f} by setting $c(j)^{-1}=2R(\rho)^{-1}R(j^{-1})$.
	This completes the proof.\end{proof}

	\bibliographystyle{abbrv}
         \bibliography{arm}

\end{document}